\documentclass[reqno, 10pt]{amsart}
\usepackage[T1]{fontenc}
\usepackage[latin1]{inputenc}
\usepackage[english]{babel}
\usepackage{hyperref}

\usepackage{cite}
\usepackage{amssymb}
\usepackage{amsmath}
\usepackage{amsthm}
\usepackage{latexsym}
\usepackage{graphicx}
\usepackage{mathrsfs}
\usepackage{mathrsfs}
\usepackage{bbm}
\usepackage{verbatim}

\usepackage[usenames,dvipsnames]{color}

\newcommand{\EEE}{\color{black}}

\definecolor{rosso}{rgb}{0.8,0,0}

\def\luca #1{{\color{blue}#1}}

\def\luca #1{#1}

\newtheorem{thm}{Theorem}[section]

\newtheorem{lem}[thm]{Lemma}

\newtheorem{defin}[thm]{Definition}
\newtheorem{remark}[thm]{Remark}

\numberwithin{equation}{section}

\def\enne{\mathbb{N}}

\def\erre{\mathbb{R}}

\def\eps{\varepsilon}

\def\genspazio #1#2#3#4#5{#1^{#2}(#5,#4;#3)}
\def\spazio #1#2#3{\genspazio {#1}{#2}{#3}T0}

\def\L {\spazio L}
\def\H {\spazio H}

\def\Lx #1{L^{#1}(\Omega)}
\def\Hx #1{H^{#1}(\Omega)}

\def\R{\mathcal{R}}
\def\Wn{{\mathcal{W}_n}}
\def\ph{\varphi}
\def\chia{\eta}
\def\<#1>{\mathopen\langle #1\mathclose\rangle}
\def\norma #1{\mathopen \| #1\mathclose \|}

\def\intQt{\int_{Q_t}}

\def\iO{\int_\Omega}

\def\ept{{\eps\tau}}

\def\Pn {\Pi_n}
\newcommand{\ov}{\overline}
\def\lhs{left-hand side}
\def\rhs{right-hand side}

\def\Accorpa #1#2 #3 {\gdef #1{\eqref{#2}--\eqref{#3}}%
  \wlog{}\wlog{\string #1 -> #2 - #3}\wlog{}}

\def\beq{\begin{equation}}
\def\eeq{\end{equation}}

\def\to{\rightarrow}
\def\wto{\rightharpoonup}
\def\wstarto{\stackrel{*}{\rightharpoonup}}
\def\embed{\hookrightarrow}

\def\norm #1{\left\|#1\right\|}

\def\sp #1#2{\left<#1,#2\right>}
\newcommand\ip\sp

\renewcommand{\d}{{\mathrm d}}

\begin{document}
\title[Non-local phase-field models for tumor growth]
{On a class of non-local phase-field models\\for tumor growth with possibly
singular potentials, chemotaxis, and active transport}

\author{Luca Scarpa}
\address[Luca Scarpa]{Faculty of Mathematics, University of Vienna, 
Oskar-Morgenstern-Platz 1, 1090 Vienna, Austria.}
\email{luca.scarpa@univie.ac.at}
\urladdr{http://www.mat.univie.ac.at/$\sim$scarpa}

\author{Andrea Signori}
\address[Andrea Signori]{Dipartimento di Matematica e Applicazioni, 
Universit\`a di Milano-Bicocca, Via Cozzi 55, 20125 Milano, Italy.}
\email{andrea.signori02@universitadipavia.it}



\begin{abstract}
This paper provides a unified mathematical analysis of a family of
non-local diffuse interface models for tumor growth \luca{describing evolutions driven by 
long-range interactions.
These integro-partial differential equations model cell-to-cell 
adhesion by a non-local term and may be seen as}
non-local variants of the corresponding local model
proposed by H.~Garcke et al.~(2016).
The model in consideration couples 
a non-local Cahn--Hilliard equation for the tumor phase variable 
with a reaction-diffusion equation for the nutrient concentration, 
and takes into account also significant mechanisms 
such as chemotaxis and active transport.
The system depends on two relaxation parameters:
a viscosity coefficient and parabolic-regularization coefficient
on the chemical potential.
The first part of the paper is devoted to the analysis of the 
system with both regularizations. Here, a
rich spectrum of results is presented.
Weak well-posedness is first addressed, also 
including singular potentials. Then, 
under suitable conditions, existence of strong
solutions enjoying the separation property is proved.
This allows also to obtain a refined stability estimate with respect to the data,
including both chemotaxis and active transport.
The second part of the paper is devoted to the 
study of the asymptotic behavior of the system
as the relaxation parameters vanish.
The asymptotics are analyzed when the parameters
approach zero both separately and jointly, 
and exact error estimates are obtained.
As a by-product,
well-posedness of the corresponding limit systems is established. 
\end{abstract}

\maketitle

\section{Introduction}
\setcounter{equation}{0}
\label{sec:intro}

In the last decades, a vivid interest has been
devoted to the challenging project
of modeling tumor growth.
\luca{The main responsible of deaths due to cancer
is  often the formation of metastases in the late stages of the pathology,
when tumor cells spread also to separate parts of the host tissue 
and give rise to secondary tumor masses.
Several clinical studies have confirmed that the 
primary mechanism leading to this process is identified in 
the ability of cells to invade adjacent tissues \cite{sporn1996war}.
Invasion and metastasis have then deserved the 
unfortunate denomination of ``hallmarks of cancer'' in \cite{hanahan2000hallmarks}.
Mathematical modeling has then become a fundamental tool in order to 
describe and possibly predict these underlying processes:
validation, analysis, and simulation are crucial steps in the direction
of designation of anti-tumoral therapies.}
Many mathematical models have been proposed to capture the complexity of
the underlying biological and
chemical phenomena: in this direction we refer
to the seminal works \cite{AST, Byrne2, Byrne, BLM, CL, Fri} and references therein.

\luca{From the continuum diffuse-interface approach 
to tumor growth, the cellular adhesion is introduced by embodying surface tension force at the tumor surface and this procedure has been successfully employed in many
instances.}
In these models, the tumor evolution is described by 
introducing an order parameter $\ph$, taking values 
between $-1$ and $1$, and
representing the local concentration of tumoral cells.
The regions $\{\varphi=1\}$ and $\{\varphi=-1\}$
represent the pure tumorous and healthy phases, respectively,
whereas the diffuse interface $\{-1<\varphi<1\}$ models 
the narrow transition layer separating them.
One of the major advantages of this modeling approach is that,
unlike free boundary models, 
it takes into account also possible delicate behaviors
such as topological changes in the tumorous regions, 
occurring for example during break-up and coalescence phenomena.
The second main variable employed in the diffuse-interface description 
of tumor dynamics is the local concentration $\sigma$ of a certain nutrient
(e.g.~oxygen, glucose), in which the tissue in consideration is embedded.
The tumor is supposed to proliferate by absorption of such nutrient, 
and reversely the evolution of the nutrient is influenced by the consumption
by the tumor cells.
The key idea behind the diffuse-interface modeling
consists then of a non-trivial coupling of a 
phase-field-type equation for $\varphi$, usually
Cahn--Hilliard equation accounting for the phase segregation, 
with a reaction-diffusion equation for $\sigma$.
The proliferation and coupling terms appearing in the system
vary from model to model, and may take into account also
further biological mechanisms exhibited by the tumor such as
apoptosis, cell-to-cell adhesion, proliferation, chemotaxis, and active transport.

The classical local Cahn--Hilliard equation can be obtained as the conserved dynamics
in the $(H^1)^*$-metric
generated by the variational derivative of the Ginzburg--Landau free energy 
$\mathcal{E_{\rm loc}}$ 
with respect to the order parameter $\ph$, where 
\begin{align*}
	\mathcal{E}_{\rm loc}(\ph):= \int_\Omega \Big( \frac 12 |\nabla \ph|^2 + F (\ph) \Big)\,.
\end{align*}
Here, $F$ is a so-called double-well potential, possessing two global minima,
with typical choices being, in the order, the regular potential, the logarithmic potential 
and the double obstacle potential defined as:
\begin{align}
  \label{Fpol}
  &F_{pol}(r):=\frac14(r^2-1)^2\,, \quad r\in\erre\,,\\
  \label{Flog}
  &F_{log}(r):=\frac{\theta}{2}\left[(1+r)\ln(1+r)+(1-r)\ln(1-r)\right]-\frac{\theta_0}{2}r^2\,, \quad r\in(-1,1)\,,
  \quad0<\theta<\theta_0\,,\\
  \label{Fdobs}
  &F_{dob}(r):=\begin{cases}
  c(1-r^2) \quad&\text{if } r\in[-1,1]\,,\\
  +\infty \quad&\text{otherwise}\,,
  \end{cases}
  \qquad c>0\,.
\end{align}

\luca{In the context of tumor growth, 
the energy $\mathcal{E}_{\rm loc}$ accounts for
cell-to-cell adhesion.
More specifically, the term in $F$
models the fact that tumor cells 
prefer to adhere to each other rather than to non-tumor cells
(hence the pure phases tend to concentrate), 
while the gradient term penalizes too scattered tumor patterns
(hence high oscillations of $\varphi$).
Despite the fact that phase segregation described by means of 
the local Cahn--Hilliard equation is widely accepted in literature, 
the local model is not effective
in capturing cell-to-cell and cell-to-matrix adhesion phenomena 
driven by long-range competitions.
In the context of tumor growth models, 
neglecting long-range interaction is an enormous drawback.
Indeed, as we have pointed out above, 
the crucial biological process responsible of
the evolution of the cancer diseases are 
tumor-cell invasion and the formation of metastases.
The spread of secondary distant tumor masses is typically 
a long-range interaction process, and cannot be captured by 
means of the local modelling approach.
One of the possible ways to include long-range competitions 
and make the model more accurate in describing cell-invasion and 
metastases-formation is to switch to a {\em non-local} model instead.
This fact has been widely confirmed in the applied literature on
biological engineering and applied analysis, for which we refer to 
the numerous contributions
\cite{APS, engwer, Chap, GC, sherr, chen, chaplain2006mathematical}
and the references therein. In particular, the mentioned results
and their subsequent developments  
agree that cell-adhesion is typically a non-local-in-space phenomenon,
and represent then the crucial milestone for the validation and simulation
of non-local models in the context of tumor growth dynamics. 
}

\luca{
In the framework of diffuse-interface modeling of tumor growth, 
long-range interactions can be incorporated by modifiying the local energy 
$\mathcal{E}_{\rm loc}$ with a non-local one.
By following the ground-breaking work done by 
G.~Giacomin and J.~L.~Lebowitz on
non-local Cahn--Hilliard equations \cite{GL,GL1,GL2} 
(see also \cite{CFG,FGG2,FGGS,FGR2,gal-gior-grass}), 
we substitute the classical local Ginzburg--Landau free energy functional 
with the corresponding non-local Helmholtz free energy given by
\begin{align*}
 \mathcal{E}_{\rm nonloc}(\ph) := 
  \frac{1}{4} \int_{\Omega\times\Omega} J(x-y)|\ph(x)-\ph(y)|^2\, \d x\, \d y + \int_\Omega F(\ph)\,.
\end{align*}
Here $J$ stands for a sufficiently fast decaying kernel,
such as the classical Bessel or Newtonian potentials. 
Let us emphasize again that a non-local free energy 
as $\mathcal{E}_{\rm nonloc}$ is not motivated here 
by merely mathematical interests, but it has 
a fundamental meaning in the modeling of tumor growth:
it previously appeared in \cite{WLFC} and is indeed
crucial in describing the longe-range 
interaction processes involved in cell-invasion 
that would otherwise be left out in a local model.
For more details on the non-local Cahn--Hilliard equation 
we refer to the introduction of \cite{gal-gior-grass}, where a 
rich description concerning the state of the art on the equation is performed.
Besides, let us mention \cite{giorgini2020well, 
della2018nonlocal, conti2020well, giorgini2018cahn} 
for mathematical results related to variations of the 
classical Cahn--Hilliard equations with possibly singular potentials,
and to the recent works \cite{DST, DST2, DRST, MRT18} dealing with 
the asymptotic convergence of non-local Cahn--Hilliard equations 
to the respective local ones when the kernel suitably peaks around zero.}

The goal of this paper is to introduce and investigate a class of
non-local phase-field models for tumor growth inspired by the work by H.~Garcke et al.~\cite{garcke}.
Let $\eps,\tau \geq 0$, $\Omega\subset\erre^3$ be a smooth bounded domain, and
$T>0$ a fixed final time horizon. We consider 
a two-parameter class of non-local models in the following form:
\begin{align}
  \label{eq1}
  &\eps\partial_t\mu + \partial_t \varphi - \Delta\mu = (P\sigma - A)h(\varphi)
  \qquad&&\text{in } (0,T)\times\Omega\,,\\
  \label{eq2}
  &\mu=\tau\partial_t\varphi + a\varphi - J*\varphi + F'(\varphi) - \chi\sigma
  \qquad&&\text{in } (0,T)\times\Omega\,,\\
  \label{eq3}
  &\partial_t\sigma - \Delta\sigma + B(\sigma-\sigma_S) + C\sigma h(\varphi)
  =-\chia \Delta\varphi 
  \qquad&&\text{in } (0,T)\times\Omega\,,\\
  \label{eq4}
  &\partial_{\bf n}\mu =
  \partial_{\bf n}(\sigma-\eta\ph) = 0
  \qquad&&\text{on } (0,T)\times\partial\Omega\,,\\
  \label{eq5}
  &\eps\mu(0) = \eps\mu_0\,, \quad
  \varphi(0)=\varphi_0\,, \quad
  \sigma(0)=\sigma_0
  \qquad&&\text{in } \Omega\,.
\end{align}
\Accorpa \sys eq1 eq5
Let us briefly review the role of the occurring symbols.
The variable $\varphi$ represents the difference in volume fractions
between tumoral and healthy cells, with 
$\{\varphi=1\}$ being the pure tumoral phase, 
and $\{\varphi=-1\}$ being the pure healthy phase.
The variable $\mu$ is the chemical potential associated to $\ph$, 
and $\sigma$
represents the concentration of 
the unknown surrounding nutrient,
with the following convention: $\sigma \simeq 1 $ represents a rich nutrient concentration,
whereas $\sigma \simeq 0$ a poor one.
Furthermore, we indicate with $\bold n$ and $\partial_{\bold n}$ the normal vector and 
the corresponding directional derivative, 
$J$ is a spatial convolution kernel, with $a:=J*1$,
while $F'$ represents the derivative of a double-well potential $F$.
Precise assumptions are given in Section~\ref{sec:main} below.

The parameter $\tau\geq0$ represents the viscosity 
coefficient associated to the Cahn--Hilliard equation, 
while $\eps\geq0$ is a relaxation coefficient 
providing a parabolic regularization on the chemical potential.
The constants $P$, $A$, $B$, and $C$ are fixed positive real numbers,
taking into account
the proliferation rate of tumoral cells
by consumption of nutrient, the apoptosis rate,
the consumption rate of the nutrient with respect 
to a pre-existing concentration $\sigma_S$, 
and the nutrient consumption rate, respectively. 
Moreover, $\chi$ and $\chia$ are fixed non-negative constants, modeling the 
chemotaxis and active transport effects, respectively. 
For further insights concerning the 
modeling aspects, 
let us refer to \cite{garcke} (see also \cite{garcke-lam, garcke-lam2}), where the authors, 
after deriving some models from thermodynamic principles, underline how it is possible 
to decouple the chemotaxis by the active transport mechanism.
It is worth mentioning that, at least formally, 
by setting $\eps=\tau=0$ and by substituting the non-locality 
$a \ph - J * \ph$ with the ``corresponding local term''
$-\Delta \ph$, we obtain exactly a particular case of the setting analyzed in \cite{garcke}.
\luca{Here, we highlight that J. A. Sherratt et al. point out in \cite{sherr} that cell adhesion is 
intrinsically a non-local in space phenomenon, whereas the chemotaxis mechanism 
is on the other hand of local nature (see \cite{paint}). 
This motivates the medical and modeling relevance of system \sys\ in 
which a non-local term is considered in equation \eqref{eq2}, capturing 
long-term interaction processes occurring in cell-invasion,
but still keeping the terms related to chemotaxis of local nature in \eqref{eq2}--\eqref{eq3}.
For a situation in which non-local chemotaxis is addressed, we refer to \cite{buttenschon2018space}.}

\luca{Let us comment further on the structure of the parameters in system \eqref{eq1}--\eqref{eq5}.
As far as the coefficients are concerned, one can identify two main classes.
The parameters $P,A,B,C, \chi, \eta$ are structural coefficients of the model itself:
they arise directly from the practical description of the tumor dynamics, 
and each one is linked to an exact undergoing biological process.
For example, $P$ and $A$ take into account proliferation and death of tumor cells, 
$B$ and $C$ calibrate diffusion of nutrient with respect to a pre-existing 
concentration $\sigma_S$, and, more importantly, 
$\chi$ and $\eta$ render the tendency of tumor cells to attract nutrient 
and to move towards regions with high levels of nutrient, respectively.
The second group of parameters, namely $\eps$ and $\tau$, 
are connected on the other hand to specific mathematical regularizations of the model,
and have to be considered as small perturbations acting on the original limit system 
(i.e., \sys\ with $\eps=\tau=0$). In this perspective, 
it is of upmost importance to stress that the introduction 
of $\eps$ and $\tau$ {\em is not} aimed at a mere mathematical technical virtuosity,
but is finalized instead to the inclusion in the model of 
specific biological/physical mechanisms that are relevant to the 
tumor growth description. 
For instance,
the coefficient $\eps>0$ is necessary 
in order to deal with possible singular potentials $F$ due to the presence
of a mass source in the Cahn--Hilliard equation. As such,
bearing in mind that singular potentials are actually more
relevant in phase--segregation, 
the choice of analyzing the case $\eps>0$ 
has then to be interpreted as an additional 
possibility to include thermodynamically--relevant potentials
in the analysis, and not as a blunt mathematical exercise.
In the same spirit, the introduction of the parameter $\tau$ in the model 
is not end in itself, but is aimed at keeping the relevant
cross--diffusion mechanisms of chemotaxis and active transport,
which otherwise could not be covered by the model.
This being stressed, the regularized system \sys\ that we propose 
is purposely very general, in order to provide a larger variety 
of frameworks that could be covered by the model
and that could adapt to different practical scenarios.
}

Up to the author's knowledge, there are still few contributions 
devoted to the mathematical analysis of
non-local tumor growth models: we recognize \cite{Lima3, FLOW, frig-lam-roc}.
By contrast, the local situation has been the subject of intensive studies.
At first, let us point out some models which neglect velocity contribution which 
are somehow variations of the model introduced by A.~Hawkins-Daarud et.~al in \cite{hawk2} 
(see also \cite{hawk, hil-zee}).
In this direction, we mention \cite{frig-grass-roc}, where the well-posedness of the system
is shown under general polynomial growth type assumptions for the involved potentials.
In \cite{col-gil-hil} (along with the related works \cite{col-gil-roc-spr, col-gil-roc-spr2}),
the authors consider some regularized version compared to \cite{frig-grass-roc},
by adding the same regularization that we have introduced here
on the viscosity and the dissipation of the chemical potential. 
Owing to these terms the authors were able to 
extend the setting of some analytic results including in the investigation also singular 
and possibly non-regular potentials like the logarithmic \eqref{Flog}
or the double obstacle one \eqref{Fdobs}. 
Moreover, the authors 
established in which sense it is possible to let these regularizations parameters to zero,
recovering some of the results already proved in \cite{frig-grass-roc}: in this sense,
our work is somehow inspired by these contributions. 
Let us also refer to \cite{col-gil-spr-frac1, col-gil-spr-frac2},
where a similar investigation was performed for fractional models.
Furthermore, in order to better emulate in-vivo tumor-growth, 
other authors have proposed to include 
fluid motion by further coupling previous systems with a velocity law of Darcy's 
or Brinkmann's-type; we refer in particular to
\cite{garcke, garcke-lam, garcke-lam2, garcke-lam3, 
garcke-lam4, eben-garcke1, eben-garcke2, eben-lam, FLOW, WLFC, JWZ, agosti2017cahn}.
We point out the recent work \cite{garcke-lam-signori} (see also \cite{Lima1,Lima2})
wrote by the second-named author in collaboration with H.~Garcke and K.~F.~Lam, 
where elasticity effects are taken into account
as physical evidence have shown that the presence of the extracellular matrix
or rigid bone can assert significant influences on tumor proliferation. 
For multi-species tumor growth models, 
we point out \cite{Dai, frig-lam-rocca-schim, garcke-lam-nur}.

Moreover, a wide number of results concerning further analysis on these models 
have been performed. In this direction, we mention the
optimal control problems addressed by \cite{SW, col-gil-roc-spr-control, 
eben-knopf1, eben-knopf2, garcke-lam-roc, col-gil-spr-fraccontrol, col-gil-mar, kahle-lam-latz}.
In particular, we mention \cite{col-signori-spre} which deals with the optimal control problem
for the corresponding local version of system \sys, and we also point out 
\cite{signori1, signori2, signori3, signori4, signori5, col-gil-spr-fraccontrol}, where 
similar relaxed models have been investigated from the optimal control viewpoint.
Let us also point out \cite{CRW, MRS, col-gil-hil}, where some long-time behavior 
for similar models is addressed.
To conclude the overview, let us mention the work \cite{orr-roc-scar},
where a phase-field model for tumor growth 
has been analyzed also taking into account possible stochastic perturbations
of the system. The paper, written by the first-named author in collaboration 
with C.~Orrieri and E.~Rocca, focuses on well-posedness and optimal 
control of treatment when two Wiener-type noises act on the 
proliferation of tumor cells and evolution of nutrient.


Let us present now the main results of the present paper.

The first part of the work is devoted to the analysis of the system \eqref{eq1}--\eqref{eq5}
when both regularizations are present, i.e.~with $\eps,\tau>0$.
In this setting, we first investigate existence of weak solutions, even 
when singular potentials as \eqref{Flog} or \eqref{Fdobs} are present,
also including chemotaxis and active transport.
Secondly, we show that without active transport 
(i.e.~$\chia=0$) continuous dependence on the data (hence uniqueness) 
holds for weak solutions. Furthermore, we investigate regularity properties of the solutions, 
and prove existence of strong solutions as well as separation results from the 
potential barriers. For strong solutions, we are finally able to refine 
the stability estimate with respect to the data, also including the case
of chemotaxis and active transport.

The second part of the work is focused on the study of the asymptotic behavior of the system
as $\eps \searrow 0$ and/or $\tau \searrow 0$.
These are performed both separately (i.e.~$\eps\searrow0$ with $\tau>0$, and
$\tau\searrow0$ with $\eps>0$) and jointly (i.e.~$\eps,\tau\searrow0$).
In each of these cases, under suitable conditions we are able to show convergence
of the system to the respective limit problem, hence also the corresponding well-posedness.
Also, we give the exact rates of convergence through precise error estimates.

Let us briefly mention here the mathematical challenges that 
we have to overcome in these asymptotics.

\noindent
{\bf Passage to the limit as $\boldsymbol{\eps \searrow 0}$.}
In this first asymptotic study the parabolic regularization on $\mu$ is ``removed'',
resulting in lack of regularity on the chemical potential.
As a consequence, due to the presence of proliferation terms in the Cahn--Hilliard equation, 
a very natural growth condition on the potential has to be required (c.f.~\eqref{pol_growth}),
allowing for any polynomial or first-order exponential potentials.
The passage to the limit, hence the existence for the limit problem with $\eps=0$,
is proved in the setting of no active transport term (i.e.~$\chia  = 0$),
due to the need of a maximum principle argument for $\sigma$.
As for the error estimate (and therefore the uniqueness for the limit system),
a rate of convergence of order $\eps^{1/4}$ is 
obtained by showing refined estimates on the solutions
and exploiting a locally-Lipschitz assumption on the potential 
(still including the classical case \eqref{Fpol} for example).

\noindent
{\bf  Passage to the limit as $\boldsymbol{\tau \searrow 0}$}
In the second passage to the limit, the viscosity of the Cahn--Hilliard equation
vanishes, and this results is a loss of regularity on the phase-variable.
The presence of $\eps>0$ still allows passing to the limit in very general
settings, such as singular potentials, chemotaxis, and active transport,
only requiring some compatibility conditions (smallness-type assumptions) on the constants.
The separation from the potential barriers is not preserved though, as it is naturally expectable.
Moreover, a corresponding error estimate showing a convergence rate of order
$\tau^{1/2}$ is obtained (and therefore the uniqueness for the limit system).

\noindent
{\bf Passage to the limit as $\boldsymbol{\eps,\tau \searrow 0}$.}
In the last passage to the limit, the parameters $\eps$ and $\tau$ vanish 
simultaneously. Here, the convergence is proved by proving some 
refined estimates on the solutions, depending on both parameters, and 
combining the assumptions above on the potential and the coefficients.
Moreover, the error estimate (and the resulting well-posedness of the limit
problem) is obtained with a rate of convergence of $\eps^{1/4}+\tau^{1/2}$,
under a suitable scaling on the two parameters.

The plan of the rest of the paper is as follows. In
Section \ref{sec:main} we set our notation and present the obtained results.
The weak and strong well-posedness of \sys\ for $\eps,\tau >0$ is addressed in Section \ref{sec:WP}.
Then, Sections \ref{sec:eps}, \ref{sec:tau} and \ref{sec:eps_tau}
are completely devoted to the asymptotic analysis of the system as $\eps$ and $\tau$
approach zero, first separately and then jointly.

\section{Setting, assumptions, and main results}
\label{sec:main}

Throughout the paper, $\Omega\subset\erre^3$ is a smooth bounded domain
and $T>0$ is a fixed final time. We set for convenience the spatiotemporal cylinders 
\[
  Q:=(0,T)\times\Omega\,, \qquad\Sigma:=(0,T)\times\partial\Omega\,, 
  \qquad Q_t:=(0,t)\times\Omega\,, \quad t\in(0,T)\,,
\]
and we introduce the functional spaces
\[
  H:=L^2(\Omega)\,, \qquad V:=H^1(\Omega)\,, \qquad
  W:=\left\{y\in H^2(\Omega):\partial_{\bf n}y=0 \text{ a.e.~on } \partial\Omega\right\}
\]
endowed with their natural norms $\norm{\cdot}_H$, 
$\norm{\cdot}_V$, and $\norm{\cdot}_W$, 
respectively. 
Likewise, we use $\norma{\cdot}_{p}$ to indicate the standard
norm of the space $\Lx p$, for all $p \in [1,\infty]$.
As usual, $H$ is identified with its dual $H^*$
through its Riesz isomorphism, so that 
\[
  W\embed V\embed H \simeq H^* \embed V^*\embed W^*\,,
\]
where all inclusions are dense, continuous, and compact. The duality 
pairing between $V^*$ and $V$, and the scalar product in $H$
will be denoted by the symbols $\ip{\cdot}{\cdot}$ and $(\cdot,\cdot)$,
respectively.

For every $f\in L^1(0,T)$ we will use the notation 
$(1\star f)(t):=\int_0^t f(s)\,\d s$, for $t \in[0,T]$.

Moreover, for every $v\in V^*$ we set 
$v_\Omega:=\frac1{|\Omega|}\ip{v}{1}$ for the 
generalised mean value of $v$. 
Let us also recall the Poincar\'e-Wirtinger inequality 
\begin{align}
  \norma{v}^2_V \leq C_{\Omega} \bigl( \norma{\nabla v}_H^2 + |v_\Omega|^2 \bigr)\,,
  \quad \forall\, v\in V\,,
  \label{poincare}
\end{align}
where the constant $C_\Omega>0$ depends only on $\Omega$.
Let us recall that the Laplace operator with homogeneous Neumann conditions may be seen
as a variational operator
\[
  -\Delta:V\to V^*\,, \qquad \ip{-\Delta v}{\zeta}:=\int_\Omega\nabla v\cdot\nabla \zeta\,, \quad\forall\,v,\zeta\in V\,.
\]
It is well know, as a consequence of the Poincar\'e-Wirtinger inequality \eqref{poincare}, 
that the restriction
of $-\Delta$ to the subspace of null-mean elements of $V$ is injective, and that it 
possesses a well defined inverse 
\[
  \mathcal N:\{v^*\in V^*:v^*_\Omega=0\}\to\{v\in V: v_\Omega=0\}\,.
\]
Lastly, let $\R=I-\Delta: V \to V^*$ be the Riesz isomorphism of $V$, i.e. the map
\begin{align*}
	\<\R u, v> := \int_\Omega (\nabla u \cdot \nabla v + u v)\,,
	\quad\forall\ u,v\in V\,.
\end{align*}
It is well-known that $\R_{|W}$ yields an isomorphism from $W$ to
$H$ with well-defined inverse $\R^{-1} : H \to W$.
In addition, for all $v \in V,$ and $v^*,w^* \in V^*$,
the following properties hold
\begin{align*}
\<{\R v},{\R^{-1}v^*}> &= \<{v^*},{v}>, \quad 
\<{v^*},{\R^{-1}w^*}> = (v^*, w^* )_{*}\,,
\end{align*}
where the symbol $(\cdot, \cdot )_{*}$ denotes the
inner product of $V^*$.
Furthermore, for every $f \in V$ it holds that
\begin{align*}
	\norma{f}_H^2 = \<{f},{f}> = \<{\R f},{\R^{-1}f}> 
	\leq \norma{f}_{V} \norma{\R^{-1} f}_{V} 
	\leq \norma{f}_{V} \norma{f}_{V^*}\,.
\end{align*}
Besides, if $v^* \in \H1 {V^*}$, we have for a.e~$t \in (0,T)$ that
\begin{align*}
\<{\partial_t v^*(t)},{\R^{-1}v^*(t)}> = \frac 12 \frac d {dt} \norma{v^*(t)}^2_{*}\,.
\end{align*}

The following structural assumptions on the data 
will be in order in the paper.
\begin{description}
  \item[A1] $P,A,B,C, \chi,\chia$ are non-negative constants. 	
  \item[A2] $h:\erre\to\mathopen[0,+\infty\mathclose)$ is bounded and 
  Lipschitz continuous.
  \item[A3] $\sigma_S\in L^\infty(Q)$ and
  \[
  0\leq \sigma_S(t,x) \leq 1 \quad\text{for a.e.~}(t,x)\in Q\,.
  \]
  \item[A4]
  \label{potential}
  $F:=F_1 + F_2 \geq0$, where 
  \[
  F_1:\erre\to[0,+\infty] \quad\text{is proper, convex, and lower semicontinuous}\,,
  \]
  and 
  \[
  F_2\in C^1 (\erre)\,, \qquad F_2':\erre\to\erre \quad\text{is Lipschitz continuous}\,,
  \qquad F_2'(0)=0\,.
  \]
  In particular, the subdifferential 
  $\partial F_1:\erre\to2^{\erre}$ is well defined in the sense 
  of convex analysis, and we assume that $0\in \partial F_1(0)$.
  The Moreau regularization of $F_1$ and the Yosida approximation of $\partial F_1$ 
  are defined, respectively, as
   \[
  F_{1,\lambda}:\erre\to\mathopen[0,+\infty\mathclose)\,,\qquad
  F_{1,\lambda}(r):=F_1(0) + \int_0^rF'_{1,\lambda}(s)\,\d s\,, \quad r\in\erre\,,
  \]
  and
  \[
  F'_{1,\lambda}:\erre\to\erre\,, \qquad
  F_{1,\lambda}':=\frac{I-(I+\lambda\partial F_1)^{-1}}{\lambda}\,, \quad\lambda>0\,,
  \]
  where $I$ stands for the identity operator.
  We recall that $F'_{1,\lambda}$ is 
  $\frac1\lambda$-Lipschitz continous and we set 
  \[
  F_\lambda:= F_{1,\lambda} + F_2\,.
  \]
  \item[A5] The kernel $J\in W^{1,1}_{loc}(\erre^3)$ is such that $J(x)=J(-x)$ for 
  a.e.~$x\in\erre^3$.
  For any measurable $v:\Omega\to\erre$ we use the notation
  \[
  (J*v)(x):=\int_\Omega J(x-y)v(y)\,\d y\,, \quad x\in\Omega\,,
  \]
  and set $a:=J*1$.
  Moreover, we assume that 
  \[
  a_*:=\inf_{x\in\Omega}\int_\Omega J(x-y)\,d y=\inf_{x\in\Omega}a(x)\geq0\,,
  \]
  \[
  a^*:= \sup_{x\in\Omega} \int_\Omega |J(x-y)|\,\d y < +\infty, \qquad
  b^*:= \sup_{x\in\Omega} \int_\Omega |\nabla J(x-y)|\, \d y < +\infty\,,
  \]
  and we set $c_a:=\max\{a^*-a_*,1\}$.
  Finally, we suppose that 
  there exists a positive constant $C_0$ such that 
  \[
  a_* + \frac{w_1-w_2}{r_1-r_2} \geq C_0 \,,
  \qquad\forall\,r_i\in D(\partial F_1)\,,\quad\forall\,w_i\in\partial F_1(r_i)+F_2'(r_i)\,,
  \quad i=1,2\,,\quad r_1\neq r_2\,.
  \]
  Note that if $F$ is of class $C^2$, the last condition is equivalent to the classical one
  \[
  a_* + F''(r)\geq C_0 \qquad\forall\,r\in D(F')\,,
  \]
  where $D(F')$ denotes the domain of $F'$.
\end{description}

\luca{For convenience, we introduce the following 
upper bounds for the coefficients $\eps$ and $\tau$
\[
  \eps_0:=\min\left\{\frac1{4c_a}, \frac{1}{\max\{1,  a^* - \min\{ a^*,C_0 \}\}},
  \frac{2C_0}{3(a^*+b^*)^2K_0^2}\right\}\,, \qquad \tau_0:=1\,,
\]
where $K_0$ denotes the norm of the continuous inclusion $H\embed V^*$.
This is only a technical requirement on the coefficients, which
is clearly not restrictive as $\eps$ and $\tau$ have to be considered as 
small perturbations.
}

The first main result deals with existence of global weak solutions to the system 
\eqref{eq1}--\eqref{eq4} under very general assumptions on the data.
In particular, any type of potential as in \eqref{Fpol}--\eqref{Fdobs}
is included in this first result.

\begin{thm}[Existence of weak solutions: $\eps,\tau>0$]
  \label{thm1}
  Assume {\bf A1}--{\bf A5}, and
  let \luca{$\eps\in(0,\eps_0)$ and $\tau\in(0,\tau_0)$}.
  Moreover, let the
  triple of initial data $(\varphi_0, \mu_0, \sigma_0)$ satisfy
  \beq\label{ip_init}
   \varphi_0\in V\,,\qquad F(\varphi_0)\in L^1(\Omega)\,,  
  \qquad \mu_0, \sigma_0\in H.  \eeq
  Then, there exists a quadruplet $(\varphi,\mu,\sigma,\xi)$ such that 
  \begin{align}
    \label{phi}
    \varphi&\in H^1(0,T; H)\cap L^\infty(0,T; V)\,,\\
    \label{mu}
  \mu, \sigma &\in H^1(0,T; V^*)\cap L^2(0,T; V)\,,\\
    \label{xi}
    \xi&\in L^2(0,T; H)\,,
 \end{align}
 where
 \beq
    \label{def_mu_xi}
   \mu=\tau\partial_t\varphi
    + a\varphi - J*\varphi + \xi + F_2'(\varphi) - \chi\sigma\,, \qquad
    \xi\in\partial F_1(\varphi) \quad\text{ a.e.~in } Q\,,
  \eeq
  with $\varphi(0)=\varphi_0$, $\mu(0)=\mu_0$, $\sigma(0)=\sigma_0$ in $H$, and such that 
  \begin{align}
  \label{var1}
    &\ip{ \partial_t (\eps \mu +  \varphi)}\zeta
    +\int_\Omega\nabla\mu\cdot\nabla\zeta 
    =\int_\Omega(P\sigma - A)h(\varphi)\zeta\,,\\
  \label{var2}
  &\ip{\partial_t\sigma}\zeta +\int_\Omega \nabla\sigma\cdot\nabla\zeta
  + \int_\Omega\left(B(\sigma-\sigma_S) + C\sigma h(\varphi) \right)\zeta
  =\chia\int_\Omega\nabla\varphi\cdot\nabla\zeta\,,
  \end{align}
  for every $\zeta\in V$, almost everywhere in $(0,T)$. Furthermore, if
  $\chia=0$ and
  \beq\label{ip_infty}
  0 \leq \sigma_0(x) \leq 1 \quad\text{for a.e.~}x\in\Omega\,,
  \eeq
  then $\sigma(t)\in L^\infty(\Omega)$ for all $t\in[0,T]$ and 
  \beq\label{max_sigma}
  0 \leq\sigma(t,x)\leq 1 \quad\text{for a.e.~}x\in\Omega\,, 
  \quad\forall\,t\in[0,T]\,.
  \eeq
\end{thm}
It is worth mentioning that, in the case of singular potentials such as \eqref{Flog}
and \eqref{Fdobs}, the assumption $F(\ph_0)\in L^1(\Omega)$ entails that 
$\ph_0\in L^\infty(\Omega)$ and that $|\ph_0(x)|\leq 1$ for 
almost every $x \in \Omega$.

The second result concerns continuous dependence of the data 
for weak solutions. This result applies again to any choice of the potential $F$,
but we are forced (so far) to restrict ourselves to the case without active transport (i.e. $\chia=0$).

\begin{thm}[Continuous dependence: $\eps,\tau>0$]
  \label{thm2}
  Assume {\bf A1}--{\bf A5}, and let $\chia=0$,
  \luca{$\eps\in(0,\eps_0)$ and $\tau\in(0,\tau_0)$}.
  Then there exists a constant $K>0$ \luca{independent of $\tau$} such that, for any pair
  of initial data $\{(\varphi_0^i, \mu_0^i, \sigma_0^i)\}_{i}$, $i=1,2$, 
  satisfying \eqref{ip_init} and \eqref{ip_infty}, 
  and
  for any respective solutions
  $\{(\varphi_i, \mu_i, \sigma_i, \xi_i)\}_i$, $i=1,2$, 
 satisfying \eqref{phi}--\eqref{max_sigma},
 it holds that
  \begin{align}
    & \notag\luca{\norm{(\eps\mu_1+\varphi_1)-(\eps\mu_2+\varphi_2)}_{L^\infty(0,T; V^*)}}
    + \norm{\mu_1-\mu_2}_{L^2(0,T; H)} \\
    &\notag
    \qquad
    + \tau^{1/2}\norm{\varphi_1-\varphi_2}_{C^0([0,T]; H)}
    + \norm{\varphi_1-\varphi_2}_{L^2(0,T; H)}
    +\norm{\sigma_1-\sigma_2}_{C^0([0,T]; H)\cap L^2(0,T; V)}
    \\ &
    \leq \luca{K\left(\norm{(\eps\mu_0^1+\varphi_0^1)-(\eps\mu_0^2+\varphi_0^2)}_{V^*}
    +\tau^{1/2} \norm{\varphi_0^1-\varphi_0^2}_{H}
    +\norm{\sigma_0^1-\sigma_0^2}_H
     \right)\,.}
     \label{cont_dep_first}
  \end{align}
\end{thm}

As a consequence of the above result, we infer the uniqueness
of the weak solution obtained in Theorem~\ref{thm1}
under the only additional requirement that $\chia =0$.
The next result deals with the regularity of weak solutions
with respect to the data.
\begin{thm}[Regularity: $\eps,\tau>0$]
  \label{thm3}
  Assume {\bf A1}--{\bf A5},
  \luca{$\eps\in(0,\eps_0)$, and $\tau\in(0,\tau_0)$}. Moreover, let the
  triple of initial data $(\varphi_0, \mu_0, \sigma_0)$ satisfy \eqref{ip_init} and also
  \beq\label{ip_init_reg}
  \exists\,\xi_0\in H:\;\xi_0\in\partial F_1(\varphi_0) \text{ a.e.~in } \Omega\,, 
  \qquad \mu_0,\sigma_0\in V\,,
  \eeq
  and suppose that $t=0$ is a Lebesgue point for $\sigma_S$ with
  \beq
    \label{ip_data_reg}
    \sigma_S(0) \in H\,.
  \eeq
  Then, the solution $(\varphi,\mu,\sigma,\xi)$ to \eqref{phi}--\eqref{var2}
  given by Theorem~\ref{thm1}
  satisfies
  \begin{align}
    \label{phi_reg}
    \varphi&\in W^{1,\infty}(0,T; H)\cap L^\infty(0,T; V)\,, \\
    \label{mu_sigma_reg}
    \mu, \sigma - \chia \ph&\in H^1(0,T; H)\cap L^\infty(0,T; V)\cap L^2(0,T; W)\,,\\
    \label{sigma_reg}
    \sigma&\in H^1(0,T; H)\cap L^\infty(0,T; V)\,.
  \end{align}
\end{thm}

Our next result is concerned with 
the separation property,
magnitude regularity, and existence of strong solutions. 
In this direction, we postulate the following assumptions for 
$F$ and $J$.

\begin{description}
  \item[A6] Setting $(-\ell,\ell):=\operatorname{Int}D(\partial F_1)$,
  with $\ell\in[0,+\infty]$, we assume that 
  \[
  F\in C^4(-\ell,\ell)\,, \qquad 
  \lim_{r\to(\pm \ell)^\mp}\left[F'(r) - \chi\chia r\right]=\pm \infty\,.
  \]
\end{description}
It is worth pointing out that {\bf A6} excludes
potentials $F$ of double-obstacle type as in \eqref{Fdobs}.
Nevertheless, the logarithmic potential \eqref{Flog} and 
any polynomial super-quadratic potential as \eqref{Fpol} is allowed.
\luca{Let us also remark that assuming the effective domain of $\partial F_1$ to be symmetric with respect to zero is mainly a matter of convenience, so to allow \eqref{Fpol}--\eqref{Flog} to be included.
In general, the symmetry condition for the domain of $\partial F_1$
it is not necessary from the analysis viewpoint, 
and one can always reconstruct this situation by renormalization of $F$.
}

As for the kernel, a natural requirement from the analytical point of view is to require 
\beq\label{Jw21}
  J\in W^{2,1}(\mathcal B_R), 
  \quad \hbox{where $\mathcal B_R:=\{x\in\erre^3:|x|< R:=\operatorname{diam}(\Omega)\}$}\,,\quad R >0.
\eeq
However, this condition prevents some relevant cases of kernels 
such as the Newtonian or the Bessel potential from being considered.
Following the ideas of \cite{FGG2,
gal-gior-grass} (see also \cite[Def.~1]{bed-rod-bert}), it is possible to 
cover also these situations by replacing the 
above condition by assuming that $J$ 
is {\it admissible} in the following sense.
\begin{defin}
\label{def_adm}
A convolution kernel $J\in W^{1,1}_{\rm loc}(\erre^3)$ is admissible if
it satisfies:
\begin{itemize}
\item $J\in C^3(\erre^3 \setminus \{0\})$.
\item $J$ is radially symmetric, i.e.~$J(\cdot)=\widetilde{J}(|\cdot|)$
	for a non-increasing $\widetilde J:\erre_+\to\erre$.
\item there exists $R_0>0$ such that
	$r\mapsto{\widetilde J}''(r)$ and $r\mapsto{\widetilde J}'(r)/r$ are monotone on $(0,R_0)$.
\item there exists $C_d>0$ such that 
	$|D^3 J(x)|\leq C_d|x|^{-4}$ for every $x\in\erre^3\setminus\{0\}$.
\end{itemize}
\end{defin}
Thus, we require
\begin{description}
  \item[A7] $J$ satisfies \eqref{Jw21}
  or it is admissible in the sense of Definition~\ref{def_adm}.
\end{description}

\begin{thm}[Existence of strong solutions, separation property: $\eps,\tau>0$]
  \label{thm4}
  Assume conditions {\bf A1}--{\bf A7}, and
  let \luca{$\eps\in(0,\eps_0)$ and $\tau\in(0,\tau_0)$}. Let the
  initial data $(\varphi_0, \mu_0, \sigma_0)$ 
  satisfy \eqref{ip_init}, \eqref{ip_init_reg}, and also
  \beq\label{ip_init_sep}
  \varphi_0\in H^2(\Omega)\,,\qquad
  \mu_0,\sigma_0\in L^\infty(\Omega)\,,\qquad
  \exists\,r_0\in(0,\ell):\;
  \norm{\varphi_0}_{L^\infty(\Omega)}\leq r_0\,.
  \eeq
  Then, the solution $(\varphi,\mu,\sigma,\xi)$ to \eqref{phi}--\eqref{var2}
  given by Theorems~\ref{thm1} and \ref{thm3} satisfies
  \begin{align}
    \label{phi_reg2}
    &\varphi\in W^{1,\infty}(0,T; V)\cap H^1(0,T; H^2(\Omega))\,, \qquad\partial_t\varphi\in L^\infty(Q)\,,
    \qquad\chia\varphi\in L^2(0,T; W)\,,\\
    \label{phi_sep}
    &\exists\,r^*\in(r_0,\ell):\quad\sup_{t\in[0,T]}\norm{\varphi(t)}_{L^\infty(\Omega)}\leq r^*\,,\\
    \label{mu_reg2}
    &\mu, \sigma \in H^1(0,T; H)\cap L^\infty(0,T; V)\cap L^2(0,T; W)\cap L^\infty(Q)\,.
  \end{align}
  In particular, equations \eqref{eq1}--\eqref{eq3} hold almost everywhere in $Q$.
\end{thm}

\begin{remark}
\noindent $(i)$ 
Note that the equation \eqref{eq2} at time $0$ reads 
\[
\mu_0=\tau\varphi_0' + a\varphi_0 - J*\varphi_0 + F'(\varphi_0) - \chi\sigma_0\,,
\]
where $\varphi_0'$ ``represents'' the initial 
value of the time-derivative of $\varphi$. 
Under the assumptions \eqref{ip_init}, \eqref{ip_init_reg}, and \eqref{ip_init_sep}
we have that $\varphi_0'\in V\cap L^\infty(\Omega)$, hence the 
improved regularities $\partial_t\varphi\in L^\infty(0,T; V)\cap L^2(0,T; H^2(\Omega))$
and $\partial_t\varphi\in L^\infty(Q)$ obtained in Theorem~\ref{thm4} are naturally expectable.

\luca{
\noindent $(ii)$ 
Let us point out that \eqref{ip_init_sep}, in the case \eqref{Flog}, prevents the initial tumor distribution $\ph_0$ to possess any region occupied by solely tumorous cells as $r_0 < \ell$.}
\end{remark}

Relying on the extra-regularity and the separation property,
we are able to show a refined continuous dependence result for strong solutions,
where the stability estimates are verified in stronger topologies.
Let us stress that in this case we are able to cover 
also the scenarios of chemotaxis and active transport,  
complementing thus the previous Theorem~\ref{thm2}. 

\begin{thm}[Refined continuous dependence: $\eps,\tau>0$]
   \label{thm5}
  Assume {\bf A1}--{\bf A7}, let \luca{$\eps\in(0,\eps_0)$, and $\tau\in(0,\tau_0)$}.
  Then for any pair
  of initial data $\{(\varphi_0^i, \mu_0^i, \sigma_0^i)\}_i$, $i=1,2$, 
  satisfying \eqref{ip_init}, \eqref{ip_init_reg}, and \eqref{ip_init_sep},
  there exists a constant $K>0$ such that,
  for any respective solutions
  $\{(\varphi_i, \mu_i, \sigma_i, \xi_i)\}_i$, $i=1,2$,
  satisfying \eqref{phi}--\eqref{var2} and \eqref{phi_reg2}--\eqref{mu_reg2},
  it holds that
  \begin{align*}
  &\norm{\mu_1-\mu_2}_{H^1(0,T; H)\cap L^\infty(0,T; V)\cap L^2(0,T; W)}+
  \norm{\varphi_1-\varphi_2}_{W^{1,\infty}(0,T; V)\cap H^1(0,T; H^2(\Omega))}\\
  &\qquad
  +\norm{\sigma_1-\sigma_2}_{H^1(0,T; H)\cap L^\infty(0,T; V)\cap L^2(0,T; W)}\\
  &\leq K\left(\norm{\mu_0^1-\mu_0^2}_V + \norm{\varphi_0^1-\varphi_0^2}_{H^2(\Omega)}
  +\norm{\sigma_0^1-\sigma_0^2}_V \right)\,,
  \end{align*}
  where $K$ only depends on $\Omega, T, 
  \eps,\tau,P,A,B,C,C_0,a_*, a^*, b^*, r^*, \norm{F}_{C^4([-r^*, r^*])}$
  and $\{(\varphi_0^i, \mu_0^i, \sigma_0^i)\}_{i=1,2}$.
\end{thm}
In particular, under the assumptions 
\eqref{ip_init}, \eqref{ip_init_reg}, and \eqref{ip_init_sep}
on the data, we deduce that the uniqueness of strong solutions in the sense
of Theorem~\ref{thm4} holds.

We now will present the results concerning the asymptotic analysis of \sys\
with respect to the parameters $\eps$ and $\tau$.
To begin with, we consider the case $\eps\searrow0$, assuming
$\tau>0$ to be fixed. In this direction, 
we need to enforce the conditions on the potential $F$. 
In fact,  proceeding with classical estimates, just a bound of 
$\nabla \mu$ in $\L2 H$ can be proved, having no information on the behavior 
of $\mu$ in $\L2 H$. 
This gap is usually bridged via the application of a Poincar\'e-type inequality,
which yields the control of $\mu$ in $\L2 V$. To this end,
some control on the spatial mean of $\mu$ is necessary:
if $\eps>0$ is fixed, this follows automatically from the estimates,
whereas in the limit $\eps\searrow0$ it has to be 
obtained from a suitable prescription on the potential.
Namely, the assumption
\begin{align}
	\label{pol_growth}
	D(\partial F_1)= \erre, \qquad 
	\exists\,C_F>0:\quad 
	|\partial F_1^0(r)| \leq C_F (F_1(r) +1)
	\quad\forall\,r\in\erre\,,
\end{align}
have to be prescribed for $F$,
where $\partial F_1^0(r)$ stands for the element of 
$\partial F_1(r)$ having minimum modulus.
This implies that for every $z\in H$ and $w\in\partial F_1(z)$ 

it holds
\begin{align*}
	\iO |w| \leq C_F \iO (F_1(z)+1)\,.
\end{align*}
Let us point out that the above requirement is met by 
all the regular potentials, everywhere defined on the real line, 
with polynomial or first-order exponential growth-rate.
The next two results deal with the asymptotic behavior as $\eps\searrow0$
and the respective error estimate: as a by-product,
these yield existence and uniqueness of solutions, as well as continuous
dependence on the data, for the system 
\eqref{eq1}--\eqref{eq5} with $\eps=0$.

\begin{thm}[Asymptotics: $\eps\searrow0$]
   \label{thm6}
  Assume {\bf A1}--{\bf A5}, \eqref{pol_growth}, and let \luca{$\tau\in(0,\tau_0)$}, and 
  $\chia=0$.
   Suppose also that
  \beq\label{init_eps0}
  \varphi_{0,\tau}\in V\,, \qquad F(\varphi_{0,\tau})\in L^1(\Omega)\,,\qquad \sigma_{0,\tau}\in H\,.
  \eeq
  For every \luca{$\eps\in(0, \eps_0)$}, let 
  the initial data $(\varphi_{0,\eps\tau}, \mu_{0,\eps\tau}, \sigma_{0,\eps\tau})$
  satisfy assumptions \eqref{ip_init} and \eqref{ip_infty},
  and denote by
  $(\ph_\ept, \mu_\ept, \sigma_\ept, \xi_\ept)$ the respective unique weak solution to the system
  \eqref{eq1}-\eqref{eq5} 
  obtained from Theorem~\ref{thm1}.
  In addition, we assume that, as $\eps\searrow0$,
  \beq
  \label{init2_eps0}
  \varphi_{0,\ept}\wto\varphi_{0,\tau} \quad\text{in } V\,, \qquad
  \sigma_{0,\ept}\to\sigma_{0,\tau} \quad\text{in } H\,,
  \eeq
  and
  \beq\label{init3_eps0}
  \exists\,M_0>0:\quad\eps^{1/2}\norm{\mu_{0,\ept}}_H 
  + \norm{F(\varphi_{0,\ept})}_{L^1(\Omega)}\leq M_0
  \quad\forall\,\luca{\eps\in \big(0,\eps_0\big)}\,.
  \eeq
  Then, there exists a quadruplet 
  $(\ph_\tau,\mu_\tau, \sigma_\tau, \xi_\tau)$, with
  \begin{align*}
	&\ph_\tau \in \H1 H \cap \L\infty V\,,
	\qquad
  	\mu_\tau \in \L2 {V}\,, \\
  	&\sigma_\tau \in \H1 {V^*} \cap \L2 V\cap L^\infty(Q)\,, \qquad
	0\leq\sigma_\tau(t,x)\leq1 \quad\text{for a.e.~}x\in\Omega\,,\quad\forall\,t\in[0,T]\,,\\
   	&\xi_\tau  \in \L2 H\,,
\end{align*}
  such that 
  \begin{align*}
  & \<\partial_t \varphi_\tau, \zeta>
   + \iO\nabla  \mu_\tau \cdot \nabla \zeta
  = \iO(P\sigma_\tau - A)h(\varphi_\tau)\zeta\,, \\
  & \<\partial_t\sigma_\tau,\zeta>
  + \iO \nabla\sigma_\tau \cdot \nabla \zeta
  + B \iO (\sigma_\tau-\sigma_S) \zeta
  + C\iO \sigma_\tau h(\varphi_\tau)\zeta
  = 0\,,
\end{align*}
for every $\zeta\in V$, almost everywhere in $(0,T)$, and
\begin{align*}
    &\mu_\tau=\tau\partial_t\varphi_\tau
    + a\varphi_\tau - J*\varphi_\tau + \xi_\tau + F_2'(\varphi_\tau) - \chi\sigma_\tau\,, \qquad
    \xi_\tau\in\partial F_1(\varphi_\tau) \qquad\text{a.e.~in } Q\,,\\
    &\varphi_\tau(0)=\varphi_{0,\tau}\,, \quad
    \sigma_\tau(0)=\sigma_{0,\tau} \qquad\text{a.e.~in } \Omega\,.
\end{align*}
  Moreover, as $\eps\searrow0$, along a non-relabelled subsequence it holds that 
\begin{align}
    \label{epstozero:1}
  	\ph_\ept  \wstarto \ph_\tau \quad &\hbox{in } \H1 H \cap \L\infty V\,,\\
  	\label{epstozero:2}
	\mu_\ept  \wto \mu_\tau \quad &\hbox{in } \L2 V\,,\\
  	\label{epstozero:3}
	\sigma_\ept  \wstarto \sigma_\tau \quad &\hbox{in } \H1 {V^*}  \cap \L2 V \cap L^\infty(Q)\,,\\
	\label{epstozero:4}
	\xi_\ept\wto\xi_\tau \quad&\text{in } L^2(0,T; H)\,,\\
  	\label{epstozero:5}
  \eps \mu_\ept  \to 0 \quad &\hbox{in } C^0([0,T]; H)\cap \L2 V\,,
\end{align}
\Accorpa\epstozero {epstozero:1} {epstozero:5}
hence in particular that
\begin{align}
	\label{epstozero:strong}
	\ph_\ept\to\ph_\tau \quad\text{in } C^0([0,T]; H)\,, \qquad
  \sigma_\ept\to\sigma_\tau \quad\text{in } C^0([0,T]; V^*)\cap L^2(0,T; H)\,.
\end{align}  
\end{thm}
\begin{thm}[Error estimate: $\eps\searrow0$]
  \label{thm7}
  In the setting of Theorem~\ref{thm6}, \luca{if also
  \beq\label{pot_reg}
  F\in C^1(\erre)\,, \qquad
  |F'(r)-F'(s)|\leq C_F(1+|r|^2 + |s|^2)|r-s| \quad\forall\,r,s\in\erre\,,
  \eeq
  then the solution $(\varphi_\tau, \mu_\tau, \sigma_\tau,\xi_\tau)$ 
  to the system \eqref{eq1}--\eqref{eq5} with $\eps=0$ is unique.
  Moreover, suppose that there exists $M_0>0$ such that 
  \beq
    \label{init4_eps0}
    \eps^{1/4}\left(\norm{\mu_{0,\ept}}_V + 
    \norm{\sigma_{0,\ept}}_V
    + \norm{F'(\varphi_{0,\ept})}_H\right)\leq M_0 \quad\forall\,
    \luca{\eps\in\big(0,\eps_0\big)}\,.
  \eeq
  Then,
  the convergences obtained in Theorem~\ref{thm6} hold along 
  the entire sequence $\eps\searrow0$, and 
  there exists $K_\tau>0$, independent of $\eps$, such that 
  the following error estimate holds:}
  \begin{align*}
	& \norma{\ph_\ept - \ph_\tau}_{\L\infty H}
	+ \norma{\mu_\ept - \mu_\tau}_{\L2 V}
	+ \norma{\sigma_\ept - \sigma_\tau}_{\L\infty H \cap \L2 V}
	\\ &
	\leq K_\tau\left(\eps^{1/4}
	+\norm{\varphi_{0,\ept}-\varphi_{0,\tau}}_H + 
	\norm{\sigma_{0,\ept}-\sigma_{0,\tau}}_H\right)\,.
  \end{align*}
\end{thm}
\begin{remark}
  Note that given $(\varphi_{0,\tau}, \sigma_{0,\tau})$ satisfying \eqref{init_eps0},
  a natural choice for the approximating sequence of initial data $(\varphi_{0,\ept}, \sigma_{0,\ept})$
  satisfying \eqref{init2_eps0}--\eqref{init3_eps0} and \eqref{init4_eps0}
  is given by the solutions to the elliptic problems
  \[
  \varphi_{0,\ept} + \eps^{1/2}\R\varphi_{0,\ept}=\varphi_{0,\tau}\,, \qquad
  \sigma_{0,\ept} + \eps^{1/2}\R\sigma_{0,\ept}=\sigma_{0,\tau}\,.
  \]
  In this case, if for example $\sigma_{0,\tau}\in V$, it is immediate to check that
  \[
  \norm{\varphi_{0,\ept}-\varphi_{0,\tau}}_H+\norm{\sigma_{0,\ept}-\sigma_{0,\tau}}_H
  \leq M_0\eps^{1/4}
  \]
  for a certain $M_0>0$, so that the rate of convergence given by Theorem~\ref{thm7}
  is exactly $1/4$.
\end{remark}

The second asymptotic study that we are going to address is the one as
$\tau\searrow0$, when $\eps>0$ is fixed. 
In this case, the presence of the parabolic regularization on $\mu$
provided by $\eps>0$ allows considering also very general 
potentials and to avoid assumptions as \eqref{pol_growth}.
The limit as $\tau\searrow0$ corresponds instead to
a vanishing viscosity argument on the system in consideration.
We expect then to lose, at the limit $\tau=0$, time regularity
on the solutions, as well as the separation principle.
The next two results deal with the asymptotic behavior as $\tau\searrow0$
and the respective error estimate: again, as a by-product,
these yield existence and uniqueness of solutions for the system 
\eqref{eq1}--\eqref{eq5} with $\tau=0$.

\begin{thm}[Asymptotics: $\tau\searrow0$]
   \label{thm8}
   Assume {\bf A1}--{\bf A5}, \luca{$\eps\in(0,\eps_0)$}, and
   \beq
   \label{ip_chi}
   \luca{
   0\leq \chi < \sqrt{c_a}\,, \qquad
   (\chi+\chia+4c_a\chi)^2<8c_aC_0+4\chi\chia\,.}
   \eeq
   Moreover, let us suppose that
  \beq\label{init_tau0}
  \varphi_{0,\eps},\mu_{0,\eps},\sigma_{0,\eps} \in H\,,\qquad F(\varphi_{0,\eps})\in L^1(\Omega)\,.
  \eeq
  For every \luca{$\tau\in(0,\tau_0)$}, let 
  the initial data $(\varphi_{0,\eps\tau}, \mu_{0,\eps\tau}, \sigma_{0,\eps\tau})$
  satisfy \eqref{ip_init},
  and denote by
  $(\ph_\ept,\mu_\ept, \sigma_\ept, \xi_\ept)$ the corresponding weak solution to
  \eqref{eq1}-\eqref{eq5} 
  obtained from Theorem \ref{thm1}.
  Suppose also that, as $\tau\searrow0$,
  \beq
  \label{init2_tau0}
  \varphi_{0,\ept}\to\varphi_{0,\eps} \quad\text{in } H\,, \qquad
  \mu_{0,\ept}\to\mu_{0,\eps} \quad\text{in } H\,, \qquad
  \sigma_{0,\ept}\to\sigma_{0,\eps} \quad\text{in } H\,,
  \eeq
  and
  \beq\label{init3_tau0}
  \exists\,M_0>0:\quad\tau^{1/2}\norm{\varphi_{0,\ept}}_V 
  + \norm{F(\varphi_{0,\ept})}_{L^1(\Omega)}\leq M_0
  \quad\luca{\forall\,\tau\in(0,\tau_0)}\,.
  \eeq
  Then, there exists a quadruplet 
  $(\ph_\eps,\mu_\eps, \sigma_\eps, \xi_\eps)$, with
  \begin{align*}
	&\ph_\eps,\mu_\eps \in L^\infty(0,T; H) \cap L^2(0,T; V)\,,
	\qquad
  	\eps\mu_\eps + \varphi_\eps \in H^1(0,T; V^*)\cap \L2 V\,, \\
  	&\sigma_\eps \in \H1 {V^*} \cap \L2 V\,,\qquad\xi_\eps  \in \L2 V\,,
\end{align*}
  such that 
  \begin{align*}
  & \<\partial_t(\eps\mu_\eps+\varphi_\eps), \zeta>
   + \iO\nabla  \mu_\eps \cdot \nabla \zeta
  = \iO(P\sigma_\eps - A)h(\varphi_\eps)\zeta\,, \\
  & \<\partial_t\sigma_\eps,\zeta>
  + \iO \nabla \sigma_\eps \cdot \nabla \zeta
  + B \iO (\sigma_\eps-\sigma_S) \zeta
  + C\iO \sigma_\eps h(\varphi_\eps)\zeta
  = \chia\int_\Omega\nabla\varphi_\eps\cdot\nabla\zeta\,,
\end{align*}
for every $\zeta\in V$, almost everywhere in $(0,T)$, and
\begin{align*}
    &\mu_\eps=a\varphi_\eps - J*\varphi_\eps + \xi_\eps + F_2'(\varphi_\eps) 
    - \chi\sigma_\eps\,, \qquad
    \xi_\eps\in\partial F_1(\varphi_\eps) \qquad\text{a.e.~in } Q\,,\\
    &\varphi_\eps(0)=\varphi_{0,\eps}\,, \quad
    \sigma_\eps(0)=\sigma_{0,\eps} \qquad\text{a.e.~in } \Omega\,.
\end{align*}
  Moreover, as $\tau\searrow0$, along a non-relabelled subsequence it holds that 
  \begin{align}
  	\label{tautozero:1}
  	\ph_\ept  \wstarto \ph_\eps \quad &\hbox{in } L^\infty(0,T; H) \cap \L2 V\,,\\
  	\label{tautozero:2}
  	\mu_\ept  \wstarto \mu_\eps \quad &\hbox{in } L^\infty(0,T; H) \cap \L2 V\,,\\
	\label{tautozero:3}
	\eps\mu_\ept+\varphi_{\ept}\wto\eps\mu_\eps+\varphi_\eps
	\quad&\text{in } H^1(0,T; V^*)\cap L^2(0,T; V)\,,\\
  	\label{tautozero:4}
  	\sigma_\ept  \wto \sigma_\eps \quad &\hbox{in } \H1 {V^*}  \cap \L2 V\,,\\
	\label{tautozero:5}
	\xi_\ept\wto\xi_\eps \quad&\luca{\text{in } L^2(0,T; H)}\,,\\
  	\label{tautozero:6}
  	\tau\varphi_\ept\to0\quad&\text{in } H^1(0,T; H)\cap L^\infty(0,T; V)\,,
  \end{align}
  \Accorpa\tautozero {tautozero:1} {tautozero:6}
  \luca{and also} that 
  \begin{align}
  \label{tautozero:strong:1}
  \ph_\ept&\to\ph_\eps \quad\text{in } L^2(0,T; H)\,, \qquad
  \mu_\ept\to\mu_\eps\quad\text{in } L^2(0,T; H)\,,\\
  \label{tautozero:strong:2}
  \sigma_\ept&\to\sigma_\eps \quad\text{in } C^0([0,T]; V^*)\cap L^2(0,T; H)\,.
  \end{align}
  Furthermore, if $\chia=0$ and $\sigma_{0,\ept}$ satisfies \eqref{ip_infty} for all $\tau>0$, 
  then the limit $\sigma_\eps$ satisfies \eqref{max_sigma} as well,
  and
  \[
  \sigma_\ept\wstarto\sigma_\eps \quad\text{in } L^\infty(Q)\,.
  \]
\end{thm}

\begin{thm}[Error estimate: $\tau\searrow0$]
\label{thm9}
  \luca{In the setting of Theorem~\ref{thm8}, suppose that $\chia=0$.}
  Then the solution $(\varphi_\eps, \mu_\eps, \sigma_\eps, \xi_\eps)$ 
  to the system \eqref{eq1}--\eqref{eq5} with $\tau=0$ is unique,
  the convergences obtained in Theorem~\ref{thm8} hold along 
  the entire sequence $\tau\searrow0$, and 
  there exists $K_\eps>0$, independent of $\tau$, such that 
  the following error estimate holds:
  \luca{
  \begin{align*}
 	& \norm{(\eps\mu_\ept+\varphi_\ept)-(\eps\mu_\eps+\varphi_\eps)}_{L^\infty(0,T; V^*)}
	+\norma{\ph_\ept - \ph_\eps}_{\L2 H}
	+ \norma{\mu_{\ept}-\mu_{\eps}}_{\L2 H}
	\\ & \qquad +\norma{\sigma_\ept - \sigma_\eps}_{\L\infty H\cap \L2 V}\\
	&\leq K_\eps\left(\tau^{1/2}+
	\norm{(\eps\mu_{0,\ept}+\varphi_{0,\ept})
	-(\eps\mu_{0,\eps}+\varphi_{0,\eps})}_{V^*}+
	\norm{\sigma_{0,\ept}-\sigma_{0,\eps}}_H\right)\,.
\end{align*}}
\end{thm}
\begin{remark}
  \luca{Note that given $(\varphi_{0,\eps}, \mu_{0,\eps}, \sigma_{0,\eps})$
  satisfying \eqref{init_tau0}, a natural choice for the approximating sequence 
  $(\varphi_{0,\ept}, \mu_{0,\ept}, \sigma_{0,\ept})$
  is given by the solutions to the elliptic problems
  \[
  \varphi_{0,\ept} + \tau\R\varphi_{0,\ept}=\varphi_{0,\eps}\,, \qquad
  \mu_{0,\ept} + \tau\R\mu_{0,\ept}=\mu_{0,\eps}\,, \qquad
  \sigma_{0,\ept} + \tau\R\sigma_{0,\ept}=\sigma_{0,\eps}\,.
  \]
  In such a case, hypotheses \eqref{init2_tau0}--\eqref{init3_tau0} are readily satisfied.
  Moreover, if for example $\varphi_{0,\eps}, \mu_{0,\eps}, \tau_{0,\eps}\in V$,
  it is immediate to check that, there is $M_0>0$, independent of $\tau$, such that 
  \[
  \norm{\varphi_{0,\ept}-\varphi_{0,\eps}}_H +
  \norm{\mu_{0,\ept}-\mu_{0,\eps}}_H+ 
  \norm{\sigma_{0,\ept}-\sigma_{0,\eps}}_H \leq M_0\tau^{1/2}\,,
  \]
  so that the rate of convergence given by Theorem~\ref{thm9} is exactly $1/2$.}
\end{remark}

The last two results we present deal with the asymptotic
study of the system \eqref{eq1}--\eqref{eq5} as the parameters 
$\eps$ and $\tau$ go to $0$ simultaneously.
Again, as a by-product, these yield existence and uniqueness 
of solutions for the limit system \eqref{eq1}--\eqref{eq5}
with $\eps=\tau=0$.

\begin{thm}[Asymptotics: $\eps,\tau\searrow0$]
\label{thm10}
Assume {\bf A1}--{\bf A5}, \eqref{pol_growth},
\eqref{ip_chi}, $\chia=0$, and suppose that
\beq\label{init1_ept}
  \varphi_0, \sigma_0 \in H\,, \qquad F(\varphi_0)\in L^1(\Omega)\,.
\eeq
\luca{For every $\eps\in(0,\eps_0)$ and $\tau\in(0,\tau_0)$,} let 
the initial data $(\varphi_{0,\ept}, \mu_{0,\ept}, \sigma_{0,\ept})$
satisfy \eqref{ip_init} and \eqref{ip_infty}, 
and denote by $(\ph_\ept,\mu_\ept, \sigma_\ept, \xi_\ept)$ 
the respective unique weak solution to the system
\eqref{eq1}-\eqref{eq5} obtained from Theorem \ref{thm1}.
Suppose also that, as $(\eps,\tau)\to(0,0)$,
  \beq
  \label{init2_ept}
  \varphi_{0,\ept}\to\varphi_0 \quad\text{in } H\,, \qquad
  \sigma_{0,\ept}\to\sigma_0 \quad\text{in } H\,,
  \eeq
  and that there exists $M_0>0$ such that 
  \beq\label{init3_ept}
  \tau^{1/2}\norm{\varphi_{0,\ept}}_V +
  \eps^{1/2}\norm{\mu_{0,\ept}}_H
  + \norm{F(\varphi_{0,\ept})}_{L^1(\Omega)}\leq M_0
  \qquad\luca{\forall\,(\eps,\tau)\in
  \big(0,\eps_0\big)
   \times(0,\tau_0)}\,.
  \eeq
  Then, there exists a quadruplet $(\varphi,\mu,\sigma,\xi)$, with 
  \begin{align*}
	&\ph \in H^1(0,T; V^*) \cap L^2(0,T; V)\,,\\
         &\mu =a\varphi - J*\varphi + \xi + F_2'(\varphi) - \chi\sigma
         \in L^2(0,T; V)\,, \\
  	&\sigma \in \H1 {V^*} \cap \L2 V\cap L^\infty(Q)\,, \qquad
	0\leq\sigma(t,x)\leq1 \quad\text{for a.e.~}x\in\Omega\,,\quad\forall\,t\in[0,T]\,,\\
	&\xi  \in \L2 V\,,\qquad
	\xi\in\partial F_1(\varphi) \qquad\text{a.e.~in } Q\,,\\
	&\varphi(0)=\varphi_{0}\,, \quad
        \sigma(0)=\sigma_0 \qquad\text{a.e.~in } \Omega\,,
\end{align*}
  such that, for every $\zeta\in V$, almost everywhere in $(0,T)$, it holds 
\begin{align*}
  & \<\partial_t\varphi, \zeta>
   + \iO\nabla  \mu \cdot \nabla \zeta
  = \iO(P\sigma - A)h(\varphi)\zeta\,, \\
  & \<\partial_t\sigma,\zeta>
  + \iO \nabla \sigma\cdot \nabla \zeta
  + B \iO (\sigma-\sigma_S) \zeta
  + C\iO \sigma h(\varphi)\zeta
  = 0\,.
\end{align*}
  Moreover, as $(\eps,\tau)\to(0,0)$, 
  along a non-relabelled subsequence it holds that 
  \begin{align}
  	\label{epstautozero:1}
  	\ph_\ept  \wstarto \ph \quad &\hbox{in } L^\infty(0,T; H) \cap \L2 V\,,\\
  	\label{epstautozero:2}
	\mu_\ept  \wto \mu \quad &\hbox{in } \L2 V\,,\\
	\label{epstautozero:3}
	\eps\mu_\ept+\varphi_{\ept}\wto\varphi
	\quad&\text{in } H^1(0,T; V^*)\cap L^2(0,T; V)\,,\\
  	\label{epstautozero:4}
  	\sigma_\ept  \wstarto \sigma \quad &\hbox{in } \H1 {V^*}  \cap \L2 V \cap L^\infty(Q)\,,\\
	\label{epstautozero:5}
	\eps\mu_\ept\to0 \quad&\text{in } C^0([0,T]; H) \cap \L2 V\,,\\
  	\label{epstautozero:6}
  	\tau\varphi_\ept\to0\quad&\text{in } H^1(0,T; H)\cap L^\infty(0,T; V)\,,
  \end{align}
  \Accorpa\epstautozero {epstautozero:1} {epstautozero:6}
  hence in particular that 
  \begin{align}
  \label{epstautozero:strong}
  &\ph_\ept\to \ph\quad\text{in } L^2(0,T; H)\,, \qquad
  \sigma_\ept\to \sigma \quad\text{in } C^0([0,T]; V^*)\cap L^2(0,T; H)\,.
  \end{align}
\end{thm}

\begin{thm}[Error estimate: $\eps,\tau\searrow0$]
\label{thm11}
\luca{In the setting of Theorem~\ref{thm10},
assume \eqref{pot_reg}, and
suppose also that there exist constants $c_F, M_0>0$ such that
\beq
  \label{pot_reg2}
  F(r)\geq c_F|r|^4 - c_F^{-1} \quad\forall\,r\in\erre\,.
\eeq
Then, the solution $(\varphi, \mu, \sigma, \xi)$ 
to the system \eqref{eq1}--\eqref{eq5} with $\eps=\tau=0$ is unique.
Moreover, let the approximating initial data verify
\eqref{init3_ept} and, for every $(\eps, \tau) \in (0,\eps_0)\times(0,\tau_0)$,
\beq
  \label{init4_ept}
  \frac{\eps^{1/4}}{\tau^{1/2}}\left(\norm{\mu_{0,\ept}}_H+\norm{F'(\varphi_{0,\ept})}_H\right)
  +\eps^{1/4}\left(\norm{\mu_{0,\ept}}_V + \norm{\sigma_{0,\ept}}_V\right) \leq M_0\,.
\eeq
Then,} the convergences obtained in 
Theorem~\ref{thm10} hold along every 
subsequence $(\eps_k,\tau_k)_k$ satisfying 
\beq\label{limsup}
  \limsup_{k\to\infty}\frac{\eps_k^{1/2}}{\tau_k}<+\infty\,,
\eeq
and in this case, there exists $K>0$, independent of $k$, such that 
the following error estimate holds:
\begin{align*}
	& \norma{\ph_{\eps_k\tau_k} - \ph}_{C^0([0,T]; V^*)\cap \L2 H}
	+ \norma{\sigma_{\eps_k\tau_k} - \sigma}_{C^0([0,T]; H) \cap \L2 V}
	\\ &
	\quad \leq K\left(\eps_k^{1/4} + \tau_k^{1/2}
	+\norm{\varphi_{0,\eps_k\tau_k}-\varphi_{0}}_{V^*} + 
	\norm{\sigma_{0,\eps_k\tau_k}-\sigma_{0}}_H\right)\,.
  \end{align*}
\end{thm}

Throughout the paper we convey 
to use the symbol $M$ to indicate constants depending only
on structural data. So, its meaning may change from line to line without 
further comments. Moreover, we will sometimes add a self-explanatory
subscript to stress its possible dependence.


\section{Analysis of the system with $\eps,\tau>0$}
\label{sec:WP}
This section is devoted to the proof of the results concerning 
the behavior of the system with $\eps,\tau>0$, namely the
existence of weak solutions contained in Theorem~\ref{thm1},
the continuous dependence result contained in Theorem~\ref{thm2},
the regularity property of Theorem~\ref{thm3}, 
the existence of strong solution and separation in Theorem~\ref{thm4},
and the refined continuous dependence result in Theorem~\ref{thm5}.
Let us recall that throughout this section $\eps,\tau>0$ are fixed. 

\subsection{The approximation}
To prove the existence of solutions we rely on an approximation procedure
based on the two parameters $n\in\enne$ and $\lambda>0$,
involving a Faedo--Galerkin approximation on the functional spaces
and the Yosida approximation on the potential (c.f. {\bf A4}), respectively.

Let $(e_j)_{j\in\enne}$ and $(l_j)_{j\in\enne}$ 
be the sequences of eigenfunctions and eigenvalues of the operator 
$-\Delta$ with homogeneous Neumann conditions, renormalized in such a way 
that $\norm{e_j}_H=1$ for all $j\in\enne$. Then it is well known that 
$(e_j)_j$ is a complete orthonormal system in $H$, and orthogonal in $V$.
For every $n\in\enne$, let $\Wn:=\operatorname{span}\{e_1,\ldots,e_n\}$,
and define $\Pn:H\to \Wn$ as the orthogonal projection on $\Wn$
with respect to the scalar product  of $H$. Then,
as $n\to\infty$, it holds that 
$\Pn v\to v$ in $H$ (resp.~$V$ or $W$) for every $v\in H$ (resp.~$V$ or $W$).
We consider the following approximated problem:
We then consider the following approximated system: find a triplet $(\ph_{\lambda,n},\mu_{\lambda,n},\sigma_{\lambda,n})$ such that
\begin{align}
  \label{eq1_app}
  &\eps\partial_t\mu_{\lambda,n} + \partial_t \varphi_{\lambda,n} - \Delta\mu_{\lambda,n}
  = {\Pn}[(P\sigma_{\lambda,n} - A)h(\varphi_{\lambda,n})]
  \quad&&\text{in } Q,\\
  \label{eq2_app}
  &\mu_{\lambda,n}=\tau\partial_t\varphi_{\lambda,n} 
  + a\varphi_{\lambda,n} - J*\varphi_{\lambda,n} + {\Pn}F_\lambda'(\varphi_{\lambda,n}) - \chi\sigma_{\lambda,n}
  \quad&&\text{in } Q,\\
  \label{eq3_app}
  &\partial_t\sigma_{\lambda,n} - \Delta\sigma_{\lambda,n} 
  + B(\sigma_{\lambda,n}-{\sigma_{{S,n}}}) + {\Pn}[C\sigma_{\lambda,n} h(\varphi_{\lambda,n})]
  =-\eta\Delta\varphi_{\lambda,n} 
  \quad&&\text{in } Q,\\
  \label{eq4_app}
  &
  \partial_{\bf n}\mu_{\lambda,n} 
   =\partial_{\bf n}(\sigma_{\lambda,n}-\eta\varphi_{\lambda,n})
  = 0
  \quad&&\text{on } \Sigma,\\
  \label{eq5_app}
  &\mu_{\lambda,n}(0) = {\Pn}\mu_0, \quad
  \varphi_{\lambda,n}(0)={\Pn}\varphi_0, \quad
  \sigma_{\lambda,n}(0)={\Pn}\sigma_0
  \quad&&\text{in } \Omega,
\end{align}
where $\sigma_{S,n}:= \Pn \sigma_S$, in the form
\begin{align*}
  \varphi_{\lambda,n}(t,x)&:=\sum_{j=1}^n\alpha^{\lambda,n}_j(t)e_j(x)\,, \quad
  \mu_{\lambda,n}(t,x):=\sum_{j=1}^n\beta^{\lambda,n}_j(t)e_j(x)\,,
  \quad 
  \sigma_{\lambda,n}(t,x):=\sum_{j=1}^n\gamma^{\lambda,n}_j(t)e_j(x)\,,
\end{align*}
for $t \in [0,T]$, $x\in\Omega$, and $j \in \{ 1,...,n\}$.
Moreover, let us introduce the vectors
\[
\boldsymbol\alpha^{\lambda,n}, 
\boldsymbol\beta^{\lambda,n}, 
\boldsymbol\gamma^{\lambda,n}:
[0,T]\to \erre^n\,,
\]
by
\begin{align*}
	\boldsymbol\alpha^{\lambda,n}:= 
	(\alpha^{\lambda,n}_1,...,\alpha^{\lambda,n}_n)^T\,, \quad
	\boldsymbol\beta^{\lambda,n}:= 
	(\beta^{\lambda,n}_1,...,\beta^{\lambda,n}_n)^T\,, \quad
	\boldsymbol\gamma^{\lambda,n}:= 
	(\gamma^{\lambda,n}_1,...,\gamma^{\lambda,n}_n)^T\,.
\end{align*}
Plugging these expression in
\eqref{eq1_app}--\eqref{eq5_app}
and taking arbitrary $e_i\in \Wn$ as test functions,
for $i=1,\ldots,n$,
we deduce that $(\varphi_{\lambda,n}, \mu_{\lambda,n}, \sigma_{\lambda,n})$ solves
the approximated system if and only if
$(\boldsymbol\alpha^{\lambda,n}, 
\boldsymbol\beta^{\lambda,n}, 
\boldsymbol\gamma^{\lambda,n})$ solves the following system of ODEs,
for $i=1,\ldots,n$:
\begin{align*}
  &
  \eps\partial_t\beta^{\lambda,n}_i
  + \partial_t \alpha^{\lambda,n}_i 
  +l_i\beta^{\lambda,n}_i
  = \int_\Omega
  \Bigl(P\sum_{j=1}^n\gamma^{\lambda,n}_je_j - A \Bigr)
  h\Bigl(\sum_{j=1}^n\alpha^{\lambda,n}_je_j\Bigr)e_i\,,\\
  &\beta^{\lambda,n}_i=
  \tau\partial_t\alpha^{\lambda,n}_i 
  + \sum_{j=1}^n\alpha^{\lambda,n}_j\int_\Omega ae_je_i
   - \sum_{j=1}^n\alpha^{\lambda,n}_j \int_\Omega(J*e_j)e_i
   + \int_\Omega F_\lambda'\Bigl(\sum_{j=1}^n\alpha^{\lambda,n}_je_j\Bigr)e_i 
   - \chi\gamma^{\lambda,n}_i\,,\\
  &\partial_t\gamma^{\lambda,n}_i
+l_i\gamma^{\lambda,n}_i 
  + B\Bigl(\gamma^{\lambda,n}_i-\int_\Omega \sigma_{S,n} e_i\Bigr) 
  + C\sum_{j=1}^n\gamma^{\lambda,n}_j
  \int_\Omega h\Bigl(\sum_{m=1}^n\alpha_m^{\lambda,n}e_m\Bigr)e_je_i
  = \chia l_i\alpha^{\lambda,n}_i\,,\\
  &\alpha^{\lambda,n}_i(0)=(\varphi_0,e_i)_H\,, \quad
  \beta^{\lambda,n}_i(0) = (\mu_0,e_i)_H\,, \quad
  \gamma^{\lambda,n}_i(0)=(\sigma_0,e_i)_H\,.
\end{align*}
Since $h,F_\lambda':\erre\to\erre$ are Lipschitz continuous and $h$ is bounded, such initial value
system can be written in the form
\[
\begin{cases}
\partial_t 
(\boldsymbol\alpha^{\lambda,n}, 
\boldsymbol\beta^{\lambda,n}, 
\boldsymbol\gamma^{\lambda,n}) &=
g_{\lambda,n}(\boldsymbol\alpha^{\lambda,n}, 
\boldsymbol\beta^{\lambda,n}, 
\boldsymbol\gamma^{\lambda,n})\,,\\
(\boldsymbol\alpha^{\lambda,n}, 
\boldsymbol\beta^{\lambda,n}, 
\boldsymbol\gamma^{\lambda,n})(0)&=((\varphi_0,e_i)_H,(\mu_0,e_i)_H,(\sigma_0,e_i)_H)\,,
\end{cases}
\]
where $g_{\lambda,n}:\erre^{3n}\to\erre^{3n}$ is locally Lipschitz continuous 
and linearly bounded. Hence, by the Cauchy--Peano theorem, the system above
admits a unique global solution 
$\boldsymbol\alpha^{\lambda,n}, 
  \boldsymbol\beta^{\lambda,n},
  \boldsymbol\gamma^{\lambda,n} \in C^1([0,T]; \erre^n)
$,
implying that
\[
  \varphi_{\lambda,n}, \mu_{\lambda,n}, \sigma_{\lambda,n} \in C^1([0,T]; \Wn)
\]
are the unique solutions to the approximated problem \eqref{eq1_app}--\eqref{eq5_app}.

\subsection{Uniform estimates}
\label{ssec:unif_est}
We prove uniform estimates independent of $\lambda$ and $n$,
still keeping $\eps,\tau>0$ fixed.

Testing \eqref{eq1_app} by $\mu_{\lambda,n}$,
\eqref{eq2_app} by $-\partial_t\varphi_{\lambda,n}$,
\eqref{eq3_app} by $\sigma_{\lambda,n}$, taking the sum
and integrating over $(0,t)$, yields
by symmetry of the kernel $J$, for every $t\in[0,T]$,
\begin{align*}
  &\frac\eps2\norm{\mu_{\lambda,n}(t)}_H^2 + \int_{Q_t}|\nabla\mu_{\lambda,n}|^2
  +\tau\int_{Q_t}|\partial_t\varphi_{\lambda,n}|^2 + 
  \frac14\int_{\Omega\times\Omega} 
  J(x-y)|\varphi_{\lambda,n}(t,x)-\varphi_{\lambda,n}(y,t)|^2\,\d x\,\d y\\
  &\qquad+\int_\Omega F_\lambda(\varphi_{\lambda,n}(t))
  +\frac12\norm{\sigma_{\lambda,n}(t)}_H^2 + \int_{Q_t}|\nabla\sigma_{\lambda,n}|^2
  +\int_{Q_t}\left(B+ Ch(\varphi_{\lambda,n})\right)|\sigma_{\lambda,n}|^2\\
  &=\frac\eps2\norm{\Pn\mu_0}_H^2 +
  \frac14\int_{\Omega\times\Omega} J(x-y)|\Pn\varphi_0(x)-\Pn\varphi_0(y)|^2\,\d x\,\d y+
  \int_\Omega F_\lambda(\Pn\varphi_0) +
   \frac12\norm{\Pn\sigma_0}_H^2\\
  &\qquad+\int_{Q_t}(P\sigma_{\lambda,n} - A )h(\varphi_{\lambda,n})\mu_{{\lambda,n}}
  +\chi\int_{Q_t}\sigma_{\lambda,n} \partial_t\varphi_{\lambda,n} 
  +B\int_{Q_t}\sigma_{S,n} \sigma_{\lambda,n} + 
  \chia\int_{Q_t}\nabla\varphi_{\lambda,n}\cdot\nabla\sigma_{\lambda,n}\,.
\end{align*}
Now, note that by assumption {\bf A5} we have
\begin{align}
  &\notag
  \frac14\int_{\Omega\times\Omega} 
  J(x-y)|\varphi_{\lambda,n}(t,x)-\varphi_{\lambda,n}(t,y)|^2\,\d x\,\d y
  =\frac12\int_\Omega\left[a(x)|\varphi_{\lambda,n}|^2 
  - (J*\varphi_{\lambda,n})\varphi_{\lambda,n}\right](t,x)\,\d x\\
  &\label{est:ca}
  \geq \frac{a_*}2\norm{\varphi_{\lambda,n}(t)}_H^2 - 
  \frac12\norm{J*\varphi_{\lambda,n}(t)}_H
  \norm{\varphi_{\lambda,n}(t)}_H\geq\frac{a_*-a^*}2\norm{\varphi_{\lambda,n}(t)}_H^2\,,
\end{align}
and similarly
\begin{align*}
  \frac14\int_{\Omega\times\Omega} J(x-y)|\Pn\varphi_0(x)-\Pn\varphi_0(y)|^2\,\d x\,\d y&=
  \frac12\int_\Omega\left[a|\Pn\varphi_{0}|^2 
  - (J*\Pn\varphi_{0})\Pn\varphi_{0}\right](x)\,\d x\\
  &
  \leq\frac{a^*+a^*}2\norm{\Pn\varphi_0}_H^2\leq a^*\norm{\varphi_0}_H^2\,.
\end{align*}
Using that $F_{\lambda}\geq0$, \eqref{est:ca}
along with the definition of $c_a$, recalling also that 
$h$ is non-negative and bounded and that $\Pn$ is a contraction on $H$,
owing to the Young inequality we infer that
\begin{align}
  \notag&\frac\eps2\norm{\mu_{\lambda,n}(t)}_H^2 + \int_{Q_t}|\nabla\mu_{\lambda,n}|^2
  +\tau\int_{Q_t}|\partial_t\varphi_{\lambda,n}|^2
  +\frac12\norm{\sigma_{\lambda,n}(t)}_H^2 + \int_{Q_t}|\nabla\sigma_{\lambda,n}|^2\\
 \notag 
 &\leq\frac\eps2\norm{\mu_0}_H^2 + a^*\norm{\varphi_0}_H^2 +
   \norm{F_\lambda(\Pn\varphi_0)}_{L^1(\Omega)} +
  \frac12\norm{\sigma_0}_H^2
  + \frac{c_a}2\norm{\varphi_{\lambda,n}(t)}_H^2
  +\int_{Q_t}(P\sigma_{\lambda,n}-A)h(\varphi_{\lambda,n})\mu_{\lambda,n}\\
  &\qquad
  +\frac 14\int_{Q_t}|\sigma_{\lambda,n}|^2
 +|Q| B^2\norm{\sigma_{S,n}}^2_{L^\infty(Q)}
 + \chi\int_{Q_t}\sigma_{\lambda,n}\partial_t\varphi_{\lambda,n}
 + \chia\int_{Q_t} \nabla\ph_{\lambda,n}\cdot\nabla\sigma_{\lambda,n}\,.
 \label{est_aux}
\end{align}
Here, we recall that $a^*-a_*\geq0$ which entails that $c_a=\max\{a^*-a_*,1\}>0$.
Then, we test equation \eqref{eq1_app} by $4c_a(\eps\mu_{\lambda,n}+\ph_{\lambda,n})$ and \eqref{eq2_app} 
by $-4c_a\Delta\varphi_{\lambda,n}$, add the resulting equalities
and integrate over $(0,t)$ and by parts, getting, thanks to assumption {\bf A5},
\begin{align*}
  &2c_a\norm{(\eps\mu_{\lambda,n}+\varphi_{\lambda,n})(t)}_H^2 + 
  4c_a\eps\int_{Q_t}|\nabla\mu_{\lambda,n}|^2+
  2c_a\tau\norm{\nabla\varphi_{\lambda,n}(t)}_H^2 + 
  4c_aC_0\int_{Q_t}|\nabla\varphi_{\lambda,n}|^2\\
  \quad &\leq2c_a \norm{ \Pn(\eps\mu_0 + \varphi_0)}_H^2
  + 2c_a\tau\norm{\nabla \Pn\varphi_0}_H^2+
  4c_a\int_{Q_t}(P\sigma_{\lambda,n}-A)
  h(\varphi_{\lambda,n})(\eps\mu_{\lambda,n} + \varphi_{\lambda,n})
  \\ & \qquad +4c_a\chi\int_{Q_t}
  \nabla\sigma_{\lambda,n}\cdot\nabla\varphi_{\lambda,n}
  +8c_a b^*
  \norm{\varphi_{\lambda,n}}_{L^2(Q_t)}
  \norm{\nabla\varphi_{\lambda,n}}_{L^2(Q_t)}\,,
\end{align*}
from which we infer, thanks to the Young inequality and the boundedness of $h$, that 
\begin{align}
  \nonumber
  &2c_a\norm{(\eps\mu_{\lambda,n}+\varphi_{\lambda,n})(t)}_H^2+
  4c_a\eps\int_{Q_t}|\nabla\mu_{\lambda,n}|^2+
  2c_a\tau\norm{\nabla\varphi_{\lambda,n}(t)}_H^2 + 
  2c_aC_0\int_{Q_t}|\nabla\varphi_{\lambda,n}|^2\\
  \nonumber
  &\leq 4c_a\eps^2\norm{\mu_0}_H^2+
  4c_a\norm{\varphi_0}_H^2 + 2c_a\tau\norm{\nabla\varphi_0}_H^2
  +4c_a\chi\int_{Q_t}\nabla\sigma_{\lambda,n}\cdot\nabla\varphi_{\lambda,n}\\
  \label{est_aux2}
  &\qquad+M\left(1
  +\int_{Q_t}|\eps\mu_{\lambda,n}+\varphi_{\lambda,n}|^2+
  \int_{Q_t}|\varphi_{\lambda,n}|^2
  +\int_{Q_t}|\sigma_{\lambda,n}|^2
  \right)\,.
\end{align}
for a constant $M>0$, independent of $\lambda$, $n$, $\eps$, and $\tau$.
Summing \eqref{est_aux} and \eqref{est_aux2}, we infer that, possibly updating $M$,
\begin{align}
  \nonumber
  &\frac\eps2\norm{\mu_{\lambda,n}(t)}_H^2 
  + (1+4c_a\eps)\int_{Q_t}|\nabla\mu_{\lambda,n}|^2
  +\tau\int_{Q_t}|\partial_t\varphi_{\lambda,n}|^2
  +\frac12\norm{\sigma_{\lambda,n}(t)}_H^2 + \int_{Q_t}|\nabla\sigma_{\lambda,n}|^2\\
  \nonumber
  &\qquad+2c_a\norm{(\eps\mu_{\lambda,n}+\varphi_{\lambda,n})(t)}_H^2+
  2c_a\tau\norm{\nabla\varphi_{\lambda,n}(t)}_H^2 + 
  2c_aC_0\int_{Q_t}|\nabla\varphi_{\lambda,n}|^2\\
    \label{est_aux3}
  &\leq \left(\frac\eps2+4c_a\eps^2\right)\norm{\mu_0}_H^2 +
  (a^*+4c_a)\norm{\varphi_0}_H^2+ 2c_a\tau\norm{\nabla\varphi_0}_H^2+
   \norm{F_\lambda(\Pn\varphi_0)}_{L^1(\Omega)}\\ 
   \nonumber
  &\qquad+\frac12\norm{\sigma_0}_H^2
  +\frac{c_a}2\norm{\varphi_{\lambda,n}(t)}_H^2
  +\chi\int_{Q_t}\sigma_{\lambda,n}\partial_t\varphi_{\lambda,n}+
  (\chia+4c_a\chi) \int_{Q_t}\nabla\sigma_{\lambda,n}\cdot\nabla\varphi_{\lambda,n}\\
\notag
 &\qquad+M\left(1+\int_{Q_t}|\eps\mu_{\lambda,n}+\varphi_{\lambda,n}|^2
  + \int_{Q_t}|\varphi_{\lambda,n}|^2
   + \int_{Q_t}|\sigma_{\lambda,n}|^2\right)
   +\int_{Q_t}(P\sigma_{\lambda,n}-A)h(\varphi_{\lambda,n})\mu_{\lambda,n}\,.
\end{align}
Note that 
\[
  \frac{c_a}2\norm{\varphi_{\lambda,n}(t)}_H^2\leq
  c_a\norm{(\eps\mu_{\lambda,n}+\varphi_{\lambda,n})(t)}_H^2
  +c_a\eps^2\norm{\mu_{\lambda,n}(t)}_H^2\,,
\]
where the two terms on the right-hand side can be incorporated 
in the left-hand side of \eqref{est_aux3}
as $2c_a-c_a=c_a>0$ and $\frac\eps2-c_a\eps^2\geq\frac\eps4$ (since $\eps\in(0,\frac1{4c_a})$).
Furthermore, using the Young inequality we have 
\begin{align*}
  & \chi\int_{Q_t}\sigma_{\lambda,n}\partial_t\varphi_{\lambda,n}+
  (\chia+4c_a\chi) \int_{Q_t}\nabla\sigma_{\lambda,n}\cdot\nabla\varphi_{\lambda,n}
  \\ &
  \leq \frac \tau 2 \int_{Q_t}|\partial_t\varphi_{\lambda,n}|^2
  + \frac {\chi^2}{2\tau} \int_{Q_t}|\sigma_{\lambda,n}|^2
  + \frac 12 \int_{Q_t}|\nabla\sigma_{\lambda,n}|^2
  +  \frac {(\chia+4c_a\chi)^2}2 \int_{Q_t}|\nabla\ph_{\lambda,n}|^2\,.
\end{align*}
Collecting the above estimates, we infer that 
\begin{align}
  \nonumber
  &\frac\eps4\norm{\mu_{\lambda,n}(t)}_H^2 
  + (1+4c_a\eps)\int_{Q_t}|\nabla\mu_{\lambda,n}|^2
  +\frac\tau2\int_{Q_t}|\partial_t\varphi_{\lambda,n}|^2
  +\frac12\norm{\sigma_{\lambda,n}(t)}_H^2 + \int_{Q_t}|\nabla\sigma_{\lambda,n}|^2\\
  \nonumber
  &\qquad+c_a\norm{(\eps\mu_{\lambda,n}+\varphi_{\lambda,n})(t)}_H^2+
  2c_a\tau\norm{\nabla\varphi_{\lambda,n}(t)}_H^2 + 
  2c_aC_0\int_{Q_t}|\nabla\varphi_{\lambda,n}|^2\\
  \nonumber
  &\leq \frac32\eps\norm{\mu_0}_H^2 +
  (a^*+4c_a)\norm{\varphi_0}_H^2+ 2c_a\tau\norm{\nabla\varphi_0}_H^2+
   \norm{F_\lambda(\Pn\varphi_0)}_{L^1(\Omega)}
   +\frac12\norm{\sigma_0}_H^2\\ 
   \nonumber
  &\qquad
  +M\left(1+\int_{Q_t}|\eps\mu_{\lambda,n}+\varphi_{\lambda,n}|^2
  + \int_{Q_t}|\varphi_{\lambda,n}|^2
   + \int_{Q_t}|\sigma_{\lambda,n}|^2\right)
   +\frac{\chi^2}{2\tau} \int_{Q_t}|\sigma_{\lambda,n}|^2\\
   \label{est_aux4}
   &\qquad
  +\frac12\int_{Q_t}|\nabla\sigma_{\lambda,n}|^2
  + \frac{(\chia+4c_a\chi)^2}2 \int_{Q_t}|\nabla\ph_{\lambda,n}|^2
  +\int_{Q_t}(P\sigma_{\lambda,n}-A)h(\varphi_{\lambda,n})\mu_{\lambda,n}\,.
\end{align}
Moreover, the last term on the \rhs\ can be easily bounded owing to Young's inequality.

Then, we fix $\lambda>0$, and since $F_\lambda$ has at most quadratic growth (depending on $\lambda$)
and $\varphi_0\in H$, we have that $\norm{F_\lambda(\Pn\varphi_0)}_{L^1(\Omega)}\leq M_\lambda$
uniformly in $n\in\enne$, for a certain $M_\lambda>0$ independent of $n$.
Therefore, Gronwall's lemma yields that 
\beq\label{est}
  \norm{\mu_{\lambda,n}}_{L^\infty(0,T; H)\cap L^2(0,T; V)}^2+
  \norm{\varphi_{\lambda,n}}_{H^1(0,T; H)\cap L^\infty(0,T; V)}^2
  + \norm{\sigma_{\lambda,n}}^2_{L^\infty(0,T; H)\cap L^2(0,T; V)}
  \leq M_\lambda\,,
\eeq
where the constant $M_\lambda$ is independent of $n$ (but not of $\tau$ and $\eps$).
Furthermore, by comparison in equations \eqref{eq1_app} and \eqref{eq3_app}, we 
deduce that
\beq
  \label{est4}
  \norm{\partial_t (\eps\mu_{\lambda,n}+\varphi_{\lambda,n})}^2_{L^2(0,T; V^*)} 
 +\norm{\partial_t\mu_{\lambda,n}}_{L^2(0,T; V^*)}^2
  + \norm{\partial_t\sigma_{\lambda,n}}_{L^2(0,T; V^*)}^2\leq M_\lambda\,.
\eeq

\subsection{Passage to the limit}
We pass now to the limit, keeping $\eps,\tau>0$ fixed, first as $n\to\infty$
and then as $\lambda\searrow0$.
From the estimates \eqref{est}--\eqref{est4} and the Aubin-Lions compactness 
theorems (see, e.g.,~\cite[Cor.~4]{simon}), we deduce that the exists a triplet
$(\varphi_\lambda,\mu_\lambda,\sigma_\lambda)$, with 
\[
  \varphi_\lambda \in H^1(0,T; H)\cap L^\infty(0,T; V)\,,\qquad
  \mu_\lambda,\sigma_\lambda\in H^1(0,T; V^*)\cap L^2(0,T; V)\,,
\]
such that, as $n\to\infty$,
\begin{align*}
  \varphi_{\lambda,n}\wstarto\varphi_{\lambda} \quad&\text{in } H^1(0,T; H)\cap L^\infty(0,T; V)\,,\qquad
  \varphi_{\lambda,n}\to\varphi_{\lambda} &&\text{in }  C^0([0,T]; H)\,,\\
  \mu_{\lambda,n}\wto\mu_{\lambda} \quad&\text{in } H^1(0,T; V^*)\cap L^2(0,T; V)\,,\qquad
  \mu_{\lambda,n}\to\mu_{\lambda} &&\text{in } C^0([0,T]; V^*)\cap L^2(0,T; H)\,,\\
  \sigma_{\lambda,n}\wto\sigma_{\lambda} \quad&\text{in } H^1(0,T; V^*)\cap L^2(0,T; V)\,,\qquad
  \sigma_{\lambda,n}\to\sigma_{\lambda} &&\text{in } C^0([0,T]; V^*)\cap L^2(0,T; H) \,.
\end{align*}
Since $F_\lambda'$ is Lipschitz continuous and $h$ is Lipschitz continuous and bounded,
it is a standard matter to pass the limit in 
the approximated problem \eqref{eq1_app}--\eqref{eq5_app} as $n\to\infty$ to
obtain, for every test function $\zeta\in V$,
\begin{align}
  \label{eq1_lam}
  &\ip{\partial_t(\eps\mu_{\lambda}+ \varphi_\lambda)}\zeta + 
  \int_\Omega\nabla\mu_{\lambda}\cdot\nabla\zeta
  = \int_\Omega(P\sigma_{\lambda} - A)h(\varphi_{\lambda})\zeta\,,\\
  \label{eq2_lam}
  &\mu_{\lambda}=
  \tau  \partial_t\varphi_{\lambda}+
  a\varphi_{\lambda} - J*\varphi_{\lambda} + F_\lambda'(\varphi_{\lambda})-
  \chi\sigma_{\lambda}\,,\\
  \label{eq3_lam}
  &\ip{\partial_t\sigma_{\lambda}}\zeta +
  \int_\Omega\nabla\sigma_{\lambda}\cdot\nabla\zeta
  + \int_\Omega\left[B(\sigma_{\lambda}-\sigma_S) + 
  C\sigma_{\lambda} h(\varphi_{\lambda})\right]\zeta
  =\chia\int_\Omega\nabla\varphi_{\lambda}\cdot\nabla\zeta\,,
\end{align}
almost everywhere in $(0,T)$, and
\beq
  \label{eq5_lam}
  \mu_{\lambda}(0) = \mu_0\,, \quad
  \varphi_{\lambda}(0)=\varphi_0\,, \quad
  \sigma_{\lambda}(0)=\sigma_0
  \qquad\text{a.e.~in } \Omega\,
\eeq
meaning that $(\varphi_\lambda, \mu_\lambda, \sigma_\lambda)$
satisfy the analogous of conditions \eqref{def_mu_xi}--\eqref{var2}
at level $\lambda$.

Clearly, by weak lower semicontinuity of the norms 
and the convex integrands, passing to the $\liminf$
as $n\to\infty$ in the estimates \eqref{est} and \eqref{est4},
and recalling that $F_\lambda\leq F$, we infer that 
there exists $M>0$, independent of $\lambda$
(but not of $\eps$ and $\tau$), such that 
\beq
  \label{est1_lam}
  \norm{\mu_{\lambda}}_{H^1(0,T; V^*)\cap L^2(0,T; V)}^2 +
  \norm{\varphi_{\lambda}}_{H^1(0,T; H)\cap L^\infty(0,T; V)}^2 +
  \norm{\sigma_{\lambda}}^2_{H^1(0,T; V^*)\cap L^2(0,T; V)}
  \leq M\,.
\eeq
Furthermore, the estimate \eqref{est1_lam}
readily implies, by comparison in \eqref{eq2_lam}, that
\beq\label{est2_lam}
  \norm{F_{1,\lambda}'(\varphi_\lambda)}_{L^2(0,T; H)}^2\leq M\,.
\eeq
Hence, there exists a quadruplet $(\varphi,\mu,\sigma,\xi)$, with 
\[
  \varphi \in H^1(0,T; H)\cap L^\infty(0,T; V)\,, \qquad
  \mu,\sigma\in H^1(0,T; V^*)\cap L^2(0,T; V)\,,\qquad
  \xi\in L^2(0,T; H)\,,
\]
such that, as $\lambda\to0$,
\begin{align*}
  \varphi_{\lambda}\wstarto\varphi \quad&\text{in } H^1(0,T; H)\cap L^\infty(0,T; V)\,,\qquad
  \varphi_{\lambda}\to\varphi &&\text{in } C^0([0,T]; H)\,,\\
  \mu_{\lambda}\wto\mu \quad&\text{in } H^1(0,T; V^*)\cap L^2(0,T; V)\,,\qquad
  \mu_{\lambda}\to\mu &&\text{in }  C^0([0,T]; V^*) \cap L^2(0,T; H)\,,\\
  \sigma_{\lambda}\wto\sigma \quad&\text{in } H^1(0,T; V^*)\cap L^2(0,T; V)\,,\qquad
  \sigma_{\lambda}\to\sigma &&\text{in } C^0([0,T]; V^*) \cap L^2(0,T; H)\,,\\
  F_{1,\lambda}'(\varphi_\lambda)\wto\xi \quad&\text{in } L^2(0,T; H)\,.
\end{align*}
\luca{The graph convergence of $F_{1,\lambda}'$ to $\partial F_1$, as $\lambda \to 0$, implies that}
$\xi\in\partial F_1(\varphi)$ almost everywhere in $Q$. Moreover, by the Lipschitz continuity 
of $F_2'$ and $h$, and the boundedness of $h$, we have that 
\[
 h(\varphi_\lambda)\to h(\varphi) \quad\text{in } L^p(Q)\quad\forall\,p\geq1\,, \qquad
   F_2'(\varphi_\lambda)\to F_2'(\varphi)\quad\text{in } L^2(0,T; H)\,.
\]
Consequently, letting $\lambda\to0$ in the variational formulation of \eqref{eq1_lam}--\eqref{eq5_lam},
we obtain exactly \eqref{def_mu_xi}--\eqref{var2} completing the proof
concerning the existence of weak solutions in Theorem~\ref{thm1}.

\subsection{Maximum principle for $\boldsymbol\sigma$} \label{sub:maxpr}
We prove here the last assertion of Theorem~\ref{thm1}, concerning 
a maximum principle for $\sigma$ under the additional requirement that 
$\chia=0$. Testing equation \eqref{var2} by 
$f_+(\sigma):=(\sigma-1)_+$, we have
\[
  \frac12 \norma{f_+(\sigma(t))}_H^2 + \int_{Q_t}f_+'(\sigma)|\nabla\sigma|^2
  +B\int_{Q_t}f_+(\sigma)(\sigma-\sigma_S) + C\int_{Q_t}f_+(\sigma)\sigma h(\varphi) =
 0,
\] 
where we have used the fact that $f_+(\sigma_0)=0$.
Since $f_+$ is non-decreasing and $h$ is non-negative, we infer that the second and
fourth terms on the left-hand side are non-negative so that
\begin{align}
	\frac12 \norma{f_+(\sigma(t))}_H^2 + 
 	 B\int_{Q_t}f_+(\sigma)(\sigma-\sigma_S)  \leq0\,.
	\label{max_pr}
\end{align}
Moreover, 
since $\sigma_S \leq 1$ by assumption {\bf A3}, 
we have that
\[
  B\int_{Q_t}f_+(\sigma)(\sigma-\sigma_S) =
  B\int_{Q_t\cap\{\sigma>1\}}(\sigma-1)(\sigma-\sigma_S)\geq0\,.
\]
Therefore, coming back to \eqref{max_pr}, we realize that
$f_+(\sigma(t))=0$ which gives us the upper bound $\sigma(t)\leq 1$ a.e~in $\Omega$, 
for every $t\in[0,T]$, as desired. The lower inequality follows by a similar argument
testing by $f_-(\sigma):=-\sigma_-$.

\subsection{Continuous dependence}
\label{ssec:cont}
Let us prove here the continuous dependence of Theorem~\ref{thm2}.
To begin with, bearing in mind the notation introduced in Theorem~\ref{thm2}, we set
$\varphi:=\varphi_1-\varphi_2$, $\mu:=\mu_1-\mu_2$, $\sigma:=\sigma_1-\sigma_2$, 
$\xi:=\xi_1-\xi_2$,
$\varphi_0:=\varphi_0^1-\varphi_0^2$, $\mu_0:=\mu_0^1-\mu_0^2$, 
$\sigma_0:=\sigma_0^1-\sigma_0^2$.
Then, we consider the difference of system \sys\ written for the two solutions to obtain
\begin{align}
  \label{cd_1}
  &\eps\partial_t\mu + \partial_t \varphi - \Delta\mu = 
  P\sigma h(\varphi_1) + (P\sigma_2-A)(h(\varphi_{1})-h(\varphi_2))
  \qquad&&\text{in } Q\,,\\
  \label{cd_2}
  &\mu=\tau\partial_t\varphi + a\varphi - J*\varphi + \xi+F_2'(\varphi_1)-F_2'(\varphi_2) - \chi\sigma
  \qquad&&\text{in } Q\,,\\
  \label{cd_3}
  &\partial_t\sigma - \Delta\sigma + B\sigma + C\sigma h(\varphi_1) = C\sigma_2(h(\varphi_2)-h(\varphi_1))
  -\chia\Delta\varphi 
  \qquad&&\text{in } Q\,,\\
  \label{cd_4}
  &
  \partial_{\bf n}\mu =
  \partial_{\bf n}(\sigma-\eta\ph) = 0 \qquad&&\text{on } \Sigma\,,\\
   \label{cd_5}
 &\mu(0) = \mu_0\,, \quad
  \varphi(0)=\varphi_0\,, \quad
  \sigma(0)=\sigma_0
  \qquad&&\text{in } \Omega\,.
\end{align}
\Accorpa\cdsys cd_1 cd_5
\luca{Next, we test the equation \eqref{cd_1} by $\mathcal R^{-1}(\eps\mu+\varphi)$,
\eqref{cd_2} by $-\varphi$, \eqref{cd_3} by $\sigma$, and take the sum to get, after integration on $[0,t]$,
\begin{align}
  &\notag\frac12\norm{(\eps\mu+\varphi)(t)}_{V^*}^2+
  \eps\int_{Q_t}|\mu|^2 
  +\frac\tau2\norm{\varphi(t)}_H^2
  +\int_{Q_t}\left[a|\varphi|^2 + \xi\varphi+
  (F_2'(\varphi_1)-F_2'(\varphi_2))\varphi\right] \\
  &\notag\qquad+\frac12\norm{\sigma(t)}_H^2 + \int_{Q_t}|\nabla\sigma|^2
  +\int_{Q_t}(B+Ch(\varphi_1))|\sigma|^2\\
  &\notag= \frac12\norm{\eps\mu_0+\varphi_0}_{V^*}^2  
  +\frac\tau2\norm{\varphi_0}_H^2  
  + \frac12\norm{\sigma_0}_H^2
  +\int_{Q_t}(\chi\sigma + J*\varphi)\varphi
  +\int_{Q_t}[C\sigma_2(h(\varphi_2)-h(\varphi_1))]\sigma\\
  &\label{last:1}
  \qquad+\int_{Q_t}\big[\mu + P\sigma h(\varphi_1)+ 
  (P\sigma_2-A)(h(\varphi_{1})-h(\varphi_2))\big]\mathcal R^{-1}(\eps\mu+\varphi)
  + \chia\int_{Q_t}\nabla\ph\cdot\nabla\sigma \,.
\end{align}}
Note that the last term on the left-hand side is non-negative 
due to the positivity of $h$. Hence, using the monotonicity of $\partial F_1$ and
recalling assumption {\bf A5}, we have
\begin{align*}
  \int_{Q_t}\left[a|\varphi|^2 + \xi\varphi+
  (F_2'(\varphi_1)-F_2'(\varphi_2))\varphi\right]
    +\int_{Q_t}(B+Ch(\varphi_1))|\sigma|^2
    \geq 
    C_0\int_{Q_t}|\varphi|^2\,.
\end{align*}
Moreover, under the assumption $\chia=0$, we have, 
owing to \eqref{max_sigma} that $\sigma_2 \in L^\infty(Q)$
with $\norm{\sigma_2}_{L^\infty(Q)}\leq 1$ and that the last term on the \rhs\
of \eqref{last:1} disappears.
Let us estimate the remaining terms on the right-hand side.
\luca{First of all, recalling that $K_0$ denotes the norm of the inclusion $H\embed V^*$,
by the Young inequality we have that, for every $\delta_1,\delta_2>0$,
\begin{align*}
  &\int_{Q_t}\big[\mu + P\sigma h(\varphi_1)+ 
  (P\sigma_2-A)(h(\varphi_{1})-h(\varphi_2))\big]\mathcal R^{-1}(\eps\mu+\varphi) \\
  &\leq\delta_1\eps\int_{Q_t}|\mu|^2 + \delta_2\int_{Q_t}|\varphi|^2 
  + \frac{P^2\norm{h}_{L^\infty(\erre)}^2}2\int_{Q_t}|\sigma|^2\\ 
  &\qquad+
  K_0^2\left(\frac{1}{4\delta_1\eps}+\frac12+\frac{(P+A)^2\norm{h'}_{L^\infty(\erre)}^2}{4\delta_2}
  \right)\int_0^t\norm{(\eps\mu+\varphi)(s)}^2_{V^*}\,\d s\,.
\end{align*}
Secondly, analogous computations yield
\[
  \chi\int_{Q_t}\sigma\varphi +\int_{Q_t}[C\sigma_2(h(\varphi_2)-h(\varphi_1))]\sigma
  \leq\delta_2\int_{Q_t}|\varphi|^2 + \frac{\chi^2+C^2\norm{h'}^2_{L^\infty(\erre)}}{2\delta_2}
  \int_{Q_t}|\sigma|^2\,.
\]
Finally, we have that
\begin{align*}
\int_{Q_t}(J*\varphi)\varphi
&\leq\int_0^t\norm{J*\varphi(s)}_V\norm{\varphi(s)}_{V^*}\,\d s\leq
(a^*+b^*)\int_0^t\norm{\varphi(s)}_H\norm{\varphi(s)}_{V^*}\,\d s\\
&\leq\delta_2\int_{Q_t}|\varphi|^2 + 
\frac{(a^*+b^*)^2}{4\delta_2}\int_0^t\norm{\varphi(s)}_{V^*}^2\,\d s\\
&\leq\delta_2\int_{Q_t}|\varphi|^2 + 
\frac{(a^*+b^*)^2}{2\delta_2}\int_0^t\norm{(\eps\mu+\varphi)(s)}_{V^*}^2\,\d s
+\frac{\eps(a^*+b^*)^2K_0^2}{2\delta_2}\left(\eps\int_{Q_t}|\mu|^2\right)\,.
\end{align*}
Rearranging the terms we deduce that 
\begin{align*}
  &\frac12\norm{(\eps\mu+\varphi)(t)}_{V^*}^2
  +\eps\int_{Q_t}|\mu|^2 
  +\frac\tau2\norm{\varphi(t)}_H^2 + C_0\int_{Q_t}|\varphi|^2
  +\frac12\norm{\sigma(t)}_H^2 + \int_{Q_t}|\nabla\sigma|^2\\
  &\leq  \frac12\norm{\eps\mu_0+\varphi_0}_{V^*}^2
  +\frac\tau2\norm{\varphi_0}_H^2
  +\frac12\norm{\sigma_0}_H^2 
  +M_{\delta_1,\delta_2,\eps}\int_0^t\left(\norm{\sigma(s)}_H^2
  +\norm{(\eps\mu+\varphi)(s)}_{V^*}^2\right)\,\d s\\
  &\qquad+\left(\delta_1+\frac{\eps(a^*+b^*)^2K_0^2}{2\delta_2}\right)
  \eps\int_{Q_t}|\mu|^2+3\delta_2\int_{Q_t}|\varphi|^2\\
\end{align*}
for some positive constant $M_{\delta_1,\delta_2,\eps}$
depending on the data of the problem and $\eps$, but independent of $\tau$.
Now, it clear that the last two terms
on the right-hand side can be incorporated in the corresponding ones on the left
provided to choose and fix $\delta_1,\delta_2>0$ such that 
\[
\delta_1+\frac{\eps(a^*+b^*)^2K_0^2}{2\delta_2}<1\,, \qquad
3\delta_2<C_0\,.
\]
An elementary computation shows that this is possible if and only if
\[
  \frac{\eps(a^*+b^*)^2K_0^2}2<\frac{C_0}3\,,
\]
which is indeed guaranteed since
$\eps<\eps_0$ and by the smallness assumption on $\eps_0$.
The thesis follows then by the Gronwall lemma.
}

\subsection{Further regularity}
\label{ssec:reg}
This section is devoted to the proof of Theorem~\ref{thm3}, concerning 
regularity of weak solutions, when $\eps,\tau>0$.
To begin with, we improve the regularity of $\varphi$ and $\sigma$
by showing that 
the approximate solutions $(\varphi_\lambda, \mu_\lambda, \sigma_\lambda)$
to the system \eqref{eq1_lam}--\eqref{eq5_lam} satisfy further estimates uniformly in $\lambda$.
We proceed formally, to avoid a further regularization on the system based on
time discretizations.
First, we analyse the system \eqref{eq1_lam}--\eqref{eq5_lam} at the initial time $t=0$
and let us claim that there exists a unique pair 
$(\varphi_{0,\lambda}', \mu_{0,\lambda}', \sigma_{0,\lambda}')
\in H\times V^*\times V^*$ such that, in $\Omega$,
\[
  \begin{cases}
  \eps\mu_{0,\lambda}' + \varphi_{0,\lambda}' - \Delta \mu_0 = (P\sigma_0 - A)h(\varphi_0)\,,\\
  \mu_0 = \tau\varphi_{0,\lambda}' + a\varphi_0 - J*\varphi_0 + F_\lambda'(\varphi_0) - \chi\sigma_0,\\
  \sigma_{0,\lambda}' - \Delta \sigma_0 + B(\sigma_0-\sigma_S(0)) + C\sigma_0h(\varphi_0)
  =-\chia\Delta\varphi_0\,.
  \end{cases}
\]
Indeed, the existence and uniqueness of $\sigma_{0,\lambda}'$ is given by the third equation
and the assumptions \eqref{ip_init}, \eqref{ip_init_reg} and \eqref{ip_data_reg}. It follows directly then 
from the second equation the unique definition for $\varphi_{0,\lambda}'$, 
and finally from the first equation the one of $\mu_{0,\lambda}'$.
Furthermore, from the second equation and 
assumption \eqref{ip_init_reg} it follows that 
$(\varphi_{0,\lambda}')_\lambda$ is uniformly bounded in $H$,
which in turn yields that $(\mu_{0,\lambda}')_\lambda$ is uniformly bounded in $V^*$.

Bearing this in mind, we test \eqref{eq1_lam} by $\partial_t\mu_\lambda$, 
the time-derivative of \eqref{eq2_lam}
by $-\partial_t \varphi_\lambda$, \eqref{eq3_lam} by 
$\partial_t (\sigma_\lambda-\chia\varphi_\lambda)$, 
and take the sum: after integrating in time we obtain
\begin{align}
  &\notag
  \eps\int_{Q_t}\!|\partial_t \mu_\lambda|^2 
  + \frac12 \norma{\nabla\mu_\lambda(t)}_H^2
  +\frac\tau2 \norma{\partial_t\varphi_\lambda(t)}_H^2 
  + \int_{Q_t}\!(a+F_\lambda''(\varphi_\lambda))|\partial_t\varphi_\lambda|^2
  +\int_{Q_t}\!|\partial_t\sigma_\lambda|^2 
  + \frac12 \norma{\nabla(\sigma_\lambda-\chia\varphi_\lambda)(t)}_H^2\\
  &\notag=\frac12\norma{\nabla\mu_0}_H^2
  + \frac\tau2 \norma{\varphi_{0,\lambda}'}_H^2
  +\frac12\norma{\nabla(\sigma_0-\chia\varphi_0)}_H^2 
  +\int_{Q_t}\!(P\sigma_\lambda-A)h(\varphi_\lambda)\partial_t\mu_\lambda\\
  &\qquad+\int_{Q_t}\!\big(J*(\partial_t\varphi_\lambda) 
  + (\chia + \chi)\partial_t\sigma_\lambda\big)\partial_t\varphi_\lambda
  +\int_{Q_t} \!\big(B(\sigma_S-\sigma_\lambda)-Ch(\varphi_\lambda)\sigma_\lambda \big)
  (\partial_t\sigma_\lambda-\chia\partial_t\varphi_\lambda)\,.
  \label{est_regularity}
\end{align}
Now, the second term on the right-hand side is uniformly bounded in $\lambda$
thanks to the remarks above, and so is the first one by assumption. Hence, recalling again {\bf A5}
we infer that 
\begin{align*}
  &\eps\int_{Q_t}|\partial_t \mu_\lambda|^2 
  + \frac12 \norma{\nabla\mu_\lambda(t)}_H^2
  +\frac\tau2 \norma{\partial_t\varphi_\lambda(t)}_H^2
  + C_0 \int_{Q_t} |\partial_t \ph|^2
  +\int_{Q_t}|\partial_t\sigma_\lambda|^2 
  + \frac12 \norma{\nabla(\sigma_\lambda-\chia\varphi_\lambda)(t)}_H^2\\
  &\leq M + \frac\eps2\int_{Q_t}|\partial_t\mu_\lambda|^2 +
  \frac1{2\eps}\int_{Q_t}|(P\sigma_\lambda-A)h(\varphi)|^2
  +\frac12\int_{Q_t}|\partial_t\sigma_\lambda|^2
  \\
  &\qquad
  +\left(a^* + (\chia + \chi)^2+\frac{\chia^2}2\right)\int_{Q_t}|\partial_t\varphi_\lambda|^2
  +\frac32\int_{Q_t}|B(\sigma_S-\sigma_\lambda)-Ch(\varphi_\lambda)\sigma_\lambda|^2\,.
\end{align*}
Taking the estimate \eqref{est1_lam} into account and using the
boundedness of $h$ and $\sigma_S$ we infer that 
\[
  \norm{\varphi_\lambda}_{W^{1,\infty}(0,T; H)}^2 + 
  \norm{\mu_\lambda}_{H^1(0,T; H)\cap L^\infty(0,T; V)}^2 +
  \norm{\sigma_\lambda}_{H^1(0,T; H)}^2 + 
  \norm{\sigma_\lambda-\chia\varphi_\lambda}_{L^\infty(0,T; V)}^2
  \leq M
\]
for some $M>0$ independent of $\lambda$.
As we already know that $(\varphi_\lambda)_\lambda$ 
is uniformly bounded in $L^\infty(0,T; V)$ by 
\eqref{est1_lam}, it is now a standard
matter to pass to the limit as $\lambda \to 0$: recalling
\eqref{phi}--\eqref{mu} and using a comparison argument
for the linear combination $\sigma - \chia \ph$, we have
\begin{align*}
 &  \varphi\in W^{1,\infty}(0,T; H)\cap L^\infty(0,T; V)\,, \qquad
  \mu,  \sigma-\chia \ph\in H^1(0,T; H)\cap L^\infty(0,T; V), \qquad 
  \\
  &\sigma\in H^1(0,T; H)\cap \L2 V.
\end{align*}
Moreover, note that \eqref{eq1} and \eqref{eq3} can be rewritten as
\begin{align}
\label{parab_mu}
 \eps\partial_t\mu- \Delta\mu &= f_\mu:= (P\sigma-A)h(\varphi) - \partial_t\varphi
   &&\hbox{in } Q\,, \\
   \label{parab_sigma}
  \partial_t ( \sigma -\chia \ph)- \Delta( \sigma -\chia \ph) 
  &= f_\sigma: = -B(\sigma-\sigma_S) -  C\sigma h(\varphi)
  - \chia \partial_t \ph
  &&\hbox{in } Q\,,    
\end{align}
endowed with homogeneus Neumann boundary conditions and 
initial data $\mu_{0}, \sigma_0-\chia\varphi_0\in V$.
Since the forcing terms and the initial data satisfy
$f_\mu, f_\sigma\in L^2(0,T; H)$,
the classical parabolic
regularity theory yields 
\begin{align*}
 \mu, \sigma-\chia\varphi\in L^2(0,T; W)\,,
\end{align*}
completing the proof of Theorem~\ref{thm3}.

\subsection{Strong solutions and separation principle}
\label{ssec:strong}
We focus here on the proof of Theorem~\ref{thm4} concerning existence of
strong solutions, separation property, and magnitude regularity,
still in the case $\eps,\tau>0$.
Let us stress that the separation result will allow us to 
exploit the regularity of the linear combination $\sigma- \chia \ph$
to derive further regularity for $\ph$ and $\sigma$.

First of all, by virtue of Theorem~\ref{thm3} we realize that \eqref{parab_mu}
consists of a parabolic equation in the variable $\mu$ with source term
$f_\mu\in\L\infty H$, and with initial datum 
$\mu_0 \in \Lx\infty$ by \eqref{ip_init_sep}.
Therefore, an application of \cite[Thm.~7.1, p.~181]{ladyz}
yields that 
\[
	\mu \in  L^\infty(Q).
\]
In a similar fashion, we notice that in \eqref{parab_sigma}
we have initial datum $\sigma_0-\chia\varphi_0\in V\cap L^\infty(\Omega)$ and
forcing term $f_\sigma\in L^\infty(0,T; H)$ by virtue of Theorem~\ref{thm3}. 
Hence, an application of \cite[Thm.~7.1, p.~181]{ladyz} yields again
\[
  \sigma-\chia\varphi \in L^\infty(Q)\,.
\]
Furthermore, we claim that from assumption {\bf A7} we can deduce 
further regularity also for the term $J*\varphi$. Indeed, 
every kernel verifying Definition~\ref{def_adm} 
satisfy the following result, whose
proof can be found, e.g., in \cite[Lemma~2]{bed-rod-bert}.
\begin{lem}
	\label{LEM_controlofdiv}
	Assume that the kernel $J$ is admissible in the
	sense of the Definition~\ref{def_adm}.
	Then, for every $p\in(1,\infty)$, there exists a positive constant $C_p$
	such that
	\begin{align}
		\label{controlofdiv}
		\norm{\nabla(\nabla J* \psi)}_{L^p(\Omega)^{3 \times 3}}
		\leq
		C_p \|\psi\|_{L^p(\Omega)}
		\quad\forall\,\psi \in {L^p(\Omega)}\,.
	\end{align}
\end{lem}
As a consequence, by taking $p=2$ in \eqref{controlofdiv}, we deduce that 
\[
  \norm{J*\varphi}_{L^\infty(0,T; H^2(\Omega))}\leq C_2\norm{\varphi}_{L^\infty(0,T;H)}\,,
\]
which readily implies, thanks to the continuous inclusion $H^2(\Omega)\embed L^\infty(\Omega)$,
that
\[
  J*\varphi \in L^\infty(0,T; H^2(\Omega))\cap L^\infty(Q) \,.
\]

We are now ready to prove the separation property. To this end, 
note that, taking these remarks into account, under the assumption 
{\bf A6} on $F$, we can rewrite equation \eqref{def_mu_xi} as
\beq
\label{eq_mu_aux}
  \tau\partial_t\varphi + a\varphi + F'(\varphi) - \chi\chia\varphi = f_\varphi:=
  \mu +\chi(\sigma - \chia\varphi) + J*\varphi\,.
\eeq
Besides, we have already proved that $f_\varphi\in L^2(0,T; H^2(\Omega))\cap L^\infty(Q)$,
so that there exists a constant $\overline M>0$ such that 
\[
  \norm{f_\varphi}_{L^\infty(Q)}\leq \overline M\,.
\]
Next, by {\bf A6} and \eqref{ip_init_sep} we infer the existence of $r^*\in(r_0,\ell)$ such that 
\[
  F'(r)-\chi\chia r \geq \overline M \quad\forall\,r\in(r^*,\ell)\,, \qquad
  F'(r)-\chi\chia r \leq -\overline M \quad\forall\,r\in(-\ell,-r^*)\,.
\]
We claim that this choice entails $\varphi(t)\leq r^*$ almost everywhere in $\Omega$, 
for all $t\in[0,T]$.
In fact, by testing \eqref{eq_mu_aux} by $(\varphi-r^*)_+$
and integrating on $[0,t]$, we immediately infer that
\begin{align*}
  &\frac\tau2 \norma{(\varphi(t)-r^*)_+}_H^2 + \int_{Q_t}a\varphi(\varphi-r^*)_+
  =
  \frac\tau2 \norma{(\varphi_0-r^*)_+}_H^2 + 
  \int_{Q_t}\left[f_\varphi
  - (F'(\varphi) - \chi\chia\varphi ) \right](\varphi-r^*)_+\,.
\end{align*}
Now,
since $r^*\in(r_0,\ell)$ and $\norm{\varphi_0}_{L^\infty(\Omega)}\leq r_0$, the first 
term on the right-hand side vanishes. 
Moreover, by definition of $\overline M$ and $r^*$ we have that 
\[
  \int_{Q_t}\left[f_\varphi - (F'(\varphi) - \chi\chia\varphi) \right](\varphi-r^*)_+
  =\int_{Q_t\cap\{\varphi>r^*\}}
  \left[f_\varphi - (F'(\varphi) - \chi\chia\varphi)\right](\varphi-r^*) \leq 0\,.
\]
Recalling also {\bf A5}, we infer that, for every $t\in[0,T]$,
\[
  \frac\tau2 \norma{(\varphi(t)-r^*)_+}_H^2+ 
  a_*\int_{Q_t\cap\{\varphi>r^*\}}\varphi(\varphi-r^*) \leq 0.
\]
Hence, since the second term on the left-hand side is non-negative, we deduce that 
\[
  (\varphi(t)-r^*)_+=0 \quad\forall\,t\in[0,T], \qquad \text{i.e.}\qquad
  \varphi(t,x)\leq r^* \quad\text{for a.e.~}x\in\Omega \quad\forall\,t\in[0,T]\,,
\]
as required. 
The other inequality $\varphi\geq-r^*$ can be deduced analogously by testing
by $-(\varphi+r^*)_-$ instead. Thus, we have shown that 
\[
  \sup_{t\in[0,T]}\norm{\varphi(t)}_{L^\infty(\Omega)}\leq r^*\,,
  \qquad \hbox{with } r^*\in(r_0,\ell)\,.
\]

Let us now
show the $L^2(0,T; W)$-regularity for $\sigma$ and $\chia\varphi$.
To this end, \luca{for an exponent $p>1$ yet to be chosen}, we test the gradient of 
\eqref{eq_mu_aux} by $|\nabla\varphi|^{p-2}\nabla\varphi$ and integrate over $Q_t$ to obtain,
by assumption {\bf A5} and the H\"older and generalized Young inequalities, that 
\begin{align*}
  &\frac\tau{p}\sup_{s\in[0,t]}\norm{\nabla\varphi(s)}_{L^p(\Omega)}^p + 
  C_0\int_{Q_t}|\nabla\varphi|^p 
  \\ & = 
  \frac\tau{p}\norm{\nabla\varphi_0}_{L^p(\Omega)}^p+
  \chi\chia\int_{Q_t}|\nabla\varphi|^p
  - \intQt (\nabla a)\ph \, |\nabla\varphi|^{p-2}\nabla\varphi 
  + \int_{Q_t}\nabla f_\varphi \cdot |\nabla\varphi|^{p-2}\nabla\varphi  \\
  &\leq\frac\tau{p}\norm{\nabla\varphi_0}_{L^p(\Omega)}^p
  +\chi\chia\int_{Q_t}|\nabla\varphi|^p+
  \frac{\tau}{2p}\sup_{s\in[0,t]}\norm{\nabla\varphi(s)}_{L^p(\Omega)}^p
  \\& \quad  
  + \frac{[4 (p-1)]^{p-1}(b^*)^p}{p\tau^{p-1}} \norm{\varphi}_{L^1(0,T; L^p(\Omega))}^p
  +\frac{[4(p-1)]^{p-1}}{p\tau^{p-1}}
  \norm{\nabla f_\varphi}_{L^1(0,T; L^p(\Omega))}^p\,.
\end{align*}
Owing to the already proved regularities $f_\varphi\in L^2(0,T; H^2(\Omega))$ and $\ph \in \L2 V$,
we deduce in particular that $\nabla f_\varphi\in L^2(0,T; V)$ so that,
using the embedding $V\embed L^6(\Omega)$, 
also $\nabla f_\varphi, \ph \in L^2(0,T; L^6(\Omega))$.
Moreover, $\varphi_0\in H^2(\Omega)$ also entails that
$\nabla\varphi_0\in L^6(\Omega)$. Choosing then $p=6$ and using the Gronwall lemma 
yields
\beq
  \label{Lp_grad}
  \varphi\in L^\infty(0,T; W^{1,6}(\Omega))\,.
\eeq
Now, for brevity we proceed formally: a rigorous argument can be 
reproduced on suitable approximations.
Applying the second-order differential operator $\partial_{x_i x_j}$
($i,j=1,2,3$) to equation \eqref{eq_mu_aux}, 
testing it by $\partial_{x_ix_j}\varphi$, and integrating on $[0,t]$
lead to
\begin{align*}
  &\frac\tau2\norm{\partial_{x_ix_j}\varphi(t)}_H^2
  +\int_{Q_t}\left(a+F''(\varphi)\right)|\partial_{x_ix_j}\varphi|^2
  =
  \frac\tau2\norm{\partial_{x_ix_j}\varphi_0}_H^2
  +\int_{Q_t}\partial_{x_ix_j} f_{\varphi}\partial_{x_ix_j}\varphi \\
  &\qquad+ \chi\chia\int_{Q_t}|\partial_{x_ix_j}\varphi|^2
  -\int_{Q_t}\left[\partial_{x_i}a\partial_{x_j}\varphi
  +\partial_{x_j}a\partial_{x_i}\varphi + (\partial_{x_ix_j} a)\varphi
  +F'''(\varphi)\partial_{x_i}\varphi\partial_{x_j}\varphi\right]\partial_{x_ix_j}\varphi\,.
\end{align*}
Now, due to the already proved separation property
$\norm{\varphi}_{L^\infty(Q)}\leq r^*<\ell$, and
recalling that $F\in C^3(-\ell,\ell)$ by {\bf A6}, we have that $F'''(\varphi)\in L^\infty(Q)$.
Hence, exploiting {\bf A5}, 
using the Young inequality, and summing on $i,j=1,2,3$ we deduce,
recalling that $\varphi\in L^2(0,T; V)$, that 
\begin{align*}
  &\frac\tau2\norm{\varphi(t)}_{H^2(\Omega)}^2
  +C_0\int_0^t\norm{\varphi(s)}_{H^2(\Omega)}^2\,\d s
  \leq\frac\tau2\norm{\varphi_0}_{H^2(\Omega)}^2
  +(2+\chi\chia)\int_0^t\norm{\varphi(s)}_{H^2(\Omega)}^2\,\d s\\
  & \qquad+\frac12\norm{f_{\varphi}}^2_{L^2(0,T; H^2(\Omega))}
  + 2\int_{Q_t}|\nabla a|^2|\nabla\varphi|^2
  +\frac12 \int_{Q_t} \sum_{i,j=1}^3|\partial_{x_ix_j}a|^2|\varphi|^2
  +\frac12\norm{F'''(\varphi)}_{L^\infty(Q)}^2
  \int_{Q_t}|\nabla\varphi|^4\,.
\end{align*}
Moreover, $\norm{\nabla a}_{L^\infty(\Omega)}\leq b^*$ by {\bf A5},
$\norm{a}_{W^{2,p}(\Omega)}\leq C_p$ for all $p\in(1,+\infty)$ by \eqref{controlofdiv}
and $\varphi\in L^\infty(0,T; V)$, so that the H\"older inequality yields
\[
  \int_{Q_t}|\nabla a|^2|\nabla\varphi|^2\leq (b^*)^2\norm{\varphi}_{L^2(0,T; V)}^2 \leq M
\]
and, by the continuous embedding $V\embed L^4(\Omega)$, also that
\[
  \int_{Q_t}\sum_{i,j=1}^3|\partial_{x_ix_j}a|^2|\varphi|^2\leq \norm{a}_{W^{2,4}(\Omega)}^2
  \norm{\varphi}_{L^4(0,T; L^4(\Omega))}^2\leq M'\norm{\varphi}_{L^4(0,T; V)}^2\leq M
\]
for certain constants $M, M'>0$. 
Using then \eqref{Lp_grad}, we are left with 
\begin{align*}
  &\frac\tau2\norm{\varphi(t)}_{H^2(\Omega)}^2
  +C_0\int_0^t\norm{\varphi(s)}_{H^2(\Omega)}^2\,\d s\leq\frac\tau2\norm{\varphi_0}_{H^2(\Omega)}^2
  +M\left(1+\int_0^t\norm{\varphi(s)}_{H^2(\Omega)}^2\,\d s\right)
\end{align*}
so that a Gronwall argument produces 
\[
  \varphi\in L^\infty(0,T; H^2(\Omega))\,.
\]
At this point, the equation for $\sigma$ can be written also as
\[
  \partial_t\sigma - \Delta\sigma = \tilde f_\sigma := -B(\sigma-\sigma_S) - C\sigma h(\varphi) - \chia\Delta\varphi \in L^\infty(0,T; H)\,,
\]
with initial datum $\sigma_0\in V\cap L^\infty(\Omega)$. Hence, by parabolic regularity
theory and again \cite[Thm.~7.1]{ladyz}, we deduce that 
\[
	\sigma \in \H1 H \cap \L\infty V \cap \L2 W \cap L^\infty(Q)\,.
\]
Since we already know that $\sigma-\chia\varphi\in L^2(0,T; W)$, 
by comparison we also infer 
\[
  \chia\varphi\in L^2(0,T; W)\,.
\]
To conclude, we go back to equation \eqref{eq_mu_aux} and note that, by difference, 
we have also the regularity 
\[
  \partial_t\varphi \in L^\infty(0,T; V)\cap L^2(0,T; H^2(\Omega))\cap L^\infty(Q)
\]
which completes the proof of Theorem~\ref{thm4}.

\subsection{Refined continuous dependence}
We prove here the refined stability estimates contained in Theorem~\ref{thm5}
which is now possible in light of the strong
regularity result established by Theorem \ref{thm4}.
It is worth pointing out that both the chemotaxis and active transport
mechanisms are now included in the analysis.
Employing the same notation of Subsection~\ref{ssec:cont}, we consider the system \cdsys\ and 
test \eqref{cd_1} by $\partial_t\mu$, the time-derivative 
of \eqref{cd_2} by $-\partial_t\varphi$, \eqref{cd_3} by $\partial_t(\sigma-\chia\varphi)$,
and integrate over $[0,t]$, to obtain
\begin{align*}
  &\eps\int_{Q_t}|\partial_t\mu|^2 
  + \frac12 \norma{\nabla\mu(t)}_H^2
  + \frac\tau2 \norma{\partial_t\varphi(t)}_H^2
  + \int_{Q_t}(a + F''(\varphi_1))|\partial_t\varphi|^2
  \luca{+}\int_{Q_t}|\partial_t\sigma|^2 + 
  \frac12 \norma{\nabla(\sigma-\chia\varphi)(t)}_H^2\\
  &= \frac12 \norma{\nabla\mu_0}_H^2 +
  \frac\tau2 \norma{\varphi_0'}^2_H
  + \frac12 \norma{\nabla(\sigma_0-\chia\varphi_0)}_H^2
  +\int_{Q_t}\left[P\sigma h(\varphi_1) 
  + (P\sigma_2-A)(h(\varphi_{1})-h(\varphi_2))\right]\partial_t\mu\\
  &\qquad+\int_{Q_t}\left[(F''(\varphi_2)-F''(\varphi_1))\partial_t\varphi_2+
  \chi\partial_t\sigma + J*\partial_t\varphi\right]\partial_t\varphi
  +\chia\int_{Q_t}\partial_t\sigma\partial_t\varphi\\
  &\qquad+\int_{Q_t}\left[C\sigma_2(h(\varphi_2)-h(\varphi_1)) - C\sigma h(\varphi_1) - B\sigma
  \right](\partial_t\sigma-\chia\partial_t\varphi)\,.
\end{align*}
First of all, notice that $\varphi_0'$ is such that 
\[
\mu_0=\tau\varphi_0' + a\varphi_0 - J*\varphi_0+F'(\varphi_0^1)-F'(\varphi_0^2) - \chi\sigma_0\,.
\]
Since the initial data satisfy \eqref{ip_init}, \eqref{ip_init_reg}, and \eqref{ip_init_sep},
for $i=1,2$ we have that $\varphi_0'\in V\cap L^\infty(\Omega)$.
Now, recalling that $F\in C^3([-r_0,r_0])$, we have
\[
  \norm{\varphi_0'}_H\leq\frac1\tau\left(\norm{\mu_0}_H+ 2a^*\norm{\varphi_0}_H
  +\norm{F''}_{C^0([-r_0,r_0])}\norm{\varphi_0}_H + \chi\norm{\sigma_0}_H\right)\,.
\]
Secondly, by the separation
property for $\varphi_1$ and $\varphi_2$,
we have $\norm{\varphi_i}_{L^\infty(Q)}\leq r^*<\ell$ for $i=1,2$
and combined with $F\in C^3([-r^*, r^*])$ we have $F''\in W^{1,\infty}(-r^*,r^*)$, so that
\[
  |F''(\varphi_1)-F''(\varphi_2)|\leq \norm{F'''}_{C^0([-r^*,r^*])}|\varphi_1-\varphi_2| \quad\text{a.e.~in } Q\,.
\]
Taking this information into account, using {\bf A5},
and exploiting the regularities $h\in W^{1,\infty}(\erre)$,
$\sigma_2\in L^\infty(Q)$, and $\partial_t\varphi_2\in L^\infty(Q)$,
we invoke the Young inequality to infer
\begin{align}
  \nonumber
  &\int_{Q_t}|\partial_t\mu|^2 
  + \norma{\nabla\mu(t)}_H^2 +
  \norma{\partial_t\varphi(t)}_H^2
  +\int_{Q_t}|\partial_t\varphi|^2+
  \int_{Q_t}|\partial_t\sigma|^2 
  + \norma{\nabla(\sigma-\chia\varphi)(t)}_H^2\\
  \label{stab_aux}
  & \leq M\left(
  \norm{\mu_0}_V^2 +\norm{\varphi_0}_H^2+\norm{\sigma_0}_H^2
  + \norm{\nabla(\sigma_0-\chia\varphi_0)}_H^2
  +\int_{Q_t}(|\sigma|^2 + |\varphi|^2 + |\partial_t\varphi|^2)\right)\,,
\end{align}
where the constant $M>0$ may depend on $\eps, \tau$ and on structural data.
Now, we take the gradient of \eqref{cd_2} and test it by $\nabla\varphi$,
getting
\begin{align*}
  \frac\tau2 \norma{\nabla\varphi(t)}_H^2
  + \int_{Q_t}(a + F''(\varphi_1))|\nabla\varphi|^2
  &=\frac\tau2 \norma{\nabla\varphi_0}_H^2
  +\int_{Q_t}\left(F''(\varphi_2)-F''(\varphi_1)\right)\nabla\varphi_2\cdot\nabla\varphi\\
  &\quad +
  \int_{Q_t}\left(\nabla\mu+\chi\nabla\sigma + (\nabla J)*\varphi - (\nabla a)\varphi\right)\cdot\nabla\varphi\,.
\end{align*}
Using {\bf A5}, along with the Lipschitz continuity of $F''$ on $[-r^*,r^*]$,
and the identity 
\begin{align*}
\chi\nabla \sigma \cdot \nabla \ph
	=\chi(\nabla (\sigma -\chia \ph) + \chia \nabla \ph)\cdot \nabla \ph,
\end{align*}
and the Young inequality lead to 
\begin{align*}
  \frac\tau2 \norma{\nabla\varphi(t)}_H^2
  + C_0\int_{Q_t}|\nabla\varphi|^2
  &\leq\frac\tau2 \norma{\nabla\varphi_0}_H^2
  +\norm{F'''}_{C^0([-r^*,r^*])}\int_{Q_t}|\varphi||\nabla\varphi_2||\nabla\varphi|\\
  &+\int_{Q_t}|\nabla\mu|^2 + 
  \chi^2\int_{Q_t}|\nabla(\sigma-\chia\varphi)|^2
  +(1+\chi\chia)\int_{Q_t}|\nabla\varphi|^2 +2(b^*)^2\int_{Q_t}|\varphi|^2\,.
\end{align*}
From the embedding $V\embed L^4(\Omega)$, H\"older's inequality
and the regularity $\varphi_2\in L^\infty(0,T; H^2(\Omega))$, we find
\[
  \int_{Q_t}|\varphi||\nabla\varphi_2||\nabla\varphi|\leq
  M'\int_0^t\norm{\varphi(s)}_{V}\norm{\varphi_2(s)}_{H^2(\Omega)}\norm{\nabla\varphi(s)}_H\,\d s
  \leq M\int_0^t\norm{\varphi(s)}_V^2\,\d s
\]
for some constants $M,M'>0$. We deduce then that, possibly updating $M$, 
for every $t\in[0,T]$,
\beq
  \label{stab_aux2}
  \norma{\nabla\varphi(t)}_H^2
  \leq
  M\left(\norma{\nabla\varphi_0}_H^2+\int_{Q_t}|\nabla\mu|^2
  +\int_{Q_t}|\nabla(\sigma-\chia\varphi)|^2
  +\int_0^t\norm{\varphi(s)}_V^2\,\d s\right)\,.
\eeq
\luca{We now combine the estimates \eqref{stab_aux} and \eqref{stab_aux2} to infer that, for all $t\in[0,T]$, }
\begin{align*}
  &\int_{Q_t}|\partial_t\mu|^2 
  + \norma{\nabla\mu(t)}_H^2 +
 \norma{\partial_t\varphi(t)}_H^2
 + \norma{\nabla\varphi(t)}_H^2
  + \int_{Q_t}|\partial_t\sigma|^2 
 + \norma{\nabla\sigma(t)}_H^2\\
  &\leq M\left(
  \norm{\mu_0}_V^2 +\norm{\varphi_0}_V^2+\norm{\sigma_0}_V^2
  +\int_0^t\left(\norm{\nabla\mu(s)}_H^2+
  \norm{\sigma(s)}_V^2+\norm{\varphi(s)}_V^2+
  \norm{\partial_t\varphi(s)}_H^2\right)\,\d s\right)\,.
\end{align*}
Since the quantities
$\norm{\sigma_2}_{L^\infty(Q)}$, $\norm{\partial_t\varphi_2}_{L^\infty(Q)}$,
and $\norm{\varphi_2}_{L^\infty(0,T; H^2(\Omega))}$ appearing implicitly in the constant $M$
can be in turn handled in terms on the norms of the initial data appearing 
in \eqref{ip_init}, \eqref{ip_init_reg}, and \eqref{ip_init_sep}, 
we can close the estimate by the Gronwall lemma.
Moreover, comparison in equation \eqref{cd_1} produces
\[
  \norm{\Delta\mu}_{L^2(0,T; H)}\leq M\left(\norm{\varphi}_{H^1(0,T; H)}
  +\norm{\partial_t\mu}_{L^2(0,T; H)}+\norm{\sigma}_{L^2(0,T; H)}\right)\,,
\]
where all the terms on the right-hand side have already been estimated.
Similarly, from \eqref{cd_2} we get 
\[
  \norm{\partial_t\varphi}_{L^\infty(0,T; V)}\leq M
  \left(\norm{\mu}_{L^\infty(0,T; V)} + \norm{\varphi}_{L^\infty(0,T; V)}
  +\norm{\sigma}_{L^\infty(0,T; V)}\right)\,,
\]
while from \eqref{cd_3} we get 
\[
  \norm{\Delta(\sigma-\chia\varphi)}_{L^2(0,T; H)}\leq 
  M\left(\norm{\sigma}_{H^1(0,T; H)} + \norm{\varphi}_{L^2(0,T; H)}\right)\,.
\]
Collecting the above estimates, along with elliptic regularity theory, we deduce that 
\begin{align}
  \nonumber
  &\norm{\mu}_{H^1(0,T; H)\cap L^\infty(0,T; V)\cap L^2(0,T; W)}^2+
  \norm{\varphi}^2_{W^{1,\infty}(0,T; V)}
  +\norm{\sigma}^2_{H^1(0,T; H)\cap L^\infty(0,T; V)}\\
  \label{stab_aux3}
  &\qquad+\norm{\sigma-\chia\varphi}_{\H1 H \cap \L\infty V \cap L^2(0,T; W)}^2
  \leq M\left(\norm{\mu_0}_V^2 + \norm{\varphi_0}_V^2
  +\norm{\sigma_0}_V^2 \right)\,.
\end{align}

To complete the proof, we need to show
a stability estimate for $\partial_t \varphi$ and $\sigma$ also in
$L^2(0,T; H^2(\Omega))$
and $L^2(0,T; W)$, respectively.
In this direction, for any $i,j=1,2,3$, we apply the 
differential operator $\partial_{x_ix_j}$ to \eqref{cd_2}
and test the obtained equation by $\partial_{x_ix_j}\varphi$, getting
\begin{align*}
  &\frac\tau2\norm{\partial_{x_ix_j}\varphi(t)}_H^2 
  + \int_{Q_t}(a+F''(\varphi_1))|\partial_{x_ix_j}\varphi|^2
  =\frac\tau2\norm{\partial_{x_ix_j}\varphi_0}_H^2+
  \int_{Q_t}\partial_{x_ix_j}(\mu + \chi(\sigma-\chia\varphi) + J*\varphi)\partial_{x_ix_j}\varphi\\
  &\qquad+\chi\chia\int_{Q_t}|\partial_{x_ix_j}\varphi|^2
  -\int_{Q_t}(\partial_{x_i}a\partial_{x_j}\varphi+
  \partial_{x_j}a\partial_{x_i}\varphi + (\partial_{x_ix_j}a)\varphi)\partial_{x_ix_j}\varphi\\
  &\qquad+\int_{Q_t}\left[(F''(\varphi_2)-F''(\varphi_1))\partial_{x_ix_j}\varphi_2
  +(F'''(\varphi_2)-F'''(\varphi_1))\partial_{x_i}\varphi_1
  \partial_{x_j}\varphi_2\right]\partial_{x_ix_j}\varphi\\
  &\qquad-\int_{Q_t}\left[F'''(\varphi_1)\partial_{x_i}\varphi_1\partial_{x_j}\varphi
  +F'''(\varphi_2)\partial_{x_i}\varphi\partial_{x_j}\varphi_2\right]\partial_{x_ix_j}\varphi\,.
\end{align*}
We recall that, due to {\bf A6}, $F\in C^4([-r^*, r^*])$, 
so that $F'''$ is Lipschitz continuous on $[-r^*, r^*]$,
and as a consequence of the separation result, also
$F'''(\varphi_i)\in L^\infty(Q)$, for $i=1,2$.
Now, we use the H\"older and Young inequalities and sum on $i,j=1,2,3$:
proceeding as in Subsection~\ref{ssec:strong} and 
exploiting assumptions {\bf A5} and {\bf A7}, we get 
\begin{align*}
  &\frac\tau2\norm{\varphi(t)}_{H^2(\Omega)}^2 
  + C_0\int_0^t\norm{\varphi(s)}_{H^2(\Omega)}^2\,\d s\\
  &\leq\frac\tau2\norm{\varphi_0}_{H^2(\Omega)}^2+
  M\left(\norm{\mu}_{L^2(0,T; W)}^2 + \norm{\sigma-\chia\varphi}^2_{L^2(0,T; W)}
   +\int_0^t\norm{\varphi(s)}_{H^2(\Omega)}^2\,\d s \right)\\
  &+M\sum_{i,j=1}^3\left(\int_{Q_t}|\varphi|^2\left(|\partial_{x_ix_j}\varphi_2|^2+
  |\partial_{x_i}\varphi_1|^2|\partial_{x_j}\varphi_2|^2\right)
  +\int_{Q_t}\left(|\partial_{x_i}\varphi_1|^2|\partial_{x_j}\varphi|^2
  +|\partial_{x_i}\varphi|^2|\partial_{x_j}\varphi_2|^2\right)\right).
\end{align*}
The first bracket on the right-hand side can be controlled using
\eqref{stab_aux3} and the Gronwall lemma, while
the sum-term can be estimated using the H\"older
inequality and the continuous inclusions $V\embed L^4(\Omega)$ 
and $H^2(\Omega)\embed L^\infty(\Omega)$ by
\begin{align*}
  &\int_0^t\norm{\varphi(s)}^2_{L^\infty(\Omega)}
  \left(\norm{\varphi_2(s)}^2_{H^2(\Omega)} 
  + \norm{\nabla\varphi_1(s)}_V^2\norm{\nabla\varphi_2(s)}_V^2\right)\,\d s\\
  &\qquad+\int_0^t\norm{\nabla\varphi(s)}_{L^4(\Omega)}^2
  \left(\norm{\nabla\varphi_1(s)}_{L^4(\Omega)}^2
  +\norm{\nabla\varphi_2(s)}_{L^4(\Omega)}^2\right)\,\d s\\
  &\leq M'\left(\norm{\varphi_2}_{L^\infty(0,T; H^2(\Omega))}^2
  +\norm{\varphi_1}^2_{L^\infty(0,T; H^2(\Omega))}
  \norm{\varphi_2}^2_{L^\infty(0,T; H^2(\Omega))}\right)
  \int_0^t\norm{\varphi(s)}^2_{H^2(\Omega)}\,\d s\\
  &\qquad+M'\left(\norm{\varphi_1}_{L^\infty(0,T; H^2(\Omega))}^2
  +\norm{\varphi_2}_{L^\infty(0,T; H^2(\Omega))}^2\right)\int_0^t\norm{\varphi(s)}_{H^2(\Omega)}^2\,\d s\,.
\end{align*}
Taking these estimates into account
and recalling the regularity $\varphi_1,\varphi_2\in L^\infty(0,T; H^2(\Omega))$, we conclude that 
\begin{align*}
  &\norm{\varphi(t)}_{H^2(\Omega)}^2
  \leq \norm{\varphi_0}_{H^2(\Omega)}^2+
  M\left(\norm{\mu}_{L^2(0,T; W)}^2 + \norm{\sigma-\chia\varphi}^2_{L^2(0,T; W)}
  +\int_0^t\norm{\varphi(s)}_{H^2(\Omega)}^2\,\d s \right)\,
\end{align*}
so that Gronwall's lemma along with the above estimates produces
\begin{align*}
	\norma{\ph}_{\L\infty {\Hx2}}^2
	\leq M\left(\norm{\mu_0}_V^2 + \norm{\varphi_0}_{\Hx2}^2
 	 +\norm{\sigma_0}_V^2 \right)\,.
\end{align*}

The stability estimate for $\sigma$ in $L^2(0,T; W)$ 
follows by comparison in \eqref{cd_3} and elliptic regularity theory.
Finally, by comparison in equation \eqref{cd_2} we also infer the stability estimate
for $\partial_t\varphi$ in $L^2(0,T; H^2(\Omega))$, 
concluding the proof of Theorem~\ref{thm5}.

\section{Asymptotics as $\eps\searrow0$}
\label{sec:eps}
This section is completely devoted to discuss the asymptotic behavior of
system \sys\ as $ \eps\searrow 0$, when $\tau>0$ is fixed.
Namely, we aim at proving Theorems~\ref{thm6} and \ref{thm7}.
Henceforth, let us assume $\tau$ to be positive and fixed.
Moreover, using the notation introduced by Theorem \ref{thm6}, we 
indicate with $(\ph_\ept, \mu_\ept, \sigma_\ept, \xi_\ept)$ the unique weak solution to
\eqref{eq1}-\eqref{eq5} with $\eps,\tau>0.$

\subsection{Uniform estimates}
\label{ssec:unif_eps}
Proceeding as in Subsection~\ref{ssec:unif_est},
we perform the analogous estimates that we used to deduce \eqref{est_aux4}.
In particular, since the implicit constant $M$ in \eqref{est_aux4} is independent 
of $\eps$ and $\tau$, recalling that we are assuming $\chia=0$, we realize that 
\begin{align}
  \nonumber
  &\frac\eps4\norm{\mu_{\ept}(t)}_H^2 
  + (1+4c_a\eps)\int_{Q_t}|\nabla\mu_{\ept}|^2
  +\frac\tau2\int_{Q_t}|\partial_t\varphi_{\ept}|^2
  +\int_\Omega F(\varphi_\ept(t))
  +\frac12\norm{\sigma_{\ept}(t)}_H^2 \\
  \nonumber
  &\qquad+ \int_{Q_t}|\nabla\sigma_{\ept}|^2
  +c_a\norm{(\eps\mu_{\ept}+\varphi_{\ept})(t)}_H^2+
  2c_a\tau\norm{\nabla\varphi_{\ept}(t)}_H^2 + 
  2c_aC_0\int_{Q_t}|\nabla\varphi_{\ept}|^2\\
  \nonumber
  &\leq \frac32\eps\norm{\mu_{0,\ept}}_H^2 +
  (a^*+4c_a)\norm{\varphi_{0,\ept}}_H^2+ 2c_a\tau\norm{\nabla\varphi_{0,\ept}}_H^2+
   \norm{F(\varphi_{0,\ept})}_{L^1(\Omega)}
   +\frac12\norm{\sigma_{0,\ept}}_H^2\\ 
   \nonumber
  &\qquad
  +M\left(1+\int_{Q_t}|\eps\mu_{\ept}+\varphi_{\ept}|^2
  + \int_{Q_t}|\varphi_{\ept}|^2
   + \int_{Q_t}|\sigma_{\ept}|^2\right)
   +\frac{\chi^2}{2\tau} \int_{Q_t}|\sigma_{\ept}|^2\\
   \label{aux_eps1}
   &\qquad
  +\frac12\int_{Q_t}|\nabla\sigma_{\ept}|^2
  + 8c_a^2\chi^2\int_{Q_t}|\nabla\ph_{\ept}|^2
  +\int_{Q_t}(P\sigma_{\ept}-A)h(\varphi_{\ept})\mu_{\ept}\,.
\end{align}
All the terms referring to the initial data on the right-hand side
are uniformly bounded in $\eps$ by virtue of assumptions \eqref{init2_eps0}--\eqref{init3_eps0}. 
Moreover, all the remaining terms can be handled 
using the Gronwall lemma, except for the last one.
To this end, note that 
by the Poincar\'e-Wirtinger inequality \eqref{poincare}, using the 
fact that $h$ is bounded, and the uniform bound $\norm{\sigma_{\ept}}_{L^\infty(Q)}\leq1$, we have 
\begin{align*}
	\int_{Q_t}(P\sigma_\ept  - A )h(\varphi_{\ept})\mu_{\ept}
 	& \leq 
 	\int_{Q_t}(P\sigma_\ept - A )h(\varphi_\ept)(\mu_\ept- (\mu_\ept)_{\Omega})
 	 + \int_{Q_t}(P\sigma_\ept - A)h(\varphi_\ept)(\mu_\ept)_{\Omega}
 	\\&	\leq   
 	\frac 12 \int_{Q_t} |\nabla \mu_{\ept}|^2 + M
 	 + (P+A)\norm{h}_{L^\infty(\erre)}t^{1/2}\norm{(\mu_\ept)_{\Omega}}_{L^2(0,t)}\,.
\end{align*}
Furthermore, noting that $(a\varphi_\ept - J*\varphi_\ept)_\Omega = 0$,
by comparison in equation \eqref{eq2} we get 
\[
  (\mu_\ept)_\Omega=\tau(\partial_t\varphi_\ept)_\Omega 
  + (\xi_\ept+F_2'(\varphi_\ept))_\Omega -\chi(\sigma_\ept)_\Omega\,,
\]
so that thanks to assumption \eqref{pol_growth} implies that 
\begin{align*}
  \norm{(\mu_\ept)_\Omega}_{L^2(0,t)} &
  \leq \tau\norm{\partial_t\varphi_\ept}_{L^2(Q_t)}
  +\norm{\xi_\ept + F_2'(\varphi_\ept)}_{L^2(0,t; L^1(\Omega))} + 
  \chi\norm{\sigma_{\ept}}_{L^2(0,t; H)}\\
  &\leq M\left(1 + \tau^2\int_{Q_t}|\partial_t\varphi_\ept|^2
  +\sup_{s\in[0,t]}\int_\Omega F(\varphi_\ept(s)) + 
  \sup_{s\in[0,t]}\norm{\sigma_\ept(s)}_H^2\right)\,,
\end{align*}
for a certain constant $M>0$, independent of $\eps$. 
Putting this information together, 
we first choose $t\in[0,T_0]$, where $T_0\in(0,T]$ is 
fixed sufficiently small so that the term corresponding to $t^{1/2}$
can be incorporated on the left-hand side, for example by picking a $T_0$ such that
\[
  (P+A)\norm{h}_{L^\infty(\erre)}T_0^{1/2}<\frac1 {\luca{2 M}}\,.
\] 
We then take supremum in $t\in[0,T_0]$ on the left-hand side of 
the inequality \eqref{aux_eps1} and rearrange the terms:
the estimate can be closed on the time interval $[0,T_0]$ using the Gronwall lemma.
As the choice of $T_0$ is independent of $\eps$, $\tau$,
and of the initial data (it only depends on $A$, $P$, $C_F$, $h$, and $\chi$),
repeating the same argument we can close the estimate also on $[T_0,2T_0]$,
and so on,
so that a classical patching argument guarantees the existence
of a constant $M>0$, independent of $\eps$, such that 
\begin{align}
  \label{est1_eps}
  \norm{\varphi_\ept}_{H^1(0,T; H)\cap L^\infty(0,T; V)}+
  \norm{\sigma_\ept}_{L^\infty(0,T; H)\cap L^2(0,T; V)} &\leq M\,,\\
  \label{est2_eps}
  \norm{(\mu_\ept)_\Omega}_{L^2(0,T)} + \norm{\nabla\mu_\ept}_{L^2(0,T; H)}+
  \eps^{1/2}\norm{\mu_\ept}_{L^\infty(0,T; H)}&\leq M\,.
\end{align}
From estimate \eqref{est2_eps}, the Poincar\'e-Wirtinger inequality yields 
\beq\label{est3_eps}
  \norm{\mu_\ept}_{L^2(0,T; V)}\leq M\,.
\eeq
Lastly, by comparison in \eqref{eq3}, we also deduce that 
\beq
  \label{est4_eps}
   \norm{\sigma_\ept}_{H^1(0,T; V^*)}\leq M\,,
\eeq
while by comparison in \eqref{eq2} we have that 
\beq
  \label{est5_eps}
  \norm{\xi_\ept}_{L^2(0,T; H)}\leq M\,.
\eeq

\subsection{Passage to the limit}
\label{ssec:lim_eps}
From the estimates \eqref{est1_eps}--\eqref{est5_eps}
and classical compactness arguments, we infer the
existence of a quadruplet $(\ph_{\tau}, \mu_{\tau}, \sigma_{\tau}, \xi_{\tau})$ with 
\begin{align*}
  &\varphi_\tau\in H^1(0,T; H)\cap L^\infty(0,T;V)\,, \qquad\mu_\tau\in L^2(0,T; V)\,,\\
  &\sigma_\tau\in H^1(0,T; V^*)\cap L^2(0,T; V)\,, \qquad \xi_\tau\in L^2(0,T; H)\,,
\end{align*}
such that, as $\eps \searrow 0$, along a non-relabelled subsequence, it holds that
the weak, weak$^*$ and strong convergences \epstozero\ and \eqref{epstozero:strong} are fulfilled.
We are then left to show that $(\ph_{\tau}, \mu_\tau, \sigma_{\tau}, \xi_{\tau})$
yields a solution to
\sys\ with $\eps=0$ in the sense of Theorem~\ref{thm6}.
In this direction, let us exploit the strong convergence of the phase variable \eqref{epstozero:strong}
along with the continuity and boundedness of $h$, and
Lebesgue convergence theorem, to deduce that, as $\eps \searrow 0 $,
\begin{align*}
	h(\ph_\ept) \to h(\ph_\tau) \quad\text{in } L^p(Q)\quad\forall\,p\geq1\,, \qquad
	F_2'(\ph_\ept) \to F_2'(\ph_\tau) 
	\quad\hbox{in } C^0([0,T]; H)\,.
\end{align*} 
Moreover, the strong-weak closure of $\partial F_1$ 
(see, e.g., \cite[Cor.~2.4, p.~41]{barbu-monot}) entails that
$\xi_\tau\in\partial F_1(\varphi_\tau)$ almost everywhere in $Q$.
Lastly, it is not difficult to pass to the limit in the weak formulation of \sys\ 
to conclude that $(\mu_{\tau}, \ph_{\tau}, \sigma_{\tau}, \xi_{\tau})$
solves \sys\ with $\eps=0$, as we claimed.
The maximum principle for $\sigma_\tau$ can be then obtained repeating
the argument of Subsection \ref{sub:maxpr} leading to $\sigma_\tau \in L^\infty(Q)$.
This concludes the proof of Theorem~\ref{thm6}.

\subsection{Error estimate}
\label{ssec:error_eps}
We focus here on the error estimate as $\eps\searrow0$ presented by Theorem~\ref{thm7}
under the additional assumptions \eqref{pot_reg}--\eqref{init4_eps0}.

First of all, we need to deduce an additional estimate on $\partial_t\mu_\ept$.
Arguing as in Subsection~\ref{ssec:reg}, by considering \eqref{est_regularity} and
multiplying it by $\eps^{1/2}$ (recall that $\chia=0$), we obtain
\begin{align*}
  &\eps^{3/2}\int_{Q_t}|\partial_t \mu_\ept|^2 
  + \frac{\eps^{1/2}}2 \norma{\nabla\mu_\ept(t)}_H^2
  +\frac{\tau\eps^{1/2}}2  \norma{\partial_t\varphi_\ept(t)}_H^2
  + C_0\eps^{1/2}\int_{Q_t}|\partial_t\varphi_\ept|^2
  \\&\qquad
  +\eps^{1/2}\int_{Q_t}|\partial_t\sigma_\ept|^2 
  + \frac{\eps^{1/2}}2 \norma{\nabla\sigma_\ept(t)}_H^2\\
  &\leq\frac{\eps^{1/2}}2\norma{\nabla\mu_{0,\ept}}_H^2+ 
  \frac{\tau\eps^{1/2}}2\norma{\varphi'_{0,\ept}}_H^2
  +\frac{\eps^{1/2}}2\norma{\nabla\sigma_{0,\ept}}_H^2
  +\eps^{1/2}\int_{Q_t}(P\sigma_\ept-A)h(\varphi_\ept)\partial_t\mu_\ept\\
  &\qquad+\eps^{1/2}\int_{Q_t}(J*(\partial_t\varphi_\ept) 
  + \chi\partial_t\sigma_\ept)\partial_t\varphi_\ept
  +\eps^{1/2}\int_{Q_t}\big (B(\sigma_S-\sigma_\ept)-Ch(\varphi_\ept)\sigma_\ept\big )
  \partial_t\sigma_\ept\,.
\end{align*}
The last two terms on the right-hand side can be easily handled as in 
Subsection~\ref{ssec:reg}, using the averaged Young inequality.
Moreover, since $\varphi'_{0,\ept}$ satisfies
\[
  \mu_{0,\ept}=\tau \varphi'_{0,\ept} + a\varphi_{0,\ept}-J*\varphi_{0,\ept}
  + F'(\varphi_{0,\ept}) - \chi\sigma_{0,\ept}\,,
\]
the first three terms on the right-hand side of the inequality above
are uniformly bounded in $\eps$ thanks
to the assumptions \eqref{init2_eps0}--\eqref{init3_eps0} and \eqref{init4_eps0}.
As for the fourth term, this can be treated using 
integration by parts in time and the boundedness of $\sigma_\ept$ in \eqref{max_sigma} as
\begin{align*}
  &-\eps^{1/2}P\int_{Q_t}\partial_t\sigma_\ept h(\varphi_\ept)\mu_\ept
  -\eps^{1/2}\int_{Q_t}(P\sigma_\ept-A)h'(\varphi_\ept)\partial_t\varphi_\ept\mu_\ept\\
  &\qquad+\eps^{1/2}\int_\Omega(P\sigma_\ept(t)-A)h(\varphi_\ept(t))\mu_\ept(t)
  -\eps^{1/2}\int_\Omega(P\sigma_{0,\ept}-A)h(\varphi_{0,\ept})\mu_{0,\ept}\\
  &\leq \frac{\eps^{1/2}}{4}\int_{Q_t}|\partial_t\sigma_\ept|^2+
  \eps^{1/2}P^2\norm{h}^2_{L^\infty(\erre)}\norm{\mu_\ept}^2_{L^2(0,T; H)}\\
  &\qquad+\eps^{1/2}(P+A)\norm{h}_{W^{1,\infty}(\erre)}
  \left(\norm{\partial_t\varphi_\ept}_{L^2(0,T; H)}
  \norm{\mu_\ept}_{L^2(0,T; H)}+
  2\norm{\mu_\ept}_{C^0([0,T]; H)}\right)\,,
\end{align*}
where the right-hand side is uniformly bounded in $\eps$
thanks to \eqref{est1_eps}--\eqref{est5_eps}.
Putting this information together, we deduce that 
\begin{align}
  \label{error1_eps}
  \eps^{3/4}\norm{\mu_\ept}_{H^1(0,T; H)} 
  + \eps^{1/4}\norm{\mu_\ept}_{L^\infty(0,T; V)}&\leq M\,,\\
  \label{error2_eps}
  \eps^{1/4}\norm{\varphi_\ept}_{W^{1,\infty}(0,T; H)}
  +\eps^{1/4}\norm{\sigma_\ept}_{H^1(0,T; H)\cap L^\infty(0,T; V)} &\leq M\,.
\end{align}

We are now ready to show the error estimate.
Taking the difference between the
unique solution $(\mu_\ept, \ph_\ept, \sigma_\ept, \xi_\ept)$
to \sys\ with $\eps,\tau >0$ and $\chia=0$
and the solution $(\mu_\tau, \ph_\tau, \sigma_\tau, \xi_\tau)$ to
\sys\ with $\eps=\chia=0$ obtained 
in Subsection~\ref{ssec:lim_eps} leads us to
\begin{align}
  \label{cd_eps1}
  &\eps\partial_t\mu_\ept + \partial_t \varphi - \Delta\mu = 
  {P\sigma}h(\varphi_\ept) + (P\sigma_\tau-A)(h(\varphi_{\ept})-h(\varphi_\tau))
  \qquad&&\text{in } Q\,,\\
  \label{cd_eps2}
  &\mu=\tau\partial_t\varphi + a\varphi - J*\varphi + F'(\varphi_\ept)-F'(\varphi_\tau) -\chi\sigma
  \qquad&&\text{in } Q\,,\\
  \label{cd_eps3}
  &\partial_t\sigma - \Delta\sigma + B\sigma + C\sigma h(\varphi_\ept) 
  = C\sigma_\tau(h(\varphi_\tau)-h(\varphi_\ept))
   \qquad&&\text{in } Q\,,\\
  \label{cd_eps4}
  &\partial_{\bf n}\mu = \partial_{\bf n}\sigma = 0
  \qquad&&\text{on } \Sigma\,,\\
   \label{cd_eps5}
 &\varphi(0)=\ph_0, \quad \sigma(0)=\sigma_0
  \qquad&&\text{in } \Omega\, ,
\end{align}
where the equations are intended in the usual variational setting, and  
where we have set
$\varphi:=\varphi_\ept-\varphi_\tau$, $\mu:=\mu_\ept-\mu_\tau$, $\sigma:=\sigma_\ept-\sigma_\tau$,
$\ph_0:= \ph_{0, \ept}-\ph_{0, \tau}, \sigma_0:= \sigma_{0, \ept}-\sigma_{0, \tau}$.
Next, we multiply 
\eqref{cd_eps1} by $\tau \mu$,
\eqref{cd_eps2} by $\mu-\ph$, \eqref{cd_eps3} by $\sigma$,
add the resulting equality and integrate over $Q_t$ to obtain,
thanks to assumption {\bf A5},
\begin{align*}
	&  \intQt |\mu|^2
	+ \tau \intQt |\nabla \mu|^2
	+ \frac \tau 2 \norma{\ph(t)}^2_H
	+ C_0 \intQt |\ph|^2
	+ \frac 12 \norma{\sigma(t)}^2_H
	+ \intQt |\nabla\sigma|^2
		+\intQt(B + Ch(\varphi_\ept)) |\sigma|^2 
	\\ & \leq 
	\frac\tau2\norm{\ph_0}_H^2 +\frac12\norm{\sigma_{0}}_H^2
	-\intQt \eps\partial_t\mu_\ept \tau \mu
	+ \intQt (\mu+\chi\sigma) \ph
	+ \intQt C\sigma_\tau(h(\varphi_\tau)-h(\varphi_\ept)) \sigma\\
	&\qquad+ \tau\intQt \big[  {P\sigma}h(\varphi_\ept) 
		+ (P\sigma_\tau-A)(h(\varphi_{\ept})-h(\varphi_\tau))\big]  \mu
	+\int_{Q_t}(a\varphi - J*\varphi + F'(\varphi_\ept)-F'(\varphi_\tau) -\chi\sigma)\mu\,.
\end{align*}
Let us estimate the terms on the right-hand side separately.
The third and fourth ones yield, 
thanks to the Young inequality and the refined estimate \eqref{error1_eps},
\begin{align*}
-\intQt \eps\partial_t\mu_\ept \tau \mu
	+ \intQt (\mu+\chi\sigma) \ph &\leq
\frac12\int_{Q_t}|\mu|^2 + \tau^2\eps^2\norm{\partial_t\mu_\ept}_{L^2(0,T; H)}^2
+2\int_{Q_t}|\varphi|^2 + \frac{\chi^2}4 \int_{Q_t}|\sigma|^2\\
&\leq M \eps^{1/2} + \frac12\int_{Q_t}|\mu|^2 +2\int_{Q_t}|\varphi|^2 
+ \frac{\chi^2}4 \int_{Q_t}|\sigma|^2\,,
\end{align*}
for a certain constant $M$ independent of $\eps$.
The fifth and sixth terms can be easily handled using the Young inequality,
the Lipschitz continuity and boundedness of $h$, and the uniform bound
$\norm{\sigma_\tau}_{L^\infty(Q)}\leq 1$, as
\begin{align*}
  &\tau\intQt \big(  {P\sigma}h(\varphi_\ept) 
  + (P\sigma_\tau-A)(h(\varphi_{\ept})-h(\varphi_\tau))\big) \mu
  + \intQt C\sigma_\tau(h(\varphi_\tau)-h(\varphi_\ept)) \sigma\\
  &\leq\frac14\int_{Q_t}|\mu|^2 + 
  \norm{h}_{W^{1,\infty}(\erre)}^2\left((2\tau^2P^2+C^2)\int_{Q_t}|\sigma|^2
  +\left(\frac14+2\tau^2(P+A)^2\right)\int_{Q_t}|\varphi|^2\right)\,.
\end{align*}
Moreover, the last term satisfies, thanks to the Young inequality and 
the growth assumption \eqref{pot_reg},
\begin{align*}
&\int_{Q_t}(a\varphi - J*\varphi + F'(\varphi_\ept)-F'(\varphi_\tau) -\chi\sigma)\mu\\
&\leq\frac18\int_{Q_t}|\mu|^2 + 12(a^*)^2\int_{Q_t}|\varphi|^2
+6\chi^2\int_{Q_t}|\sigma|^2
+C_F\int_{Q_t}(1+|\varphi_\ept|^2+|\varphi_\tau|^2)|\varphi||\mu|\,,
\end{align*}
where, thanks to the inclusion $V\embed L^6(\Omega)$ and the H\"older inequality,
\begin{align*}
  \int_{Q_t}(1+|\varphi_\ept|^2+|\varphi_\tau|^2)|\varphi||\mu|&\leq
  \int_0^t\left(|\Omega|^{1/3} + \norm{\varphi_\ept}_{L^6(\Omega)}^2
  +\norm{\varphi_\tau}_{L^6(\Omega)}^2\right)
  \norm{\varphi(s)}_H\norm{\mu(s)}_{L^6(\Omega)}\,\d s\\
  &\leq M\left(1+\norm{\varphi_\ept}_{L^\infty(0,T; V)}^2
  +\norm{\varphi_\tau}_{L^\infty(0,T; V)}^2\right)\int_0^t\norm{\varphi(s)}_H\norm{\mu(s)}_V\,\d s\,,
\end{align*}
which yields, thanks to the estimate \eqref{est1_eps} and again the Young inequality, that 
\[
\int_{Q_t}(1+|\varphi_\ept|^2+|\varphi_\tau|^2)|\varphi||\mu| \leq
\min\{1/16, \tau/2\}
\norm{\mu}^2_{L^2(0,t; V)} + M_\tau\int_{Q_t}|\varphi|^2
\]
for a certain constant $M_\tau>0$ independent of $\eps$.
Hence, collecting the above estimates we obtain
\begin{align*}
	&  \min\{1/16, \tau/2\}\norm{\mu}^2_{L^2(0,t; V)}
	+ \frac \tau 2 \norma{\ph(t)}^2_H
	+ \frac 12 \norma{\sigma(t)}^2_H
	+ \intQt |\nabla\sigma|^2
	\\ & \leq M\left(\eps^{1/2} +
	\frac\tau2\norm{\ph_0}_H^2 +\frac12\norm{\sigma_{0}}_H^2+
	 \int_{Q_t}|\varphi|^2 + \int_{Q_t}|\sigma|^2\right)\,,
\end{align*}
where the updated constant $M$ depends on $\tau$, and the initial data 
$(\varphi_{0,\tau}, \sigma_{0,\tau})$. The error estimate follows then
by the Gronwall lemma.

Finally, it is not difficult to check that 
\luca{exactly the same argument performed here
yields uniqueness of the solution $(\varphi_\tau, \mu_\tau, \sigma_\tau,\xi_\tau)$
for the system \eqref{eq1}--\eqref{eq5} at $\eps=0$,
even without assumption \eqref{init4_eps0}.}
This reality implies then 
that the convergences as $\eps\searrow0$ hold along the entire sequence $\eps$
which completes the proof of Theorem~\ref{thm7}.

\section{Asymptotics as $\tau\searrow0$}
\label{sec:tau}
Let us now investigate the behavior of system \sys\ as $\tau \searrow 0$
by proving Theorems \ref{thm8} and \ref{thm9}.
Proceeding as before, notice that throughout this section we assume 
\luca{$\eps\in(0,\eps_0)$} to be fixed.

\subsection{Uniform estimates}
\label{ssec:unif_tau}
Performing the same estimates as in Subsection~\ref{ssec:unif_est},
and noting that the constant $M$ in \eqref{est_aux3} is independent of 
$\tau$, $\eps$, $\lambda$, and $n$, we infer that 
\begin{align}
  \nonumber
  &\frac\eps2\norm{\mu_{\ept}(t)}_H^2 
  + (1+4c_a\eps)\int_{Q_t}|\nabla\mu_{\ept}|^2
  +\tau\int_{Q_t}|\partial_t\varphi_{\ept}|^2
  +\frac12\norm{\sigma_{\ept}(t)}_H^2 + \int_{Q_t}|\nabla\sigma_{\ept}|^2\\
  \nonumber
  &\qquad+2c_a\norm{(\eps\mu_{\ept}+\varphi_{\ept})(t)}_H^2+
  2c_a\tau\norm{\nabla\varphi_{\ept}(t)}_H^2 + 
  2c_aC_0\int_{Q_t}|\nabla\varphi_{\ept}|^2\\
    \notag
  &\leq \frac32\eps\norm{\mu_{0,\ept}}_H^2 +
  (a^*+4c_a)\norm{\varphi_{0,\ept}}_H^2+ 2c_a\tau\norm{\nabla\varphi_{0,\ept}}_H^2+
   \norm{F(\varphi_{0,\ept})}_{L^1(\Omega)}\\ 
   \nonumber
  &\qquad+\frac12\norm{\sigma_{0,\ept}}_H^2
  +\frac{c_a}2\norm{\varphi_{\ept}(t)}_H^2
  +\chi\int_{Q_t}\sigma_{\ept}\partial_t\varphi_{\ept}+
  (\chia+4c_a\chi) \int_{Q_t}\nabla\sigma_{\ept}\cdot\nabla\varphi_{\ept}\\
	\label{aux1_tau}
 &\qquad+M\left(1+\int_{Q_t}|\eps\mu_{\ept}+\varphi_{\ept}|^2
  + \int_{Q_t}|\varphi_{\ept}|^2
   + \int_{Q_t}|\sigma_{\ept}|^2\right)
   +\int_{Q_t}(P\sigma_{\ept}-A)h(\varphi_{\ept})\mu_{\ept}\,.
\end{align}
First of all, note that all the terms on the right-hand side
referring to the initial data are uniformly bounded in $\tau$ due to 
assumptions \eqref{init2_tau0}--\eqref{init3_tau0}.
Moreover, since $\eps\in(0,\frac1{4c_a})$ we have
a bound from below on the left-hand side in the form
\begin{align}
  \notag 2c_a\norm{(\eps\mu_{\ept}+\varphi_{\ept})(t)}_H^2
  +\frac\eps2\norm{\mu_{\ept}(t)}_H^2&\geq
  2c_a\norm{(\eps\mu_{\ept}+\varphi_{\ept})(t)}_H^2
  +2c_a\eps^2\norm{\mu_{\ept}(t)}_H^2\\
  &\geq (c_a-\rho)\norm{\varphi_{\ept}(t)}_H^2 + 2\rho\eps^2\norm{\mu_\ept(t)}_H^2
  \label{unif_est_uno}
\end{align}
for every $\rho\in(0,c_a)$.
Hence the corresponding term $\frac{c_a}2\norm{\varphi_{\ept}(t)}_H^2$
on the right-hand side can be incorporated on the left-hand side of \eqref{aux1_tau},
provided we choose $\rho<c_a/2$.
Furthermore, from the boundedness of $h$
the last term in \eqref{aux1_tau}
can be easily handled using the Young inequality and the Gronwall lemma.
Hence, we only need to estimate the terms involving $\chi$ and $\chia$.
To this end, we first use integration by parts and the equation \eqref{var2} to deduce,
thanks to the Young inequality and the boundedness of $h$, that
\begin{align}
  \notag\chi\int_{Q_t}\sigma_{\ept}\partial_t\varphi_{\ept}&=
  -\chi\int_0^t\ip{\partial_t\sigma_{\ept}(s)}{\varphi_{\ept}(s)}\,\d s
  +\chi\int_\Omega\sigma_\ept(t)\varphi_\ept(t)-\chi\int_\Omega\sigma_{0,\ept}\varphi_{0,\ept}\\
  &\notag=\chi\int_{Q_t}\nabla\sigma_\ept\cdot\nabla\varphi_\ept+
  \chi\int_{Q_t}(B(\sigma_\ept-\sigma_S)+ C \sigma_\ept h(\varphi_\ept))\varphi_\ept
  -\chi\chia\int_{Q_t}|\nabla\varphi_\ept|^2\\
  &\notag\qquad+\chi\int_\Omega\sigma_\ept(t)\varphi_\ept(t)
  -\chi\int_\Omega\sigma_{0,\ept}\varphi_{0,\ept}\\
  &\notag\leq\chi\int_{Q_t}\nabla\sigma_\ept\cdot\nabla\varphi_\ept
  -\chi\chia\int_{Q_t}|\nabla\varphi_\ept|^2
  +M\left(1+\int_{Q_t}|\varphi_\ept|^2 + \int_{Q_t}|\sigma_\ept|^2\right)\\
  &\qquad+\delta\chi^2\norm{\varphi_\ept(t)}_H^2 + \frac1{4\delta}\norm{\sigma_\ept(t)}_H^2\,,
  \label{last:2}
\end{align}
for every $\delta>0$. Now, 
it is immediate to check that
assumption \eqref{ip_chi} yields $\frac12<\frac{c_a}{2\chi^2}$
(with the convention that $\frac{1}{\chi}=+\infty$ if $\chi=0$): hence
\begin{align}
	\label{deltabar}
	\exists \, \bar\delta \in \Big(\frac12, \frac{c_a}{2\chi^2}\Big) \quad \hbox{such that} \quad
  \bar\delta\chi^2<\frac{c_a}{2}\,, \qquad 
  \frac1{4\bar\delta}<\frac12\,,
\end{align}  
so that we can
incorporate the last two terms on the right-hand side 
of \eqref{last:2} on the left-hand side of \eqref{aux1_tau}.
Taking these remarks into account, we are left with 
\begin{align}
  \nonumber
  &2\rho\eps^2\norm{\mu_{\ept}(t)}_H^2 +
  \left(\frac{c_a}2  - \rho- \bar\delta\chi^2  \right)\norm{\varphi_\ept(t)}_H^2
  + (1+4c_a\eps)\int_{Q_t}|\nabla\mu_{\ept}|^2
  +\tau\int_{Q_t}|\partial_t\varphi_{\ept}|^2\\
  \nonumber
  &\qquad+\left(\frac12-\frac1{4\bar\delta}\right)
  \norm{\sigma_{\ept}(t)}_H^2
  + \int_{Q_t}|\nabla\sigma_{\ept}|^2+
  2c_a\tau\norm{\nabla\varphi_{\ept}(t)}_H^2 + 
  (2c_aC_0 + \chi\chia)\int_{Q_t}|\nabla\varphi_{\ept}|^2\\
  \label{aux2_tau}
  &\leq M\left(1 + \int_{Q_t}|\mu_{\ept}|^2
  + \int_{Q_t}|\varphi_{\ept}|^2
   + \int_{Q_t}|\sigma_{\ept}|^2\right)
  +(\chi+\chia+4c_a\chi) \int_{Q_t}\nabla\sigma_{\ept}\cdot\nabla\varphi_{\ept}\,,
\end{align}
which holds for every $\rho\in(0,c_a/2)$. By choosing $\bar\delta$ such that \eqref{deltabar} 
are fulfilled, it is also possible to choose and fix $\bar\rho\in(0,c_a/2)$ such that 
\[
  \frac{c_a}2 - \bar\rho - \bar\delta\chi^2 >0\,.
\]
Next, we use again the averaged Young inequality to obtain,
for every $\kappa>0$,
\[
(\chi+\chia+4c_a\chi) \int_{Q_t}\nabla\sigma_{\ept}\cdot\nabla\varphi_{\ept}\leq
\kappa\int_{Q_t}|\nabla\sigma_\ept|^2 
+ \frac{(\chi+\chia+4c_a\chi)^2}{4\kappa}\int_{Q_t}|\nabla\varphi_\ept|^2\,
\]
where the two terms on the \rhs\ can be incorporated on the left-hand side of \eqref{aux2_tau}
provided to choose $\kappa$ such that 
\[
  \kappa<1\,, \qquad
  \frac{(\chi+\chia+4c_a\chi)^2}{4\kappa}<2c_aC_0 + \chi\chia\,.
\]
Easy computations show that this is possible if and only if
\[
  \frac{(\chi+\chia+4c_a\chi)^2}{4(2c_aC_0 + \chi\chia)}<1
\]
which is verified owing to \eqref{ip_chi}.

Therefore, after rearranging the terms and using the Gronwall lemma, we infer that 
there exists a constant $M>0$, which may depend on $\eps$, but it is independent of $\tau,$
such that
\begin{align}
 \label{est1_tau}
 \norm{\varphi_\ept}_{L^\infty(0,T; H)\cap L^2(0,T; V)} + 
 \norm{\mu_\ept}_{L^\infty(0,T; H)\cap L^2(0,T; V)}+
 \norm{\sigma_\ept}_{L^\infty(0,T; H)\cap L^2(0,T; V)}&\leq M\,,\\
 \label{est2_tau}
 \tau^{1/2}\norm{\varphi_\ept}_{H^1(0,T; H)}+
 \tau^{1/2}\norm{\varphi_\ept}_{L^\infty(0,T; V)} &\leq M\,,
\end{align}
yielding in turn, by comparison in equations \eqref{eq1} and \eqref{eq3},
\beq
  \label{est3_tau}
  \norm{\eps\mu_\ept+\varphi_\ept}_{H^1(0,T; V^*)} + \norm{\sigma_\ept}_{H^1(0,T; V^*)} \leq M\,.
\eeq
Testing equation \eqref{eq2} by $\xi_\ept$ and using the estimate \eqref{est1_tau},
its is a standard matter to deduce also that 
\beq
  \label{est4_tau}
  \norm{\xi_\ept}_{L^2(0,T; H)}\leq M\,.
\eeq

\subsection{Passage to the limit}
The estimates \eqref{est1_tau}--\eqref{est4_tau} 
and classical compactness results
(see, e.g., \cite[Sec.~8, Cor.~4]{simon}) ensure that there exists
a quadruplet $(\varphi_\eps, \mu_\eps, \sigma_\eps, \xi_\eps)$ with 
\begin{align*}
  &\varphi_\eps, \mu_\eps\in L^\infty(0,T; H)\cap L^2(0,T; V)\,, \qquad
  \lambda_\eps:=\eps\mu_\eps+\varphi_\eps\in H^1(0,T; V^*) \cap L^2(0,T; V)\,,\\
  &\sigma_\eps\in H^1(0,T; V^*)\cap L^2(0,T; V)\,, \qquad
  \xi_\eps\in L^2(0,T; H)\,,
\end{align*}
such that, as $\tau\searrow0$ (on a subsequence) it holds that \tautozero\ and
\eqref{tautozero:strong:1}--\eqref{tautozero:strong:2}
are satisfied, and also that
  \[
	\lambda_\ept\wto \lambda_\eps \quad\hbox{in } H^1(0,T; V^*) \cap \L2 V\,,\qquad
	\lambda_\ept \to \lambda_\eps \quad\hbox{in } C^0([0,T]; V^*)\cap L^2(0,T; H)\,.
  \]
Moreover, let us claim that the above strong convergences imply the strong convergences 
\begin{align}
	\label{strong_conv}
	\mu_\ept \to \mu_\eps
	\quad \hbox{in } \L2 H\,,
	\qquad
	\ph_\ept \to \ph_\eps
	\quad \hbox{in } \L2 H\,.
\end{align}
To this end, we argue as in \cite[Sec.~3]{col-gil-roc-spr2}, checking that the sequence
$\{ \lambda_\ept \}_\tau$ is a Cauchy sequence in $L^2(0,T; H)$.
Let us pick two arbitrary $\tau,\tau'>0$ and take the difference of 
the corresponding equation \eqref{eq2} for $\tau$ and $\tau'$.
Next, we multiply the resulting equation by $\eps$, add to both sides
$\ph_\ept-\ph_{\ept'}$, test the resulting equation by $\ph_\ept-\ph_{\ept'}$, and integrate over $Q_t$
to obtain
\begin{align*}
	& \int_{Q_t}|\ph_\ept-\ph_{\ept'}|^2
	+ \eps\int_{Q_t} 
	\left(a(\ph_\ept-\ph_{\ept'}) +
	\xi_\ept-\xi_{\ept'}+F_2'(\ph_\ept)-F_2'(\ph_{\ept'})\right)
	(\ph_\ept-\ph_{\ept'})
	\\ &\leq
	\int_{Q_t} \big( (\lambda_\ept-\lambda_{\ept'})
	-\eps(\tau \partial_t \ph_\ept - \tau' \partial_t \ph_{\ept'}) 
	+ \eps\chi (\sigma_\ept - \sigma_{\ept'})
	\big)(\ph_\ept-\ph_{\ept'})\\
	&\qquad+\eps\int_{Q_t}J * (\ph_\ept-\ph_{\ept'})(\ph_\ept-\ph_{\ept'})\,.
\end{align*}
Owing to \tautozero\ and \eqref{tautozero:strong:1}--\eqref{tautozero:strong:2} 
we easily infer that the first term on the \rhs\ goes to zero as $\tau,\tau' \to 0$.
Moreover, on the left-hand side we have, thanks to assumption {\bf A5},
\[ 
\int_{Q_t} \left(a(\ph_\ept-\ph_{\ept'}) +\xi_\ept-\xi_{\ept'}
		+F_2'(\ph_\ept)-F_2'(\ph_{\ept'}) \right)(\ph_\ept-\ph_{\ept'})
		\geq  C_0\int_{Q_t}|\ph_\ept-\ph_{\ept'}|^2\,,
\]
while the last term on the right-hand side satisfies 
\[
  \int_{Q_t}J * (\ph_\ept-\ph_{\ept'})(\ph_\ept-\ph_{\ept'})\leq 
  a^*\int_{Q_t}|\ph_\ept-\ph_{\ept'}|^2\,.
\]
Rearranging the terms leads us to 
\begin{align*}
	& (1+(C_0-a^*)\eps)\int_{Q_t}|\ph_\ept-\ph_{\ept'}|^2
	\leq
	\int_{Q_t} \big( (\lambda_\ept-\lambda_{\ept'})
	-\eps(\tau \partial_t \ph_\ept - \tau' \partial_t \ph_{\ept'}) 
	+ \chi (\sigma_\ept - \sigma_{\ept'})
	\big)(\ph_\ept-\ph_{\ept'})
\end{align*}
where the right-hand side converges to $0$ as $\tau\searrow0$.
\luca{Since $\eps a^*<\eps C_0 + 1$ as a consequence
of the smallness assumption on $\eps_0$,}
this yields the second of 
\eqref{strong_conv} and by comparison also the first one follows, as we claimed.

With the strong convergence of the phase variable at disposal
it is now straightforward to infer by combining the boundedness of $h$ and
the Lebesgue convergence theorem that, as $\tau \searrow 0$, 
\begin{align*}
	h(\ph_\ept) \to h(\ph_\eps) \quad\text{in } L^p(Q) \quad\forall\,p\geq1\,,
	\qquad 
	F_2'(\ph_\ept) \to F_2'(\ph_\eps) 
	\quad \hbox{in } \L2 H\,.
\end{align*} 
Hence, since $\xi_\eps\in\partial F_1(\varphi_\eps)$
by the strong-weak closure of $\partial F_1$,
it is a standard matter to pass to the limit as $\tau\searrow0$
in the weak formulation of \sys\, and 
deduce that the limit 
$(\mu_\eps, \ph_\eps, \sigma_\eps, \xi_\eps)$
yields a solution to \sys\ with $\tau =0$.
\luca{Notice in particular that by difference in the limit 
equation \eqref{eq2} we deduce the further regularity $\xi_\eps\in L^2(0,T; V)$,}
while the last assertion of Theorem \ref{thm8} follows as before by repeating the computations
of Subsection \ref{sub:maxpr} completing the proof of Theorem~\ref{thm8}.

\subsection{Error estimate}
The last result of this section follows with few changes from 
the proof of the continuous dependence estimate
\eqref{cont_dep_first} established in Theorem~\ref{thm2}.

Indeed, we can repeat almost the same computations performed 
in Subsection~\ref{ssec:cont} with the choices 
\begin{align*}
	(\varphi_1, \mu_1, \sigma_1, \xi_1):=(\varphi_\ept, \mu_\ept, \sigma_\ept, \xi_\ept)\,,
	\quad
	(\varphi_2, \mu_2, \sigma_2, \xi_2):=(\varphi_\eps, \mu_\eps, \sigma_\eps, \xi_\eps)\,.
\end{align*}
Moreover, by setting $\varphi:=\varphi_\ept-\varphi_\eps$, $\mu:=\mu_\ept-\mu_\eps$, 
$\sigma:=\sigma_\ept-\sigma_\eps$, $\varphi_0:=\varphi_{0,\ept}-\varphi_{0,\eps}$,
$\mu_0:=\mu_{0,\ept}-\mu_{0,\eps}$, and $\sigma_{0}:=\sigma_{0,\ept}-\sigma_{0,\eps}$,
\luca{recalling that we are assuming $\chia=0$
we infer from \eqref{last:1} that
\begin{align*}
  &\frac12\norm{(\eps\mu+\varphi)(t)}_{V^*}^2+\eps\int_{Q_t}|\mu|^2 
  + C_0\int_{Q_t}|\varphi|^2
  +\frac12\norm{\sigma(t)}_H^2 + \int_{Q_t}|\nabla\sigma|^2\\
  &\leq -\tau\int_{Q_t}\partial_t\varphi_\ept\varphi +
  \frac12\norm{\eps\mu_0+\varphi_0}_{V^*}^2  
  + \frac12\norm{\sigma_0}_H^2
  +\int_{Q_t}(\chi\sigma + J*\varphi)\varphi
  +\int_{Q_t}[C\sigma_{\eps}(h(\varphi_\eps)-h(\varphi_\ept))]\sigma\\
  &\qquad+
  \int_{Q_t}\big[\mu + P\sigma h(\varphi_\ept)+ 
  (P\sigma_\eps-A)(h(\varphi_{\ept})-h(\varphi_\eps))\big]\mathcal R^{-1}(\eps\mu+\varphi)  \,.
\end{align*}}
All the terms on the right-hand side, except the first one, 
can be handled in exactly the same way as in Subsection~\ref{ssec:cont}.
As for the first one, we use the Young inequality and estimate \eqref{est2_tau} to infer,
\luca{for every $\delta>0$,
\begin{align*}
	-\tau \int_{Q_t} \partial_t\varphi_\ept \ph
	\leq \delta \int_{Q_t} |\ph|^2
	+ \frac {\tau^2} {4\delta} \int_{Q_t}|\partial_t\varphi_\ept|^2
	\leq \delta  \int_{Q_t} |\ph|^2
	+ M_\delta \tau\,,
\end{align*}
so that the first term on the \rhs\ can be absorbed on the left
provided to choose again $\delta$ small enough, which is indeed
possible as we noted in Subsection~\ref{ssec:cont}.}
We can now argue as before and conclude using Gronwall's lemma.
Moreover, the same argument on the limit problem yields
uniqueness of solution for the system with $\tau=0$, hence also 
that the convergences hold along the entire sequence and the proof of Theorem~\ref{thm9}
is concluded.

\section{Asymptotics as both $\eps,\tau\searrow0$}
\label{sec:eps_tau}

The last issue we are going to address here concerns the joint asymptotic limit as both $\eps,\tau\searrow0$.
Let us recall that in this section we are supposing that 
$\chia=0$.

\subsection{Uniform estimates}
Let us come back to estimate \eqref{est_aux3}
with $\chia=0$. We have
\begin{align}
  \nonumber
  &\frac\eps2\norm{\mu_{\ept}(t)}_H^2 
  + (1+4c_a\eps)\int_{Q_t}|\nabla\mu_{\ept}|^2
  +\tau\int_{Q_t}|\partial_t\varphi_{\ept}|^2
  +\int_\Omega F(\varphi_\ept(t))
  +\frac12\norm{\sigma_{\ept}(t)}_H^2\\
   \nonumber
  &\qquad+ \int_{Q_t}|\nabla\sigma_{\ept}|^2
  +2c_a\norm{(\eps\mu_{\ept}+\varphi_{\ept})(t)}_H^2+
  2c_a\tau\norm{\nabla\varphi_{\ept}(t)}_H^2 + 
  2c_aC_0\int_{Q_t}|\nabla\varphi_{\ept}|^2\\
    \notag
  &\leq \frac32\eps\norm{\mu_{0,\ept}}_H^2 +
  (a^*+4c_a)\norm{\varphi_{0,\ept}}_H^2+ 
  2c_a\tau\norm{\nabla\varphi_{0,\ept}}_H^2+
   \norm{F(\varphi_{0,\ept})}_{L^1(\Omega)}\\ 
   \nonumber
  &\qquad+\frac12\norm{\sigma_{0,\ept}}_H^2
  +\frac{c_a}2\norm{\varphi_{\ept}(t)}_H^2
  +\chi\int_{Q_t}\sigma_{\ept}\partial_t\varphi_{\ept}+
  4c_a\chi \int_{Q_t}\nabla\sigma_{\ept}\cdot\nabla\varphi_{\ept}\\
	\label{aux1_ept}
	 &\qquad+M\left(1+\int_{Q_t}|\eps\mu_{\ept}+\varphi_{\ept}|^2
  + \int_{Q_t}|\varphi_{\ept}|^2
   + \int_{Q_t}|\sigma_{\ept}|^2\right)
   +\int_{Q_t}(P\sigma_{\ept}-A)h(\varphi_{\ept})\mu_{\ept}\,,
\end{align}
where the constant $M>0$ is independent of both $\eps$ and $\tau$.
Now, all the terms on the right-hand side 
referring to the initial data are uniformly bounded in both
$\eps$ and $\tau$ thanks to assumptions \eqref{init2_ept}--\eqref{init3_ept}.
Moreover, as done in \eqref{unif_est_uno},
on the left-hand side we have 
\[
  2c_a\norm{(\eps\mu_{\ept}+\varphi_{\ept})(t)}_H^2
  +\frac\eps2\norm{\mu_{\ept}(t)}_H^2
  \geq (c_a-\rho)\norm{\varphi_{\ept}(t)}_H^2 + 2\rho\eps^2\norm{\mu_\ept(t)}_H^2
\]
for every $\rho\in(0,c_a/2)$, so that the term on the right-hand side of the above inequality
can be absorbed on the \lhs\ of \eqref{aux1_ept}. Furthermore, proceeding again as in 
Subsection~\ref{ssec:unif_tau} and recalling that here $\chia=0$, we have
\begin{align*}
  &\chi\int_{Q_t}\sigma_{\ept}\partial_t\varphi_{\ept}=
  -\chi\int_0^t\ip{\partial_t\sigma_{\ept}(s)}{\varphi_{\ept}(s)}\,\d s
  +\chi\int_\Omega\sigma_\ept(t)\varphi_\ept(t)-\chi\int_\Omega\sigma_{0,\ept}\varphi_{0,\ept}\\
  &\qquad\leq\chi\int_{Q_t}\nabla\sigma_\ept\cdot\nabla\varphi_\ept
  +M\left(1+\int_{Q_t}|\varphi_\ept|^2 + \int_{Q_t}|\sigma_\ept|^2\right)
  +\delta\chi^2\norm{\varphi_\ept(t)}_H^2 + \frac1{4\delta}\norm{\sigma_\ept(t)}_H^2\,,
\end{align*}
for every $\delta>0$.
Moreover,
we can choose $\bar \delta $ such that \eqref{deltabar} are satisfied,
so that the corresponding two terms on the right-hand side can be incorporated on the left.
The remaining terms on the right-hand side of \eqref{aux1_ept}
 containing $\chi$ can be handled as,
for every $\kappa>0$,
\[
(\chi+4c_a\chi) \int_{Q_t}\nabla\sigma_{\ept}\cdot\nabla\varphi_{\ept}\leq
\kappa\int_{Q_t}|\nabla\sigma_\ept|^2 
+ \frac{(\chi+4c_a\chi)^2}{4\kappa}\int_{Q_t}|\nabla\varphi_\ept|^2\,.
\]
Again, the two terms on the right can be incorporated on the left-hand side of \eqref{aux1_ept}
provided that we choose $\kappa$ such that 
\[
  \kappa<1\,, \qquad
  \frac{(\chi+4c_a\chi)^2}{4\kappa}<2c_aC_0\,,
\]
which is indeed possible since \eqref{ip_chi} and the fact that $\chia=0$ yield $\frac{(\chi+4c_a\chi)^2}{8c_aC_0}<1$.
To close the estimate, 
we only need to handle the last term on the right-hand side
of \eqref{aux1_ept}: this can be done 
exactly in the same way as in Subsection~\ref{ssec:unif_eps}.
Indeed, on the right-hand side we have, thanks to the boundedness of $h$
and the fact that $\norm{\sigma_{\ept}}_{L^\infty(Q)}\leq1$,
\[
	\int_{Q_t}(P\sigma_\ept  - A )h(\varphi_{\ept})\mu_{\ept}
 	\leq  
 	\frac 12 \int_{Q_t} |\nabla \mu_{\ept}|^2 + M'\left(1+
 	 T_0^{1/2}\norm{(\mu_\ept)_{\Omega}}_{L^2(0,t)}\right)
\]
for every $t\in[0,T_0]$ and $T_0<T$, 
where $M'$ only depends on $P$, $A$, and $h$.
Furthermore, by comparison in equation \eqref{eq2} and
thanks to \eqref{pol_growth}, since $\tau\in(0,1)$, we have
\[
  |(\mu_\ept(t))_\Omega| 
  \leq M''\left(1 + \tau\int_{Q_t}|\partial_t\varphi_\ept(t)|^2
  +\sup_{s\in[0,t]}\int_\Omega F(\varphi_\ept(s)) + \sup_{s\in[0,t]}\norm{\sigma_\ept(s)}_H^2\right)\,,
\]
where $M''>0$ only depends on $C_F$ and $\chi$. Hence, 
using a patching argument as in Subsection~\ref{ssec:unif_eps}, we deduce the following 
uniform estimates
\begin{align}
  \label{est1_ept}
  \norm{\varphi_\ept}_{L^\infty(0,T; H)\cap L^2(0,T; V)} + 
  \norm{\mu_\ept}_{L^2(0,T; V)} + \norm{\sigma_\ept}_{L^\infty(0,T; H)\cap L^2(0,T; V)}
  &\leq M\,,\\
  \label{est1_bis_ept}
  \norm{F(\varphi_\ept)}_{L^\infty(0,T; L^1(\Omega))} &\leq M\,,\\
  \label{est2_ept}
  \eps^{1/2}\norm{\mu_\ept}_{L^\infty(0,T; H)}+
  \tau^{1/2}\norm{\varphi_\ept}_{H^1(0,T; H)\cap L^\infty(0,T; V)} &\leq M\,.
\end{align}
Comparison then in the system gives us in particular that 
\beq
  \label{est3_ept}
  \norm{\xi_\ept}_{L^2(0,T; H)} + \norm{\sigma_\ept}_{H^1(0,T; V^*)}+
  \norm{\eps\mu_\ept + \varphi_\ept}_{H^1(0,T; V^*)} \leq M\,,
\eeq
as well as
\beq
  \label{est4_ept}
  \eps\tau^{1/2}\norm{\mu_\ept}_{H^1(0,T; V^*)}\leq M\,.
\eeq
The uniform bound for $\sigma_\ept$ in $L^\infty(Q)$ can be obtained as before
using Subsection \ref{sub:maxpr}.

\subsection{Passage to the limit}
The estimates \eqref{est1_ept}--\eqref{est4_ept} ensure, thanks to the classical 
compactness results, that there exists a quadruplet $(\varphi, \mu, \sigma, \xi)$,
with 
\begin{align*}
  &\varphi\in H^1(0,T; V^*)\cap L^2(0,T; V)\,, \qquad\mu\in L^2(0,T; V)\,,\\
  &\sigma\in H^1(0,T; V^*)\cap L^2(0,T; V)\cap L^\infty(Q)\,, \qquad
  0\leq\sigma(t,x)\leq1 \quad\text{for a.e.~}x\in\Omega\,,\quad\forall\,t\in[0,T]\,,\\
  &\xi\in L^2(0,T; H)\,,
\end{align*}
such that, as $(\eps,\tau)\searrow0$ it holds that, along a non-relabelled subsequence,
\epstautozero\ and \eqref{epstautozero:strong} are fulfilled. In addition, setting 
$\lambda_\ept:=\eps\mu_\ept+\varphi_\ept$, we have
\begin{align*}
  \lambda_\ept\wto\varphi \quad&\text{in } H^1(0,T; V^*)\cap L^2(0,T; V)\,, \qquad
  \lambda_\ept\to\varphi \quad \text{in } C^0([0,T]; V^*)\cap L^2(0,T; H)\,,\\
  \xi_\ept\wto \xi \quad&\text{in } L^2(0,T; H)\,.
\end{align*}
In particular, by difference we deduce that 
\[
  \varphi_\ept=\lambda_\ept - \eps\mu_\ept \to \varphi \quad\text{in } L^2(0,T; H)
\]
which readily implies that $\xi\in\partial F_1(\varphi)$ almost everywhere in $Q$, and that 
\[
  h(\varphi_\ept)\to h(\varphi) \quad\text{in } L^p(Q) \quad\forall\,p\geq1\,, \qquad
  F_2'(\varphi_\ept)\to F_2'(\varphi) \quad\text{in } L^2(0,T; H)\,.
\]
It is then a standard matter to let $(\eps,\tau)\searrow0$ in the weak formulation 
of \eqref{eq1}--\eqref{eq5} to conclude. Note in particular that by difference in the limit 
equation \eqref{eq2} we deduce the further regularity $\xi\in L^2(0,T; V)$,
which concludes the proof of Theorem~\ref{thm10}.

\subsection{Error estimate}
In this last subsection we prove the error estimate as 
both $\eps$ and $\tau$ go to zero.

The idea is to adapt the argument presented in Subsection~\ref{ssec:error_eps}.
First of all, we need to prove a refined estimate:
proceeding as in Subsection~\ref{ssec:error_eps}, we know that 
\begin{align}
  \nonumber
  &\eps^{3/2}\int_{Q_t}|\partial_t \mu_\ept|^2 
  + \frac{\eps^{1/2}}2 \norma{\nabla\mu_\ept(t)}_H^2
  +\frac{\tau\eps^{1/2}}2\norma{\partial_t\varphi_\ept(t)}_H^2
  + C_0\eps^{1/2}\int_{Q_t}|\partial_t\varphi_\ept|^2\\
  \notag
  &\qquad
  +\eps^{1/2}\int_{Q_t}|\partial_t\sigma_\ept|^2 
  + \frac{\eps^{1/2}}2\norma{\nabla\sigma_\ept(t)}_H^2\\
  \nonumber
  &\leq\frac{\eps^{1/2}}2\norma{\nabla\mu_{0,\ept}}_H^2 + 
  \frac{\tau\eps^{1/2}}2\norma{\varphi'_{0,\ept}}_H^2
  +\frac{\eps^{1/2}}2\norma{\nabla\sigma_{0,\ept}}_H^2
  +\eps^{1/2}\int_{Q_t}(P\sigma_\ept-A)h(\varphi_\ept)\partial_t\mu_\ept\\
	\label{aux2_ept}
	  &\qquad+\eps^{1/2}\int_{Q_t}(J*(\partial_t\varphi_\ept) 
  + \chi\partial_t\sigma_\ept)\partial_t\varphi_\ept
  +\eps^{1/2}\int_{Q_t}\left(B(\sigma_S-\sigma_\ept)-Ch(\varphi_\ept)\sigma_\ept\right)
  \partial_t\sigma_\ept\,.
\end{align}
The first and third terms on the right-hand side are uniformly bounded in 
$\eps$ and $\tau$ due to assumptions \eqref{init3_ept} and \eqref{init4_ept}.
As for the second term on the right-hand side, using \eqref{eq2} we realize that
\[
  \mu_{0,\ept}=\tau \varphi'_{0,\ept}+ a\varphi_{0,\ept}-J*\varphi_{0,\ept}
  + F'(\varphi_{0,\ept}) - \chi\sigma_{0,\ept}\,,
\]
so that, multiplying both sides by $\eps^{1/4}/\tau^{1/2}$ and squaring,
\[
  \tau\eps^{1/2}\norm{\varphi'_{0,\ept}}_H^2\leq
  5\frac{\eps^{1/2}}{\tau}\left(\norm{\mu_{0,\ept}}_H^2
  +2(a^*)^2\norm{\varphi_{0,\ept}}_H^2 + \norm{F'(\varphi_{0,\ept})}_H^2
  +\chi^2\norm{\sigma_{0,\ept}}_H^2\right)\,,
\]
from which we deduce by \eqref{init4_ept} that the second
term on the right-hand side of \eqref{aux2_ept} is uniformly bounded in $\eps$ and $\tau$.
Let us focus on the fourth term on the right-hand side:
proceeding as in Subsection~\ref{ssec:error_eps}, this can be bounded
using integration by parts and the Young inequality by the quantity 
\begin{align*}
  &\frac{\eps^{1/2}}{4}\int_{Q_t}|\partial_t\sigma_\ept|^2+
  M\eps^{1/2}\left(\norm{\mu_\ept}^2_{L^2(0,T; H)}+
  \norm{\partial_t\varphi_\ept}_{L^2(0,T; H)}^2+
  \eps^{1/2}\norm{\mu_\ept}_{C^0([0,T]; H)}\right)
\end{align*}
for a positive constant $M$ independent of $\eps$ and $\tau$.
The first term can be then incorporated on the left-hand side, 
and the remaining others are
uniformly bounded in $\eps$ and $\tau$ thanks to 
the estimates \eqref{est1_ept}, \eqref{est2_ept}, and condition \eqref{limsup}
on $(\eps,\tau)$.
\luca{Finally, noting that
\[
\eps^{1/2}\int_{Q_t}(J*\partial_t\varphi_\ept)\partial_t\varphi_\ept\leq
(a^*+b^*)\eps^{1/2}\int_0^t\norm{\partial_t\varphi_\ept(s)}_H
\norm{\partial_t\varphi{_\ept}(s)}_{V^*}\,\d s\,
\leq  M\eps^{1/2}\norm{\partial_t\varphi_\ept}_{L^2(0,T; H)}^2\,,
\]
the remaining terms on the right-hand side of \eqref{aux2_ept}
can be handled similarly, using the averaged Young inequality,
estimate \eqref{est1_ept}--\eqref{est2_ept}, and condition \eqref{limsup}.}
Thus, there exists $M>0$, independent of both $\eps$ and $\tau$, such that 
\begin{align}
  \label{error1_ept}
  \eps^{3/4}\norm{\mu_\ept}_{H^1(0,T; H)} 
  + \eps^{1/4}\norm{\mu_\ept}_{L^\infty(0,T; V)}&\leq M\,,\\
  \label{error2_ept}
  \tau^{1/2}\eps^{1/4}\norm{\varphi_\ept}_{W^{1,\infty}(0,T; H)}
  +\eps^{1/4}\norm{\sigma_\ept}_{H^1(0,T; H)\cap L^\infty(0,T; V)} &\leq M\,.
\end{align}

We are now ready to show the error estimate.
Setting $\ov\varphi:=\varphi_\ept-\varphi$, $\ov\mu:=\mu_\ept-\mu$,
$\ov\sigma:=\sigma_\ept-\sigma$, $\ov\varphi_0:=\varphi_{0,\ept}-\varphi_0$,
and $\ov\sigma_{0}:=\sigma_{0,\ept}-\sigma_0$, 
we write the difference of the system \eqref{eq1}--\eqref{eq5} with $\chia=0$
at $\eps,\tau>0$ and $\eps=\tau=0$ to find that
\begin{align}
  \label{cd_epstau1}
  &\eps\partial_t\mu_\ept +\partial_t \ov\varphi -\Delta \ov\mu = 
  P\ov\sigma h(\varphi_\ept) + (P \sigma -A)(h(\varphi_{\ept})-h(\varphi))
  \qquad&&\text{in } Q\,,\\
  \label{cd_epstau2}
  &\ov\mu=\tau\partial_t\varphi_\ept + a\ov\varphi - J*\ov\varphi + 
  F'(\varphi_\ept)-F'(\varphi) - \chi\ov\sigma
  \qquad&&\text{in } Q\,,\\
  \label{cd_epstau3}
  &\partial_t\ov\sigma -\Delta\ov\sigma + B\ov\sigma + C\ov\sigma h(\varphi_\ept) =
  C\sigma(h(\varphi)-h(\varphi_\ept))
   \qquad&&\text{in } Q\,,\\
  \label{cd_epstau4}
  &\partial_{\bf n}\ov\mu = \partial_{\bf n}\ov\sigma = 0
  \qquad&&\text{on } \Sigma\,,\\
   \label{cd_epstau5}
 & \ov\varphi(0)=\ov\varphi_0\,, \quad \ov\sigma(0)=\ov\sigma_0
  \qquad&&\text{in } \Omega\, ,
\end{align}
where the equations have to be intended in the usual variational framework.
We test \eqref{cd_epstau1} by $\mathcal N(\ov\varphi-(\ov\varphi)_\Omega)$, 
\eqref{cd_epstau2} by $\ov\varphi-(\ov\varphi)_\Omega$, 
\eqref{cd_epstau3} by $\ov\sigma$, integrate over $Q_t$, add the resulting equalities and use {\bf A5} to get
\begin{align}
\nonumber
	& \frac 12 \norma{(\ov\varphi-(\ov\varphi)_\Omega)(t)}_{V^*}^2
	+ \luca{C_0\intQt |\ov\ph|^2}
	+ \frac 12 \norma{\ov\sigma(t)}_H^2
	+ \intQt |\nabla\ov\sigma|^2
	+\int_{Q_t}(B+Ch(\varphi_\ept))|\ov\sigma|^2
	\\ 
	\nonumber
	&= \frac12\norm{\ov\ph_0-(\ov\ph_0)_\Omega}_{V^*}^2
	+\frac12\norm{\ov\sigma_0}_H^2
	-\eps\int_{Q_t}\partial_t\mu_\ept\mathcal{N}(\ov\varphi-(\ov\varphi)_\Omega)
	+\int_{Q_t}\ov\varphi(\chi\ov\sigma-\tau\partial_t\varphi_\ept)\\
	\nonumber
	&\qquad+\int_{Q_t}\ov\mu(\ov\varphi)_\Omega
	+\luca{\int_{Q_t}(J*\ov\varphi)\ov\varphi}
	+C\int_{Q_t}\sigma (h({\varphi})-h(\varphi_\ept)) \ov\sigma\\
	&\label{aux3_ept}
	\qquad+\int_{Q_t}\Big(P\ov\sigma h(\varphi_\ept) 
		+ (P\sigma -A)(h(\varphi_{\ept})-h(\varphi))
		\Big)\mathcal N(\ov\varphi-(\ov\varphi)_\Omega)\,.
\end{align}
Now, note that the Young inequality and the estimates \eqref{est2_ept}
and \eqref{error1_ept} yield
\begin{align*}
  &-\eps\int_{Q_t}\partial_t\mu_\ept\mathcal N(\ov\varphi-(\ov\varphi)_\Omega)
  +\int_{Q_t}\ov\varphi(\chi\ov\sigma-\tau\partial_t\varphi_\ept)\\
  &\leq \eps^2\norm{\partial_t\mu_\ept}_{L^2(0,T; H)}^2 
  +\frac14\int_0^t\norm{\mathcal N(\ov\varphi-(\ov\varphi)_\Omega)(s)}_{H}^2\,\d s
  +\luca{\frac{C_0}{4}\intQt |\ov\ph|^2}\\
  &\qquad+\luca{\frac{2}{C_0}}\left(\tau^2\norm{\partial_t\varphi_\ept}^2_{L^2(0,T; H)}+
  \chi^2\int_{Q_t}|\ov\sigma|^2\right)\\
  &\leq  \luca{\frac{C_0}{4}}\intQt |\ov\ph|^2
  + M\left(\eps^{1/2} + \tau + \int_0^t\norm{(\ov\varphi-(\ov\varphi)_\Omega)(s)}_{V^*}^2\,\d s
  +\int_{Q_t}|\ov\sigma|^2\right)
 \,,
\end{align*}
for a certain constant $M>0$ independent of $\eps$ and $\tau$.
Furthermore, using the boundedness 
and Lipschitz continuity of $h$, and the fact that $\norm{\sigma}_{L^\infty(Q)}\leq 1$,
the last two terms in \eqref{aux3_ept} can be handled again by the Young inequality as
\begin{align*}
  &C\int_{Q_t}\sigma (h(\varphi)-h(\varphi_\ept)) \ov\sigma
  +\int_{Q_t}\Big(P\ov\sigma h(\varphi_\ept) 
  + (P\sigma -A)(h(\varphi_{\ept})-h(\varphi))
  \Big)\mathcal N(\ov\varphi-(\ov\varphi)_\Omega)\\
  &\leq\luca{\frac{C_0}{4}}\int_{Q_t}|\ov\varphi|^2
  +M\left(\int_{Q_t}|\ov\sigma|^2 + 
  \int_0^t\norm{(\ov\varphi-(\ov\varphi)_\Omega)(s)}^2_{V^*}\,\d s\right)
\end{align*}
for a certain $M>0$ independent of $\tau$ and $\eps$,
\luca{and similarly we have the estimate 
\[
  \int_{Q_t}(J*\ov\varphi)\ov\varphi\leq (a^*+b^*)
  \int_0^t\norm{\ov\varphi(s)}_H\norm{\ov\varphi(s)}_{V^*}\,\d s
  \leq\frac{C_0}8\int_{Q_t}|\ov\varphi|^2 + M\int_0^t\norm{\ov\varphi(s)}^2_{V^*}\,\d s\,.
\]}
Finally, as for the fifth term on the right-hand side of \eqref{aux3_ept} we have,
for a positive $\delta$ yet to be chosen,
\begin{align*}
  \int_{Q_t}\ov\mu(\ov\varphi)_\Omega=
  |\Omega|\int_0^t(\ov\varphi(s))_\Omega(\ov\mu(s))_\Omega\,\d s\leq
  \delta\int_0^t|(\ov\mu(s))_\Omega|^2\,\d s + 
  \frac{\luca{|\Omega|^2}}{4\delta}\int_0^t|(\ov\varphi(s))_\Omega|^2\,\d s\,,
\end{align*}
where, by comparison in equation \eqref{cd_epstau2} and by
using the estimate \eqref{est2_ept},
\begin{align*}
  \int_0^t|(\ov\mu(s))_\Omega|^2\,\d s\leq M\left(
  \tau+
  \int_0^t\norm{F'(\varphi_\ept(s))-F'(\varphi(s))}^2_{L^1(\Omega)}\,\d s
  +\chi^2\int_{Q_t}|\ov\sigma|^2\right)\,.
\end{align*}
Next, owing to \eqref{pot_reg} and the H\"older inequality, we infer that 
\begin{align*}
  &\int_0^t\norm{F'(\varphi_\ept(s))-F'(\varphi(s))}^2_{L^1(\Omega)}\,\d s
  \leq C_F^2\int_0^t\norm{(1+|\varphi_\ept|^2+|\varphi|^2)\ov\varphi}_{L^1(\Omega)}^2\,\d s\\
  &\qquad\leq M'\left(1+\luca{\norm{\varphi_\ept}^4_{L^\infty(0,T; L^4(\Omega))}
  +\norm{\varphi}^4_{L^\infty(0,T; L^4(\Omega))}}\right)\int_0^t\norm{\ov\varphi(s)}_H^2\,\d s
\end{align*}
for some $M'>0$ independent of $\eps$ and $\tau$, which yields in turn, 
due to assumption \eqref{pot_reg2} and to the previous estimates,
\[
  \int_0^t\norm{F'(\varphi_\ept(s))-F'(\varphi(s))}^2_{L^1(\Omega)}\,\d s
  \leq M^*\int_{Q_t}|\ov\varphi|^2
\]
for a constant $M^*>0$ independent of $\eps$ and $\tau$.
Thus, collecting the above estimates and rearranging the terms,
we see that choosing $\delta>0$ sufficiently small,
\luca{for example $\delta=\frac{C_0}{4M{M}^*}$,}
we are left with 
\begin{align}
\nonumber
	& \frac 12 \norma{(\ov\varphi-(\ov\varphi)_\Omega)(t)}_{V^*}^2
	+ \luca{\frac{C_0}8\intQt |\ov\ph|^2}
	+ \frac 12 \norma{\ov\sigma(t)}_H^2
	+ \intQt |\nabla\ov\sigma|^2
	\\ 
	&\leq \frac12\norm{\ov\ph_0}_{V^*}^2
	+\frac12\norm{\ov\sigma_0}_H^2
	\label{aux4_ept}
	+M\left(\eps^{1/2} + \tau + 
	\int_0^t\norm{(\ov\varphi-(\ov\varphi)_\Omega)(s)}_{V^*}^2\,\d s
	+\int_{Q_t}|\ov\sigma|^2
	+\int_0^t\norm{\ov\varphi(s)}^2_{V^*}\,\d s\right)\,.
\end{align}

In order to conclude, we only need to handle the last terms on the right-hand side of
\eqref{aux4_ept}. To this end, note that integrating equation \eqref{cd_epstau1} 
on $\Omega$ and testing by $(\ov\varphi)_\Omega$ yields,
using the estimate \eqref{error1_ept}, the Young inequality,
 the boundedness of $\sigma$ and $h$, and the Lipschitz continuity of $h$,
\begin{align}
  \nonumber
  \frac12|(\ov\varphi(t))_\Omega|^2 &= \frac12|(\ov\varphi_0)_\Omega|^2
  -\eps\int_0^t(\ov\varphi(s))_\Omega(\partial_t\mu_\ept(s))_\Omega\,\d s\\
  \nonumber
  &\qquad+\int_0^t\Big(P\ov\sigma(s) h(\varphi_\ept(s)) +
  (P \sigma(s) -A)(h(\varphi_{\ept}(s))-h(\varphi(s)))\Big)_\Omega
  (\ov\varphi(s))_\Omega\,\d s\\
  \label{aux5_ept}
  &\leq\frac12\norm{\ov\varphi_0}_{V^*}^2 + M\left(\eps^{1/2} +
  \int_0^t|(\ov\varphi(s))_\Omega|^2\,\d s + \int_{Q_t}|\ov\sigma|^2\right)+
  \luca{\frac{C_0}{16}}\int_{Q_t}|\ov\varphi|^2\,.
\end{align}

Summing then \eqref{aux4_ept} and \eqref{aux5_ept}, 
we infer that 
\begin{align*}
\nonumber
	& \norma{\ov\varphi(t)}_{V^*}^2
	+ \intQt \!\!|\ov\ph|^2
	+ \norma{\ov\sigma(t)}_H^2
	+ \intQt\!\! |\nabla\ov\sigma|^2
	\leq M\left(
	\norm{\ov\ph_0}_{V^*}^2
	+\norm{\ov\sigma_0}_H^2
	+\eps^{1/2} + \tau + 
	\int_0^t\norm{\ov\varphi(s)}_{V^*}^2\,\d s
	+\int_{Q_t}\!\!|\ov\sigma|^2
	\right)
\end{align*}
for a certain constant $M$, independent of $\eps$ and $\tau$. 
Therefore, we invoke the Gronwall lemma to conclude.
\luca{It is not difficult to check that the same argument performed here yields
uniqueness of solutions for the limit problem, even
without assuming \eqref{init3_ept} and \eqref{init4_ept}.}
This concludes the proof of Theorem~\ref{thm11}.

\luca{
\section*{Conclusions}
Large part of the applied literature on tumor growth modeling agrees
that relevant biological mechanisms such as cell-to-cell adhesion
are typically a non-local process. In this spirit, 
in this paper we introduce and investigate 
from a mathematical perspective a wide class non-local models of
tumor growth capturing long-range interactions in cell-invasion.
The analyzed model contains two regularization coefficients
$\eps$ and $\tau$, which allow the investigation in very broad scenarios
such as the thermodynamically--relevant potentials 
and crucial mechanisms of chemotaxis and active transport.
Then, we perform a complete asymptotic analysis 
showing how the parameters $\eps$ and $\tau$ may
approach zero, both separately and jointly, allowing us to establish well-posedness 
of the limiting systems obtained by formally setting to zero those coefficients.}

\section*{Acknowledgement}
L.S.~gratefully acknowledges financial support from
the Austrian Science Fund (FWF) project M 2876-N. 
\luca{The authors are very grateful to the anonymous referees 
and to the Editors for the constructive criticisms and insightful comments 
that led to a more complete and thorough presentation of the results.}


\bibliographystyle{abbrv}

\def\cprime{$'$}

\end{document}